\pdfoutput=1

\documentclass[11pt,a4paper]{article}

\usepackage[T1]{fontenc}
\usepackage[utf8]{inputenc}
\usepackage{amsmath,amssymb,amsthm}
\usepackage{tikz-cd}
\usepackage{booktabs,array}
\usepackage{url}
\usepackage[margin=1in]{geometry}
\usepackage[hidelinks]{hyperref}

\theoremstyle{definition}
\newtheorem{theorem}{Theorem}[subsection]
\newtheorem{proposition}[theorem]{Proposition}
\newtheorem{lemma}[theorem]{Lemma}
\newtheorem{definition}[theorem]{Definition}
\newtheorem{corollary}[theorem]{Corollary}
\newtheorem{example}[theorem]{Example}
\newtheorem{problem}[theorem]{Problem}
\newtheorem{remark}[theorem]{Remark}

\title{Decomposition of completeness and the universality of the Daniell extension}
\author{Hayato Saito}
\date{\today}

\begin{document}

\maketitle

\begin{abstract}
We treat Daniell-style integration systems, introduce the category
$\mathbf{Int}$ whose morphisms are pullbacks along maps of the base spaces,
and locate the Daniell extension as a structure of this category.

The main result is that the Daniell extension defines a functor
$\mathcal{D}:\mathbf{Int}\to\mathbf{DInt}$ and is left adjoint to the inclusion
of the full subcategory $\mathbf{DInt}$ of Daniell closed integration systems.
That is, $\mathbf{DInt}$ is a reflective subcategory of $\mathbf{Int}$ and the
Daniell extension is its reflector.

The paper contributes the following four points.

\begin{enumerate}
  \item \textbf{Decomposing Daniell closedness into two closure conditions.}
    An integration system is Daniell closed (that is, $L^1=L$) if and only if it
    satisfies the monotone convergence theorem and is closed under modification
    on null sets. Neither condition can be dropped. The two conditions on the
    right occur independently in measure theory; the second one corresponds to
    completeness of the measure.
  \item \textbf{Determining the second reflection.}
    Requiring Daniell closedness of the codomain is not a dispensable
    restriction: the Daniell extension itself is \emph{not} left adjoint to the
    inclusion of the category of integration systems satisfying the monotone
    convergence theorem. That inclusion does have a left adjoint, but its
    reflector is \emph{not} the Daniell extension. Under the Stone condition
    this second reflection is determined as the set of all elements of $L^1$
    that are measurable with respect to the $\sigma$-algebra generated by $L$.
    This characterizes, in terms of $(X,L,I)$ itself, when the two reflections
    agree. For the integration system of step functions it consists of the
    Borel measurable integrable functions. Moreover, for an arbitrary
    integration system, $L^1$ splits as the sum of its $\sigma(L)$-measurable
    part and the functions supported on null sets. Hence \textbf{the gap between
    the two reflections amounts exactly to the functions supported on null
    sets}, which in the case of step functions appears as the difference between
    Borel and Lebesgue, that is, as the completion of the measure. Neither
    separation of points by $L$ nor the Stone condition is sufficient for the
    two reflections to agree.
  \item \textbf{Assembling the whole into a category.}
    We define identities and composition, verify that the Daniell extension
    satisfies the functor laws, and make the uniqueness of the factorization
    explicit. This part is formal. Analytic content enters the uniqueness only
    in the category obtained by dividing morphisms by null functions and in the
    category of operator morphisms below.
  \item \textbf{Making explicit what the category of operator morphisms sees.}
    We detach morphisms from maps of the base spaces and redefine them as
    positive linear maps between Riesz spaces. Operator morphisms extend
    uniquely to a codomain whose quotient is a Banach space, and the Daniell
    extension again gives a reflection in this category. However, morphisms of
    this category see only the quotient $\widetilde{L}=L/\mathcal{N}_I$: the
    assignment $(X,L,I)\mapsto(\widetilde{L},\widetilde{I})$ is faithful and
    full, so neither the base space nor the realization as a function space is
    distinguished. What the Daniell extension gives in this category is the
    Banach completion. The subcategory to be taken as the codomain of the main
    theorem is the class of objects whose quotient is a Banach space, not the
    class of Daniell closed objects; the latter is not isomorphism-closed, and
    the former is its isomorphism closure. Consequently the distinction between
    an abstract completion and a concrete realization is retained as categorical
    structure only on the side where morphisms are maps of the base spaces.
\end{enumerate}

An appendix treats the general theory of adjoint functors used in the proofs,
and the idea of forgetting the base space.

The construction of the Daniell extension by series representations, the fact
that the extension is again an integration system, that it satisfies the
monotone convergence theorem, that it is idempotent, the minimality of the
extension for a fixed base space, and the notion corresponding to a morphism of
integration systems are all known and are not claimed as new here. Among the
four points above, (1) and (2) are the parts with analytic content, but both lie
in the neighbourhood of classical Daniell theory, and the claim is only that,
\textbf{to the author's knowledge, they have not been stated explicitly}.
Point (3) is formal and point (4) is an observation about the choice of category.

Taken together, (1) and (2) split the Daniell extension into a stage that takes
the $\sigma(L)$-measurable integrable functions and a stage that completes the
measure. Both are classical operations. The position of this paper is therefore
not a new completion theorem, but rather
\textbf{the Daniell counterpart of the completion of a measure space, organized
in categorical language}.
Section~\ref{chap:intro} discusses the relation to earlier work in detail.
\end{abstract}

\tableofcontents

\section{Introduction}\label{chap:intro}

\subsection{Background}

Daniell-style integration theory \cite{Daniell1918} does not construct a measure
first, but starts from a positive continuous linear functional on a Riesz space
of functions.
From this standpoint an integral is strictly an operation assigning real numbers
to functions on a base space, and it requires no additional structure, such as a
$\sigma$-algebra, on the base space.
The advantage of this standpoint is that it separates out how much of what is
needed to construct the Lebesgue integral belongs to the side of the base space
and where the side of the function space and the functional begins.

This advantage is, however, usually discussed as an argument internal to a
single integration system.
How integration systems on different base spaces are to be related to one
another, and how the Daniell extension behaves with respect to that relation,
are not normally made an explicit subject.
The purpose of this paper is to organize that part in the language of
categories.

Note that requiring no structure on the base space was also the historical
motivation for this standpoint.
Daniell himself applied this construction to infinite-dimensional base spaces
\cite{Daniell1919}.

\subsection{Organization and contributions of this paper}

In this paper the examples are not collected into an independent section but
placed immediately after the definition each of them tests.
Examples of integration systems are
Example~\ref{ex:counting} and Example~\ref{ex:step-functions}; that the Daniell
extension recovers the Lebesgue integral is
Example~\ref{ex:lebesgue}; that operator morphisms include stochastic matrices
is Example~\ref{ex:stochastic-matrix}.
In Subsection~\ref{sec:measure} we also confirm, via the Stone condition and the
Daniell--Stone representation theorem, that integration systems recover measure
theory.
The results of this paper do not use that condition.
That is, this paper \textbf{does not assume the Stone condition and does not
take representation by a measure as a premise of the theory}.

Section~\ref{chap:preliminaries} is preliminary.
It defines integration systems $(X,L,I)$, integrability via series
representations, the Daniell extension $(X,L^1,I^1)$, and the monotone
convergence theorem and Daniell closedness, and lists the statements needed
later.
\textbf{Nothing in this section is new, and as a rule we refer to earlier work
for the proofs as well.}

Section~\ref{chap:functoriality} is the main subject of this paper.
We introduce the category $\mathbf{Int}$ of integration systems and their
morphisms, the full subcategory $\mathbf{MInt}$ of integration systems
satisfying the monotone convergence theorem, and the further full subcategory
$\mathbf{DInt}$ of Daniell closed ones, and show that the Daniell extension
defines a functor
\[
  \mathcal{D}:\mathbf{Int}\to\mathbf{DInt}
\]
and that this is left adjoint to the inclusion functor $\iota_{D}$.
That is, $\mathbf{DInt}$ is a reflective subcategory of $\mathbf{Int}$ and the
Daniell extension is its reflector.

We then show, by a counterexample on a finite set, that the restriction imposing
Daniell closedness on the codomain in this theorem cannot be dispensed with.
That is, the functor $\iota_M\circ\mathcal{D}:\mathbf{Int}\to\mathbf{MInt}$
defined by the Daniell extension is \emph{not} left adjoint to the inclusion
functor $\iota:\mathbf{MInt}\to\mathbf{Int}$.
We then show that $\iota$ does have a left adjoint, but that its reflector is
\emph{not} the Daniell extension
(Theorem~\ref{thm:mint-inclusion-has-left-adjoint}).
We investigate when these two reflections agree.
Using the $\sigma$-algebra $\sigma(L)$ generated by $L$, the second reflection
$L^{\mathrm{M}}$ is contained in the $\sigma(L)$-measurable part of $L^1$
(Proposition~\ref{prop:LM-in-Lsigma}), and under the Stone condition the two
coincide (Theorem~\ref{thm:LM-equals-Lsigma}).
This characterizes the agreement of the two reflections as a property of
$(X,L,I)$ itself (Corollary~\ref{cor:when-coincide-sigma}).
As a corollary, for the integration system consisting of the step functions on
$[0,1)$ and the Riemann integral, $L^{\mathrm{M}}$ is the set of all Borel
measurable integrable functions and $L^1$ is the set of all Lebesgue integrable
functions, so the two reflections do not agree
(Corollary~\ref{cor:lebesgue-reflections-differ}).
Furthermore, for an arbitrary integration system, $L^1$ splits as the sum of its
$\sigma(L)$-measurable part and the part supported on null sets
(Proposition~\ref{prop:l1-decomposition}).
Hence \textbf{the gap between the two reflections amounts exactly to the
functions supported on null sets}, and in the case of step functions this
appears as the difference between Borel and Lebesgue, that is, as the completion
of the measure.
Moreover, \textbf{neither separation of points by $L$ nor the Stone condition is
a sufficient condition for the two reflections to agree}.
The remaining question is whether the Stone condition is necessary
(Problem~\ref{prob:when-reflections-coincide}).
Finally we prove a main theorem of the same shape in the category obtained by
dividing morphisms by null functions (Theorem~\ref{thm:main-ae}).

In Section~\ref{chap:operator} we detach morphisms from maps of base spaces and
redefine them as positive linear maps between Riesz spaces.
This is based on the operators introduced in \cite{Yamagami} in the context of
iterated integrals, but with the order continuity imposed there removed.
The reason for removing it is given at the beginning of
Section~\ref{chap:operator} and in
Remark~\ref{rem:order-continuity-redundant}.
We show that in this category an operator morphism extends uniquely to a
codomain whose quotient is a Banach space
(Theorem~\ref{thm:operator-extension}), and from this we obtain the
operator-morphism version of the main theorem
(Theorem~\ref{thm:operator-main-theorem}).
\textbf{Here uniqueness does not follow automatically from the agreement of maps
of base spaces, but is obtained by way of contractivity and density.}
This is a difference from the non-ae version in
Section~\ref{chap:functoriality}.
The ae version of that section
(Subsection~\ref{sec:main-theorem-ae}) uses the same density.
We also see that, since the extension is determined only up to null functions,
it becomes natural to divide the codomain by null functions (we do not, however,
prove that ``the extension is undetermined unless one divides'';
Problem~\ref{prob:no-lift}).

Finally we investigate the relation between the two sections.
The functor from the category whose morphisms are maps of base spaces to the
category of operator morphisms becomes faithful once the ae-quotient is taken,
but is not full.
Using this, one finds that the \emph{uniqueness} part of the main theorem of
Section~\ref{chap:functoriality} follows from the theorem of
Section~\ref{chap:operator}, while \emph{existence} does not.
We also show that the morphisms of the category of operator morphisms see only
the quotient $\widetilde{L}$
(Proposition~\ref{prop:Q-full-faithful}).
Hence that category distinguishes neither the base space nor the realization as
a function space, and what the Daniell extension gives is the Banach completion
of a pre-AL space
(Proposition~\ref{prop:D-op-is-banach-completion}).
The subcategory to be taken as the codomain of the main theorem in this category
is the class of all objects whose quotient is a Banach space, not the class of
all Daniell closed objects.
The latter is not isomorphism-closed, and the former is its isomorphism closure
(Proposition~\ref{prop:BInt-characterization}).
\textbf{That is, only the side of Section~\ref{chap:functoriality} retains the
distinction between an abstract completion and a concrete realization as
categorical structure.}
This is the reason for keeping the two sections apart.

\medskip
\noindent
\textbf{The categories and functors of this paper.}
On the side where morphisms are maps of the base spaces there are two
reflections.
\[
  \begin{tikzcd}[column sep=4.5em, row sep=3.6em]
    \mathbf{DInt}
      \arrow[r, hook, "\iota_M"]
    & \mathbf{MInt}
      \arrow[r, hook, "\iota"]
    & \mathbf{Int}
      \arrow[l, bend right=42, "(-)^{\mathrm{M}}"']
      \arrow[ll, bend left=26, "\mathcal{D}"]
  \end{tikzcd}
\]
Both curved arrows are left adjoints:
$\mathcal{D}\dashv\iota_{D}$ (where $\iota_{D}=\iota\circ\iota_M$) and
$(-)^{\mathrm{M}}\dashv\iota$, the latter being $A\mapsto A^{\mathrm{M}}$
(Theorem~\ref{thm:mint-inclusion-has-left-adjoint}).
\textbf{These two reflections do not agree}, and that is the starting point of
Subsection~\ref{sec:why-complete}
(Remark~\ref{rem:D-I-not-left-adjoint}).

Dividing morphisms by null functions and then replacing them by operator
morphisms gives the following diagram.
\[
  \begin{tikzcd}[column sep=5em, row sep=3.4em]
    \mathbf{Int}_{\mathrm{ae}}
      \arrow[r, "\mathcal{D}_{\mathrm{ae}}"]
      \arrow[d, "\mathcal{E}_{\mathrm{ae}}"']
    & \mathbf{DInt}_{\mathrm{ae}}
      \arrow[d, "\mathcal{E}_{\mathrm{ae}}"]
    \\
    \mathbf{Int}_{\mathrm{op},\mathrm{ae}}
      \arrow[r, "\mathcal{D}_{\mathrm{op}}"']
      \arrow[d, "Q"']
    & \mathbf{BInt}_{\mathrm{op},\mathrm{ae}}
    \\
    \mathbf{AL}_0
  \end{tikzcd}
\]
The category $\mathbf{Int}_{\mathrm{ae}}$ in the upper left is obtained by
dividing the morphisms of $\mathbf{Int}$ by ae-equivalence
(Definition~\ref{def:int-ae}).
Both horizontal arrows are reflections
(Theorem~\ref{thm:main-ae}, Theorem~\ref{thm:operator-main-theorem}), and the
square commutes (Proposition~\ref{prop:E-commutes-with-D}).
The last functor $Q$ is faithful and full
(Proposition~\ref{prop:Q-full-faithful}), which is what the statement above
that the bottom row sees only the quotient $\widetilde{L}$ amounts to.

\medskip
The appendix treats the general theory of adjoint functors used in the proofs
and the idea of forgetting the base space.

\subsection{Relation to earlier work}\label{sec:prior-work}

\subsubsection{On the content of the preliminary section}

The content of Section~\ref{chap:preliminaries} follows \cite{Yamagami}.
The definition of integrability by series representations, the fact that $L^1$ is
a Riesz space, that $I^1$ is an integral, that the Daniell extension satisfies
the monotone convergence theorem, and the idempotence of the Daniell extension
are all in that book.
As that book itself states, these procedures follow \cite{BM}, and \cite{BM} in
turn points out that the construction of the complete Daniell space by this
method was already presented in \cite{M1989}.
The origin of the method itself is \cite{Mikusinski}.

That is, the lineage of the content stated in
Section~\ref{chap:preliminaries} is
\[
  \text{\cite{JM1974}},\ \text{\cite{Mikusinski}}
  \ \to\ \text{\cite{M1989}}
  \ \to\ \text{\cite{BM}}
  \ \to\ \text{\cite{Yamagami}}
  \ \to\ \text{this paper} ,
\]
and \textbf{this paper claims no novelty whatsoever for this section.}
Note that \cite{JM1974} has the same title as \cite{M1989} and is fifteen years
earlier, and may be the oldest source for the construction by series, but the
relation between the two could not be confirmed here.

Note also that the series representation adopted in this paper is the form of
\cite{Yamagami}, in which the dominating functions $\varphi_n$ are made
explicit, and differs on its face from the form without dominating functions
adopted in \cite{BM}.
That the two are equivalent is shown in
Lemma~\ref{lem:series-rep-equiv}.

\subsubsection{On arguments with the base space fixed}

The behaviour of the Daniell extension with the base space $X$ fixed, namely
that the extension is idempotent, that it is monotone with respect to inclusion
of extensions, and that it is the ``smallest complete extension'' (in the
terminology of this paper the smallest Daniell closed extension; see
Remark~\ref{rem:why-daniell-closed}), is stated as a corollary in
\cite{Yamagami}.
\cite{BM} makes remarks to the same effect.
There is also, under the name of the ``closure'' of a Daniell space,
\cite{CWS} and \cite{CFW}, which develop this kind of construction
systematically.

\textbf{Hence this paper claims nothing about minimality with the base space
fixed.}
What this paper treats is the part concerning what can be said when the base
space varies.

\subsubsection{On morphisms and categories}

A notion corresponding to a morphism between integration systems in fact also
appears in \cite{Yamagami}.
In the context of discussing iterated integrals, that book introduces linear
maps $\Lambda$ from a Riesz space $L$ on a set $S$ to a Riesz space $M$ on a set
$T$ that are positive and order continuous, and points out that these include
probability kernels, Markov operators and transition probabilities.
It also contains the fact that $\Lambda$ extends to $L^1$ for integrable
functions (although determined only up to null sets).

\textbf{What remains for this paper is the part that assembles these into a
category, and the part that decomposes Daniell closedness into two closure
conditions.}
Specifically, the following four points.

\begin{enumerate}
  \item Defining identities and composition, and verifying that the Daniell
    extension satisfies the functor laws.
  \item \textbf{Uniqueness of the factorization}. What earlier work supplies is
    the \emph{existence} of the extension, not that the morphism factoring
    through a given morphism is unique.
    Uniqueness is an assertion that acquires meaning only once a class of
    morphisms is fixed, and cannot even be posed as a question without
    introducing a category.

    Whether analytic content enters the uniqueness, however, depends on the
    choice of morphisms.
    When the maps of base spaces are taken as morphisms as they stand,
    uniqueness follows formally from the agreement of the maps of base spaces
    and has no content.
    When morphisms are divided by null functions
    (Subsection~\ref{sec:main-theorem-ae}) and when operator morphisms are
    adopted (Section~\ref{chap:operator}), an argument by way of the density of
    $L$ in $L^1$ becomes necessary.
  \item \textbf{Decomposing Daniell closedness into two closure conditions}
    (Theorem~\ref{thm:complete-intrinsic}).
    As long as $\mathbf{DInt}$ is defined by the condition $L^1=L$, the main
    theorem says no more than that the subcategory of fixed points of
    $\mathcal{D}$ has $\mathcal{D}$ as its reflector, and carries no new
    content.
    Through the equivalence of Daniell closedness with
    \begin{quote}
      satisfying the monotone convergence theorem and being closed under
      modification on null sets,
    \end{quote}
    the main theorem takes the form ``the reflector into the subcategory cut out
    by two conditions that occur independently in measure theory is given by the
    concrete construction of series representations''
    (Remark~\ref{rem:main-theorem-not-tautology}).
    This is not a matter of assembling a category but a part with analytic
    content.
    The line of argument ``it has content because it can be written without
    mentioning $L^1$'' does not, however, work.
    The reason is given in
    Remark~\ref{rem:not-about-mentioning-L1}.
  \item On the basis of the above, stating that $\mathbf{DInt}$ is a reflective
    subcategory.
\end{enumerate}

The four results treated in
Subsection~\ref{sec:why-complete}, namely the counterexample showing that the
Daniell extension is not left adjoint to
$\iota:\mathbf{MInt}\to\mathbf{Int}$, the construction of the reflection into
$\mathbf{MInt}$, the theorem determining that reflection by the
$\sigma$-algebra generated by $L$
(Theorem~\ref{thm:LM-equals-Lsigma}), and the decomposition of $L^1$
(Proposition~\ref{prop:l1-decomposition}), were obtained in this paper
independently of the above.
Taking the last two together, it follows that the gap between the two
reflections amounts exactly to the functions supported on null sets, and that for
the integration system of step functions this is the difference between Borel
and Lebesgue.
These, together with
Theorem~\ref{thm:complete-intrinsic}, are the parts of this paper with technical
content.
Both lie, however, in the neighbourhood of classical Daniell theory, and the
claim is only that \textbf{to the author's knowledge they have not been stated
explicitly}.

\subsubsection{Neighbouring work}

There is a body of work attempting to treat integration theory categorically.
\cite{Hartig} restates the Riesz representation theorem on compact spaces as a
natural equivalence between the functor giving the Banach dual of the space of
continuous functions and the functor giving the space of signed Baire measures.

In particular \cite{BorgerHabil}, and \cite{Borger2010}, which treats its
subject in a general category, construct the integral itself by a universal
property, starting from Boolean algebras.
For a Boolean algebra $\mathcal{B}$, the assignment of indicator functions is the
universal charge into the space of simple functions, and any charge factors
through it uniquely; the way it factors is what corresponds to ``the integral of
simple functions''.
The functor taking the algebra of simple functions is left adjoint to the functor
assigning to a ring the Boolean algebra of its idempotents.

As far as we have read the text of the published version \cite{Borger2010}, this
is a different theorem from the present one.
The starting point differs (a Boolean algebra, versus a pair consisting of a
Riesz space and an integral), the content of the universal property differs
(constructing an integral, versus extending an integration system already
given), and the destination differs as well (the space of simple functions,
versus the space $L^1$ of integrable functions).
The universal property in this paper does not produce an integral; it is a
reflection sending a given integration system to a Daniell closed one.

Moreover, the boundary between the two is made explicit within that paper.
What $LB$ represents is the space of simple functions, and immediately after its
Theorem 2.1 that paper states that, given suitable completeness conditions, the
integral can be extended from simple functions to more general measurable
functions, without carrying out that extension.
Concerning the monotone convergence theorem, it goes no further than remarking in
the introduction that convergence theorems such as monotone convergence and
dominated convergence ought to mean continuity in some topology.
\textbf{What this paper treats is precisely that side, the extension.}
Neither the Daniell integral, nor Riesz spaces, nor $L^1$ appears in that paper,
and there is no assertion about reflective subcategories or completions.
Its MSC classification, too, does not include 28C05, the subject of this paper.

We have not, however, obtained the text of \cite{BorgerHabil} (219 pages).
\cite{Borger2010} states that a preprint gives, on the basis of Linton's method,
a treatment in a special category similar to that of the
Habilitationsschrift, from which one can read off that that work treats the same
subject as above.
Even so, we record that \textbf{the possibility that results of the same kind are
contained there cannot be entirely ruled out}.

Work treating measure theory and probability from the categorical side has also
been active in recent years.
\cite{FritzLorenzin} regards Boolean $\sigma$-algebras and monotone
$\sigma$-complete C*-algebras as abstract measurable spaces and organizes
functors, adjunctions and equivalences with the category of measurable spaces.
What that paper treats, however, are measurable spaces, not objects carrying a
measure or an integral.
\cite{DahlqvistKozen} constructs a symmetric monoidal closed category of ordered
Banach spaces and regular operators and uses it for the semantics of
probabilistic programs.
\cite{EhrhardGeoffroy}, motivated by the semantics of probabilistic programming
languages, develops a theory of integration on cones modelled on the Pettis
integral and shows that the category whose morphisms are the integral-preserving
linear maps is a model of linear logic.
Its objects and morphisms differ from those of this paper, but it makes the same
choice in taking preservation of the integral as the condition on morphisms.
This reference was pointed out to the author by T.~Fritz; the author has read
only its abstract and has made no detailed comparison with its contents.

Also, the normed space obtained in this paper by dividing an object of
$\mathbf{Int}$ by the null functions is the pre-completion form of a space whose
norm is additive on the positive cone, that is, of an abstract Lebesgue space
(AL-space).
The operation of building this completion from a Riesz space and a positive
linear functional is itself standard and can be found in, for instance,
\cite{AbramovichAliprantis}.

Indeed, assertions of the form ``the completion is a reflector'' are a
well-established type in the theory of ordered vector spaces.
For example \cite{BMO} treats the Dedekind completion functorially, and
\cite{AvilesBilokopytov} states it explicitly as a reflector.
\textbf{What these treat, however, is the Dedekind completion, a completion
determined by the order structure alone, and it is a different completion from
the Daniell extension of this paper, which depends on the integral $I$.}
For the integration system of step functions, the Dedekind completion widens the
space in the direction of the order, whereas what the Daniell extension gives is
the space of Lebesgue integrable functions.

On the measure-theoretic side, all the ingredients corresponding to the
reflection of this paper are classically in place.
\cite{FremlinMT2} constructs the completion of a measure space in item 212C and
states in 212D that it is complete, that it extends the original measure, and
that it coincides with the original space if and only if that measure was
already complete.
That this extension is unique on the completed domain is 212E(c), and that the
integrals agree is 212F(b).
As for the class of morphisms, the definition of an inverse-measure-preserving
function is 234A, that the completion preserves this class is 234B(a), that the
class is closed under composition is 234B(g), and the transformation of
integrals along such a map is 235G.

\textbf{That book does not, however, assemble these into a functor, an
adjunction or a reflection.}
Indeed, none of the words ``functor'', ``adjoint'', ``reflect'' or ``category''
occurs anywhere in the relevant Chapters 21 and 23.
That is, the position taken in
Subsection~\ref{sec:our-position} --- \textbf{the ingredients are classical and
the assembly is this paper's contribution} --- is confirmed from the
measure-theoretic side as well.

Note that whether a completion is a reflector depends strongly on the choice of
the class of morphisms.
\cite{BMO} points out that the usual Dedekind completion is not functorial as it
stands and adds structure, and the reflection theorem of
\cite{AvilesBilokopytov} is likewise stated in a category of morphisms
preserving existing suprema.
It is for the same reason that this paper prepares two choices of morphisms and
discusses their relation in
Subsection~\ref{sec:two-chapters}.
What the main theorem of this paper states is not that \emph{an abstract
completion exists}, but that \textbf{the concrete construction of the Daniell
extension realizes that completion as a function space on the same base set, and
does so including the extension of morphisms}.
This distinction is meaningful in the case where Daniell closedness decomposes
into two conditions that occur independently in measure theory, the monotone
convergence theorem and modification on null sets
(Theorem~\ref{thm:complete-intrinsic}).
Under that decomposition the main theorem takes the form ``the reflector into the
subcategory cut out by these two conditions is given by the concrete
construction of series representations''
(Remark~\ref{rem:main-theorem-not-tautology}).
The formulation of Section~\ref{chap:functoriality} is given so that this point
is not lost.

\subsubsection{The position of this paper}\label{sec:our-position}

On the basis of the above, let us state the position of this paper honestly.

The results of Subsection~\ref{sec:why-complete} split the Daniell extension
into the following two stages
(Remark~\ref{rem:D-factors}).
\[
  L
  \quad\subset\quad
  L^{\mathrm{M}}=L^{\sigma}
  \quad\subset\quad
  L^1=L^{\mathrm{M}}+\mathcal{Z}_I
\]
The first stage is the operation of taking the integrable functions measurable
with respect to the $\sigma$-algebra generated by $L$, and the second is the
operation of completing the measure.
\textbf{Both are classical operations, and moreover it is this paper's own
theorems that show them to be so.}

Hence this paper is not ``a new completion theorem distinct from an abstract
completion''. In a phrase, it is

\begin{quote}
  \textbf{the Daniell counterpart of the completion of a measure space,
  organized in the language of categories}.
\end{quote}

\noindent
The four contributions listed up to
Subsection~\ref{sec:prior-work} are the breakdown of that organization.
The decomposition of Daniell closedness
(Theorem~\ref{thm:complete-intrinsic}) supplies the division into the first and
second stages; the determination of the second reflection
(Theorem~\ref{thm:LM-equals-Lsigma},
Proposition~\ref{prop:l1-decomposition}) identifies the destination of each
stage; the assembly into a category restates the composite as a reflector; and
the observations on the operator category show how far one can abstract before
the realization becomes invisible.

It is more accurate to put this position in the foreground than to hide it, and
it also makes it easier for the reader to see what is classical and what is not.

\subsubsection{On the scope of the literature survey}

In writing this paper we searched, via zbMATH Open, the intersection of
MSC 28C05 (integration theory via linear functionals) with 18-XX (category
theory).
Within that range we found no work identifying the Daniell extension as a
reflection into the Daniell closed integration systems.
The works obtained were \cite{Hartig} and \cite{BorgerHabil}, listed above, and
work on section spaces of Banach lattices, whose subject differs from this
paper's.

We also searched for the lineage described above, and confirmed from the review
of \cite{M1989} that that paper is the presentation of the construction by
series itself and contains no minimality, universality or categorical content.
We further searched for work treating abstract Lebesgue spaces and completions of
Riesz spaces as adjunctions or reflections; everything obtained concerned the
Dedekind completion, and we found nothing taking the Daniell extension as its
object.

For \cite{CWS} we consulted the book itself. What we checked is
Chapter 1, Section 4 (closed Daniell spaces), and the outcome is recorded in
Remark~\ref{rem:extended-real-collapse}.

\textbf{We have not, however, been able to carry out a more exhaustive survey.}
We state explicitly what has not been checked.

\begin{itemize}
  \item Of \cite{CWS}, \textbf{Chapter 4 (extensions of a Daniell space and the
    integral), Chapter 5 Sections 3 and 4 (measurability and integrability),
    and the historical remarks at the end of the volume have not been
    checked.}
    As far as one can judge from the table of contents and the index,
    Chapter 4 treats families of extensions and their order structure, which
    places it close to the content of
    Subsection~\ref{sec:why-complete}.
    In particular its Theorem 4.1.10 is, according to the index heading, a main
    theorem about the operation of adjoining exceptional functions, and may be
    an assertion of the same kind as
    Proposition~\ref{prop:l1-decomposition}.
    Likewise Chapter 5 Section 4 treats the relation between measurability and
    integrability and is the place to be compared with
    Theorem~\ref{thm:LM-equals-Lsigma}.
  \item \cite{CFW} is held by few libraries in Japan and we have not been able
    to consult its text.
  \item We have been unable to obtain the texts of
    \cite{BorgerHabil}, \cite{M1989} and \cite{JM1974}, and our grasp of their
    contents rests on the zbMATH reviews or is unconfirmed.
\end{itemize}

\noindent
Accordingly, \textbf{for the results of
Subsection~\ref{sec:why-complete}, and in particular for
Proposition~\ref{prop:l1-decomposition} and
Theorem~\ref{thm:LM-equals-Lsigma}, the comparison with earlier work is not
complete.}
If all or part of the results of this paper are already known, the author would
be grateful to be informed.

\section{Preliminaries: integration systems and the Daniell extension}\label{chap:preliminaries}

Section~\ref{chap:preliminaries} is a preliminary section whose purpose is to fix
the definitions, notation and statements used below.
Its content follows \cite{Yamagami}, and within it the construction of the
Daniell extension follows the procedure of \cite{BM}, as that book itself states.
\textbf{Nothing in Section~\ref{chap:preliminaries} is new.}

To keep the length down, Section~\ref{chap:preliminaries} states as a rule only
definitions and assertions, and we refer to the literature above for proofs.
There are two exceptions.

\begin{itemize}
  \item Lemma~\ref{lem:series-rep-equiv} (equivalence of the definition with
    dominating functions and the one without).
    The form of series representation that this paper adopts, following
    \cite{Yamagami}, differs on its face from the form adopted in \cite{BM}, so
    we give a proof as a bridge between the two.
  \item Lemma~\ref{lem:daniell-ext-is-complete} (Daniell closedness of the
    Daniell extension).
    It is used essentially in the main theorem of
    Section~\ref{chap:functoriality}, so we give a proof.
\end{itemize}

\subsection{Notation}

\begin{definition}
  When a sequence $\{a_n\}\subset\mathbb{R}$ is monotonically nonincreasing and
  $\lim_{n\to\infty}a_n$ converges to a real number $a$, we write
  \[
    a_n\downarrow a .
  \]
  Likewise, when a sequence $\{a_n\}\subset\mathbb{R}$ is monotonically
  nondecreasing and $\lim_{n\to\infty}a_n$ converges to a real number $a$, we write
  \[
    a_n\uparrow a .
  \]
\end{definition}

\begin{definition}
  Let $X$ be a set and $A$ a subset of $X$.
  The function
  \[
    1_A:X \to\mathbb{R};x\mapsto\begin{cases} 1 & x\in A \\ 0 & x\notin A \end{cases}
  \]
  is called the \textbf{indicator function} of the subset $A$\footnote{
    If $A=\emptyset$ then $x\notin\emptyset$ for every $x\in X$, so
    $1_\emptyset$ is the function constantly equal to $0$.
  }.
\end{definition}

\begin{definition}
  Let $X$ be a set.
  We define an order on the family $\mathbb{R}^X$ of functions on $X$ by
  \[
    f\le g \iff \forall x\in X  \; f(x)\le g(x)
  \]
  for all $f,g\in\mathbb{R}^X$.
\end{definition}

\begin{definition}
  We define operations $\land, \lor$ on the reals by
  \begin{align*}
    a\land b&:=\min\{a,b\}\\
    a\lor b&:=\max\{a,b\} ,
  \end{align*}
  and, for a set $X$ and $f,g\in\mathbb{R}^X$, we define the functions
  $f\land g$ and $f\lor g$ by
  \begin{align*}
    (f\land g)(x)&:=f(x)\land g(x) \\
    (f\lor g)(x)&:=f(x)\lor g(x) .
  \end{align*}
\end{definition}

\subsection{Riesz spaces and Riesz systems}

\begin{definition}
  Let $X$ be a set and let $L\subset\mathbb{R}^X$ be an
  $\mathbb{R}$-vector space\footnote{
    We allow $X=\emptyset$. In that case $\mathbb{R}^X=\{\emptyset\}$ is the
    zero-dimensional vector space, and $(\emptyset,\{0\},0)$ is the only
    integration system on the empty set.
    The definitions and assertions of this paper remain meaningful in this case,
    and they can be checked directly.
  }.
  We say that $L$ is a \textbf{Riesz space} on $X$\footnote{
    What this paper calls a Riesz space is called a vector lattice in some books.
    We use ``Riesz space'' throughout.
  } if
  \begin{align*}
    f\land g &\in L \\
    f\lor g &\in L
  \end{align*}
  holds for all $f,g\in L$.
\end{definition}

\begin{remark}[on terminology]\label{rem:riesz-terminology}
  The definition above is for subspaces of $\mathbb{R}^X$, and it requires
  closure under the lattice operations with respect to the pointwise order.
  Such an $L$ is also a Riesz space in the abstract sense (an ordered vector
  space that is a lattice), with the pointwise order as its order.
  That is, a ``Riesz space on $X$'' in this paper is a special case of an
  abstract Riesz space.

  In this paper, \textbf{whenever we attach ``on $X$'' we always mean a
  sublattice of $\mathbb{R}^X$}.
  Riesz spaces in the abstract sense appear only in
  Appendix~\ref{sec:forget}, where we call them \textbf{abstract Riesz spaces}
  to keep them apart.

  This distinction is not a formality, because not every abstract Riesz space can
  be realized as a Riesz space on some $X$.
  Indeed, if $L$ is a sublattice of $\mathbb{R}^X$ then evaluation
  $\mathrm{ev}_x$ at a point gives a Riesz homomorphism from $L$ to $\mathbb{R}$,
  and these separate the elements of $L$.
  But by Proposition~\ref{prop:spec-collapses} there is no nonzero Riesz
  homomorphism $L\to\mathbb{R}$ on $L^1[0,1]$.
  Hence $L^1[0,1]$ cannot be realized as a Riesz space on any set $X$.
\end{remark}

\begin{definition}
  Let $X$ be a set and let $L$ be a Riesz space on $X$.
  The pair
  \[
    (X , L)
  \]
  is called a \textbf{Riesz system} on the set $X$.
  The set $X$ in a Riesz system $(X,L)$ is called the \textbf{base space} of the
  Riesz system $(X,L)$.
\end{definition}

\begin{definition}
  Let $(X,L)$ and $(Y,M)$ be Riesz systems on sets $X$ and $Y$
  respectively.
  A \textbf{morphism of Riesz systems}
  \[
    F:(X,L)\to(Y,M)
  \]
  is a map $\underline{F}:Y\to X$ such that
  \[
    F^*f:=f\circ\underline{F}\in M
  \]
  holds for every $f\in L$.\footnote{Strictly, a morphism is the triple
  consisting of the domain, the codomain and the map of base spaces. Below we fix
  the domain and the codomain and present a morphism by the map of base spaces.}
\end{definition}

\begin{lemma}
  \begin{enumerate}
    \item For every Riesz system $(X,L)$ on a set $X$, the identity map
      $\mathrm{id}_X=\underline{\mathrm{id}_X}:X\to X$ induces a morphism of
      Riesz systems $\mathrm{id}_{(X,L)}:(X,L)\to (X,L)$.
    \item For any two morphisms of Riesz systems $F:(X,L)\to(Y,M)$ and
      $G:(Y,M)\to(Z,N)$, the map $\underline{F}\circ\underline{G}:Z\to X$ induces
      a morphism of Riesz systems $G\circ F:(X,L)\to(Z,N)$.
  \end{enumerate}
\end{lemma}
\begin{proof}
  \begin{enumerate}
    \item Immediate, since $f\circ \operatorname{id}_X=f\in L$ for every $f\in L$.
    \item For every $f\in L$ we have $F^*f\in M$ by hypothesis. Hence, again by
      hypothesis,
      $(G\circ F)^*(f)=f\circ(\underline{F}\circ\underline{G})=(f\circ \underline{F})\circ\underline{G}\in N$,
      which is the assertion.
  \end{enumerate}
\end{proof}

\begin{definition}
  The category whose objects are Riesz systems and whose morphisms are morphisms
  of Riesz systems is written
  \[
    \mathbf{RieszSys}
  \]
  and called the \textbf{category of Riesz systems}.
\end{definition}

\subsection{Integrals and integration systems}

\begin{definition}
  Let $X$ be a set and let $(X,L)$ be a Riesz system on $X$.
  We say that an $\mathbb{R}$-linear map
  \[
    I:L\to\mathbb{R}
  \]
  is an \textbf{integral} on the Riesz system $(X,L)$ if the following two
  conditions hold.
  \begin{itemize}
    \item For every $f\in L$, if $f\ge0$ then $I(f)\ge0$.
    \item For every sequence of functions $\{f_n\}\subset L$, if
      $f_n(x)\downarrow0$ for all $x\in X$ then $\lim_{n\to\infty}I(f_n)=0$.
  \end{itemize}
\end{definition}

\begin{definition}
  Let $(X,L)$ be a Riesz system on a set $X$ and let $I:L\to\mathbb{R}$
  be an integral on $(X,L)$.
  The triple
  \[
    (X,L,I)
  \]
  is called an \textbf{integration system} on the set $X$.
\end{definition}

\begin{example}[counting]\label{ex:counting}
  Let $X=\mathbb{N}$ and let $L$ be the set of all functions that vanish outside
  finitely many $n$.
  Clearly $L$ is a Riesz space on $X$.
  Putting
  \[
    I(f):=\sum_{n=1}^{\infty}f(n) ,
  \]
  the right-hand side is a finite sum and $I$ is an integral on $(X,L)$.
  Indeed, positivity is immediate, and if $f_n\downarrow0$ then the supports of
  the $f_n$ lie in the support of $f_1$, a finite set, so $I(f_n)\to0$.
\end{example}

\begin{example}[step functions and the Riemann integral]\label{ex:step-functions}
  Let $X=[0,1)$ and let $L$ be the set of all finite linear combinations of
  indicator functions of half-open intervals $[a,b)\subset[0,1)$, that is, the
  set of all step functions on $X$.
  Then $L$ is a Riesz space on $X$.
  Taking $I$ to be the Riemann integral
  \[
    I(f):=\int_0^1f(x)\,dx ,
  \]
  the triple $(X,L,I)$ is an integration system.

  Positivity is immediate. Continuity, that is, $I(f_n)\to0$ whenever
  $f_n(x)\downarrow0$ holds for every $x\in X$, follows from additivity of length
  over countable partitions of an interval.
  For a proof see \cite{BM}, \S2.
  The Daniell extension of this example is described in
  Example~\ref{ex:lebesgue}.
\end{example}

\begin{remark}\label{rem:why-half-open}
  There is a reason for taking the base space to be $[0,1)$ rather than $[0,1]$.
  If one takes $X=[0,1]$ and builds $L$ from indicator functions of half-open
  intervals, then every element of $L$ vanishes at $x=1$.
  Then $\sum_n\varphi_n(1)=0<+\infty$ for every series representation, so
  condition (3) of a series representation is imposed at $x=1$ always, and every
  element of $L^1$ satisfies $f(1)=0$.
  This differs from the set of all Lebesgue integrable functions on $[0,1]$.

  This phenomenon is a consequence of the fact that condition (3) of a series
  representation is shaped so as to constrain values only at the points where the
  sum of the dominating functions converges.
  If the base space has a point at which all elements of $L$ are forced to take
  the same value, then the Daniell extension becomes needlessly thin there.
\end{remark}

\begin{definition}
  Let $(X,L,I)$ and $(Y,M,J)$ be integration systems on sets $X$ and
  $Y$ respectively.
  A \textbf{morphism of integration systems}
  \[
    F:(X,L,I)\to(Y,M,J)
  \]
  is a morphism of Riesz systems $F:(X,L)\to(Y,M)$ such that
  \[
    J(F^*f)=I(f)
  \]
  holds for every $f\in L$.
\end{definition}

\begin{lemma}
  \begin{enumerate}
    \item For every integration system $(X,L,I)$ on a set $X$, the
      identity map $\operatorname{id}_X:X\to X$ induces a morphism of
      integration systems
      $\operatorname{id}_{(X,L,I)}:(X,L,I)\to (X,L,I)$.
    \item For any two morphisms of integration systems
      $F:(X,L,I)\to(Y,M,J)$ and $G:(Y,M,J)\to(Z,N,K)$, the map
      $\underline{F}\circ\underline{G}:Z\to X$ induces a morphism of integration
      systems $G\circ F:(X,L,I)\to(Z,N,K)$.
  \end{enumerate}
\end{lemma}
\begin{proof}
  We have already seen that the identity map induces a morphism of Riesz systems
  $\operatorname{id}_{(X,L)}$.
  Moreover $f\circ\underline{\operatorname{id}_X}=f\in L$ and
  $I(f\circ\underline{\operatorname{id}_X})=I(f)$ for every $f\in L$, so the
  identity map clearly induces a morphism of integration systems
  $\operatorname{id}_{(X,L,I)}:(X,L,I)\to (X,L,I)$.

  We have already seen that $\underline{F}\circ\underline{G}:Z\to X$ determines a
  morphism of Riesz systems from $(X,L)$ to $(Z,N)$.
  Moreover, by the hypothesis on morphisms of integration systems,
  \[
    K((G\circ F)^*(f))=K(G^*(F^*(f)))=J(F^*(f))=I(f) .
  \]
  Hence $\underline{F}\circ\underline{G}:Z\to X$ induces a morphism of
  integration systems $G\circ F:(X,L,I)\to(Z,N,K)$.
\end{proof}

\begin{definition}
  The category whose objects are integration systems and whose morphisms are
  morphisms of integration systems is written
  \[
    \mathbf{Int}
  \]
  and called the \textbf{category of integration systems}.
\end{definition}

\begin{remark}[on the direction of morphisms]\label{rem:opposite-convention}
  Note that the map of base spaces underlying a morphism
  $F:(X,L,I)\to(Y,M,J)$ is $\underline{F}:Y\to X$, so that the direction is
  reversed.
  In the language of measure theory, the condition
  $J(f\circ\underline{F})=I(f)$ says that $\underline{F}$ pushes $J$ forward to
  $I$, that is, that $\underline{F}$ is a measure-preserving map
  $(Y,J)\to(X,I)$.

  Hence \textbf{$\mathbf{Int}$ is the opposite of the usual category whose
  morphisms are measure-preserving maps}.
  This paper adopts this direction because pullback of functions is
  contravariant and we wish to treat the Daniell extension as a covariant
  functor.
  A reader who prefers to read everything in the usual convention should replace
  ``reflection'' below by ``coreflection'' and ``left adjoint'' by
  ``right adjoint''.
\end{remark}

\subsection{Integration systems satisfying the monotone convergence theorem}

\begin{definition}\label{def:satisfies-mct}
  We say that an integration system $(X,L,I)$ on a set $X$
  \textbf{satisfies the monotone convergence theorem} if, for every function
  $f:X\to\mathbb{R}$ and every monotonically nondecreasing sequence
  $\{f_n\}_{n=1}^{\infty}\subset L$,
  \[
    f_n(x)\uparrow f(x)\quad(\forall x\in X),\qquad\sup_n I(f_n)<+\infty
  \]
  implies
  \[
    f\in L,\qquad I(f)=\lim_{n\to\infty} I(f_n) .
  \]
\end{definition}

\noindent
This condition is a restatement in axiomatic form of the classical theorem of
Beppo Levi (1906), and \cite{CWS} calls a triple satisfying it \textbf{closed}
(its Definition 1.4.1).
The correspondence of terms with that book is collected in
Remark~\ref{rem:extended-real-collapse}.

\begin{definition}
  The full subcategory of the category $\mathbf{Int}$ of integration systems
  consisting of the integration systems that satisfy the monotone convergence
  theorem is written
  \[
    \mathbf{MInt}
  \]
  and called the \textbf{category of integration systems satisfying the monotone
  convergence theorem}.
\end{definition}

\subsection{The Daniell extension}

The content of this subsection follows \cite{Yamagami}. As that book states, the
facts that $L^1$ is a Riesz space and that $I$ extends to $L^1$ follow the
procedure of \cite{BM}.
\textbf{Nothing in this subsection is new; it is here to fix the definitions,
notation and assertions used in later subsections.}
Accordingly we state only definitions and assertions below, and refer to the
literature above for proofs.

\subsubsection{Definition of the Daniell extension}

\begin{definition}\label{def:series-representation}
  Let $(X,L,I)$ be an integration system on a set $X$.
  Given a function $f:X\to\mathbb{R}$, if there exist sequences of functions
  $\{f_n\}\subset L$ and $\{\varphi_n\}\subset L$ satisfying the following, then
  the pair of sequences $\{f_n\},\{\varphi_n\}$ is called a \textbf{series
  representation} of the function $f$, and a function admitting a series
  representation is said to be \textbf{integrable} with respect to the
  integration system $(X,L,I)$.
  \begin{enumerate}
    \item $|f_n|\le\varphi_n$ for every $n$;
    \item $\sum_{n=1}^\infty I(\varphi_n)<+\infty$;
    \item for every $x\in X$, if $\sum_{n=1}^\infty \varphi_n(x)<+\infty$ then
      $f(x)=\sum_{n=1}^\infty f_n(x)$.
  \end{enumerate}
  We write a series representation $\{f_n\},\{\varphi_n\}$ of a function $f$ as
  \[
    f\overset{\varphi_n}{\simeq}\sum f_n .
  \]
\end{definition}

\begin{definition}
  Let $(X,L,I)$ be an integration system on a set $X$.
  We denote the set of all functions integrable with respect to the integration
  system $(X,L,I)$ by
  \[
    L^1:=\{f\in\mathbb{R}^X \mid\text{$f$ is integrable with respect to the integration system $(X,L,I)$}\}
  \]
  \footnote{
    $L^1$ depends on $I$ as well, but in this paper we never put several
    integrals on one Riesz space $L$.
    We therefore write $L^1$, as something determined by $L$.
    Where this could be ambiguous we write it out, as in $(L^{\mathrm{M}})^1$.
    When referring to the object itself we use the functorial notation
    $\mathcal{D}(X,L,I)=(X,L^1,I^1)$.
  }.
\end{definition}

\medskip
Definition~\ref{def:series-representation} takes the form in which the
dominating functions $\varphi_n$ are made explicit, but this is equivalent to the
form without dominating functions adopted in \cite{BM}.
We will not use this below, but we state it in order to make the correspondence
with the literature clear.

\begin{lemma}[equivalence of the definitions with and without dominating functions]\label{lem:series-rep-equiv}
  Let $(X,L,I)$ be an integration system on a set $X$ and let
  $f:X\to\mathbb{R}$.
  Then the following are equivalent.
  \begin{enumerate}
    \item[(i)] $f$ is integrable in the sense of
      Definition~\ref{def:series-representation}.
      That is, there exist $\{f_n\},\{\varphi_n\}\subset L$ with
      $|f_n|\le\varphi_n$ and $\sum_n I(\varphi_n)<+\infty$ such that
      $f(x)=\sum_n f_n(x)$ holds at every $x$ with
      $\sum_n\varphi_n(x)<+\infty$.
    \item[(ii)] There exists $\{g_m\}\subset L$ with
      $\sum_m I(|g_m|)<+\infty$ such that $f(x)=\sum_m g_m(x)$ holds at every
      $x$ with $\sum_m|g_m(x)|<+\infty$.
  \end{enumerate}
\end{lemma}
\begin{proof}
  For $(ii)\Rightarrow(i)$ take $f_n:=g_n$ and $\varphi_n:=|g_n|$.
  Since $L$ is a Riesz space we have $|g_n|\in L$, and
  $\sum_n\varphi_n(x)<+\infty$ is literally $\sum_n|g_n(x)|<+\infty$.

  We show $(i)\Rightarrow(ii)$. Put $\psi_n:=\varphi_n-|f_n|\ge0$, so that
  $\psi_n\in L$.
  Define the sequence $\{g_m\}$ by
  \[
    f_1,\ \psi_1,\ -\psi_1,\ f_2,\ \psi_2,\ -\psi_2,\ f_3,\ \psi_3,\ -\psi_3,\ \dots .
  \]
  First, $I(|f_n|)\ge0$ because $I$ preserves order, and
  \[
    \sum_m I(|g_m|)
    =\sum_n\bigl(I(|f_n|)+2I(\psi_n)\bigr)
    =\sum_n\bigl(2I(\varphi_n)-I(|f_n|)\bigr)
    \le2\sum_n I(\varphi_n)<+\infty.
  \]
  Next, if $x$ satisfies $\sum_m|g_m(x)|<+\infty$, then
  \[
    \sum_m|g_m(x)|
    =\sum_n\bigl(|f_n(x)|+2\psi_n(x)\bigr)
    \ge\sum_n\bigl(|f_n(x)|+\psi_n(x)\bigr)
    =\sum_n\varphi_n(x)
  \]
  gives $\sum_n\varphi_n(x)<+\infty$. Hence $f(x)=\sum_n f_n(x)$ by condition
  (i).
  On the other hand $\sum_m g_m(x)$ converges absolutely, so rearrangement and
  regrouping are free, and $\psi_n$ cancels against $-\psi_n$, giving
  \[
    \sum_m g_m(x)=\sum_n f_n(x)=f(x) .
  \]
  Therefore (ii) holds.
\end{proof}

\subsubsection{Construction of the extended integral}

\begin{lemma}\label{lem:equiv-condition-for-riesz-space}
  Let $X$ be a set and let $L\subset\mathbb{R}^X$ be an
  $\mathbb{R}$-vector space.
  Then the condition that $|f|\in L$ for every $f\in L$ is equivalent to $L$
  being a Riesz space on $X$.
\end{lemma}

\begin{lemma}\label{lem:series-exp-is-abs-conv}
  Let $(X,L,I)$ be an integration system on a set $X$ and let
  $f\overset{\varphi_n}{\simeq}\sum f_n$.
  Then, for every $x\in X$, if $\sum_{n=1}^\infty \varphi_n(x)<+\infty$ then
  $f(x)=\sum_{n=1}^\infty f_n(x)$ converges absolutely.
\end{lemma}

\begin{theorem}\label{thm:l1-is-riesz-space}
  Let $(X,L,I)$ be an integration system on a set $X$.
  Then the family $L^1$ of functions integrable with respect to the integration
  system $(X,L,I)$ is a Riesz space on $X$.
\end{theorem}

\begin{lemma}\label{lem:integration-preserve-order}
  Let $(X,L,I)$ be an integration system on a set $X$.
  Then $I(f)\le I(g)$ for all $f,g\in L$ with $f\le g$.
\end{lemma}

\begin{lemma}\label{lem:triangle-ineq-for-integration}
  Let $(X,L,I)$ be an integration system on a set $X$.
  Then $|I(f)|\le I(|f|)$ for every $f\in L$.
\end{lemma}

\begin{lemma}\label{lem:abs-conv-for-int-of-series-rep}
  Let $(X,L,I)$ be an integration system on a set $X$.
  For every integrable function $f\in L^1$ with respect to the integration
  system $(X,L,I)$ and every series representation
  $f\overset{\varphi_n}{\simeq}\sum f_n$ of it, the sum
  $\sum_{n=1}^\infty I(f_n)$ converges absolutely.
\end{lemma}

\begin{lemma}[continuity of the integral]\label{lem:integration-is-continuous}
  Let $(X,L,I)$ be an integration system on a set $X$, let $\{f_n\}$ be
  a sequence of functions in $L$, and let $f\in L$.
  Then, if $f_n(x)\uparrow f(x)$ for every $x\in X$, we have
  $I(f_n)\uparrow I(f)$.
\end{lemma}

\begin{lemma}\label{lem:ext-int-of-zero-func}
  Let $(X,L,I)$ be an integration system on a set $X$.
  Then, given a series representation $0\overset{\eta_n}{\simeq}\sum h_n$ of the
  function constantly equal to $0$,
  \[
    \sum_{n=1}^\infty I(h_n)=0 .
  \]
\end{lemma}

\begin{lemma}\label{lem:ext-int-well-defined}
  Let $(X,L,I)$ be an integration system on a set $X$.
  Then, for any two series representations
  \begin{align*}
    f&\overset{\varphi_n}{\simeq}\sum f_n \\
    f&\overset{\psi_n}{\simeq}\sum g_n
  \end{align*}
  of an integrable function $f\in L^1$ with respect to the integration system
  $(X,L,I)$, we have $\sum_{n=1}^\infty I(f_n)=\sum_{n=1}^\infty I(g_n)$.
\end{lemma}

\begin{definition}
  Let $(X,L,I)$ be an integration system on a set $X$.
  For each integrable function $f\in L^1$ with respect to the integration system
  $(X,L,I)$, choose a series representation
  $f\overset{\varphi_n}{\simeq}\sum f_n$.
  Then the absolutely convergent sum
  (Lemma~\ref{lem:abs-conv-for-int-of-series-rep})
  \[
    \sum_{n=1}^\infty I(f_n)
  \]
  does not depend on the choice of series representation, by
  Lemma~\ref{lem:ext-int-well-defined}.
  We define this value to be $I^1(f)$.
\end{definition}

\begin{lemma}\label{lem:ext-int-is-linear}
  Let $(X,L,I)$ be an integration system on a set $X$.
  Then the map $I^1:L^1\to\mathbb{R}$ is an $\mathbb{R}$-linear map on the Riesz
  space $L^1$ on $X$.
\end{lemma}

\begin{lemma}[truncation estimate for a series representation of a positive function]\label{lem:positive-series-representation-truncation}
  Let $(X,L,I)$ be an integration system on a set $X$ and let
  $f\in L^1$ be an integrable function with $f\ge0$.
  Choose a series representation $f\overset{\varphi_n}{\simeq}\sum f_n$ of $f$
  and, for each $N\ge1$, put
  \[
    S_N:=\sum_{k=1}^Nf_k,\qquad h_N:=S_N\vee0,\qquad q_N:=S_N\wedge0 .
  \]
  Then
  \[
    I(|q_N|)\le\sum_{k=N+1}^\infty I(\varphi_k),
    \qquad
    |I^1(f)-I(h_N)|\le2\sum_{k=N+1}^\infty I(\varphi_k)
  \]
  hold. In particular $I(h_N)\longrightarrow I^1(f)$ as $N\to\infty$.
\end{lemma}

\begin{lemma}\label{lem:ext-int-is-positive}
  Let $(X,L,I)$ be an integration system on a set $X$ and let
  $f\in L^1$.
  Then $f\ge0$ implies $I^1(f)\ge0$.
\end{lemma}

\begin{lemma}\label{lem:positive-good-representation}
  Let $(X,L,I)$ be an integration system on a set $X$ and let
  $f\in L^1$.
  Then, if $f\ge0$, for every $\varepsilon>0$ there is a series representation
  $f\overset{\varphi_n}{\simeq}\sum f_n$ of $f$ satisfying
  \[
    \sum_{n=1}^\infty I(\varphi_n)\le I^1(f)+\varepsilon .
  \]
\end{lemma}

\begin{lemma}[flattening series representations into one sequence]\label{lem:flatten-series-representations}
  Let $(X,L,I)$ be an integration system on a set $X$.
  Suppose that for each $n\ge1$ a series representation
  $u_n\overset{\varphi_{n,k}}{\simeq}\sum_{k=1}^\infty g_{n,k}$ of $u_n\in L^1$
  is given, and that
  \[
    \sum_{n=1}^\infty\sum_{k=1}^\infty I(\varphi_{n,k})<+\infty .
  \]
  Put
  \[
    A:=\left\{x\in X\mathrel{}\middle|\mathrel{}\sum_{n=1}^\infty\sum_{k=1}^\infty\varphi_{n,k}(x)<+\infty\right\} .
  \]
  Then $\sum_{n=1}^\infty u_n(x)$ converges absolutely for every $x\in A$.

  Suppose further that a function $u:X\to\mathbb{R}$ satisfies
  $u(x)=\sum_{n=1}^\infty u_n(x)$ for $x\in A$.
  Choose a bijection $\sigma:\mathbb{N}\to\mathbb{N}\times\mathbb{N}$,
  $\sigma(m)=(n_m,k_m)$, and put $\widetilde{g}_m:=g_{n_m,k_m}$ and
  $\widetilde{\varphi}_m:=\varphi_{n_m,k_m}$.
  Then $u\overset{\widetilde{\varphi}_m}{\simeq}\sum_{m=1}^\infty\widetilde{g}_m$
  is a series representation of $u$.
  In particular $u\in L^1$ and
  \[
    I^1(u)=\sum_{n=1}^\infty I^1(u_n)
  \]
  holds, the right-hand side converging absolutely.
\end{lemma}

\begin{theorem}\label{thm:ext-int-is-integral}
  Let $(X,L,I)$ be an integration system on a set $X$.
  Then the map $I^1:L^1\to\mathbb{R}$ is an integral on the Riesz space $L^1$ on
  $X$.
\end{theorem}

\begin{lemma}\label{lem:daniell-ext-contains-source}
  Let $(X,L,I)$ be an integration system on a set $X$ and let $f\in L$.
  Then $f\in L^1$ and $I^1(f)=I(f)$.
\end{lemma}

\begin{definition}
  Let $(X,L,I)$ be an integration system on a set $X$.
  The integration system $(X,L^1,I^1)$ is called the \textbf{Daniell extension}
  of the integration system $(X,L,I)$.
\end{definition}

\begin{example}[the Lebesgue integral]\label{ex:lebesgue}
  Consider the integration system of Example~\ref{ex:step-functions}, that is,
  the integration system $(X,L,I)$ consisting of the step functions on $[0,1)$
  and the Riemann integral.
  Then its Daniell extension $(X,L^1,I^1)$ coincides with the set of all
  Lebesgue integrable functions on $[0,1)$ together with the Lebesgue integral.

  That is, the construction of Section~\ref{chap:preliminaries} recovers the
  Lebesgue integral.
  This is the origin of the definition of integrability by series
  representations; for details see \cite{Mikusinski} and \cite{BM}, \S2.
\end{example}

\begin{corollary}
  Let $(X,L,I)$ be an integration system on a set $X$ and let
  $(X,L^1,I^1)$ be its Daniell extension.
  Then the identity map $\operatorname{id}_X$ induces a morphism of integration
  systems $i_X:(X,L,I)\to(X,L^1,I^1)$.
\end{corollary}

\subsubsection{The Daniell extension satisfies the MCT}

\begin{lemma}[termwise integration of a nonnegative series]\label{lem:positive-series-flattening}
  Let $(X,L,I)$ be an integration system on a set $X$ and let
  $(X,L^1,I^1)$ be its Daniell extension.
  Suppose $\{h_k\}_{k=1}^{\infty}\subset L^1$ with $h_k\ge0$ satisfies
  $\sum_k I^1(h_k)<\infty$.
  Then, putting
  \[
    h(x):=\begin{cases}
      \sum_n h_n(x) & \text{if $\sum_n h_n(x)$ converges} \\
      0 & \text{otherwise,}
    \end{cases}
  \]
  we have $h\in L^1$ and $I^1(h)=\sum_k I^1(h_k)$.
\end{lemma}

\begin{lemma}\label{lem:daniell-ext-satisfies-mct}
  Let $(X,L,I)$ be an integration system on a set $X$.
  Then its Daniell extension $(X,L^1,I^1)$ satisfies the monotone convergence
  theorem (Definition~\ref{def:satisfies-mct}).
\end{lemma}

\subsection{The category of Daniell closed integration systems}

\begin{definition}\label{def:daniell-closed}
  We say that an integration system $(X,L,I)$ on a set $X$ is
  \textbf{Daniell closed} if
  \[
    L^1=L .
  \]
\end{definition}

\begin{remark}[on terminology]\label{rem:why-daniell-closed}
  This condition is often called ``completeness''.
  We assign a different word to it because
  \textbf{``complete'' occurs in this paper in three senses}:
  \begin{itemize}
    \item $L^1=L$ (the present definition);
    \item completeness with respect to the norm $\|\cdot\|_I$
      (Theorem~\ref{thm:l1-norm-complete});
    \item completion of a measure (Subsection~\ref{sec:why-complete},
      Remark~\ref{rem:separation-not-sufficient}).
  \end{itemize}
  These are three distinct conditions, and in this paper all three live in the
  same subsection.
  We therefore call the first one \textbf{Daniell closed}.
  The name follows \cite{CWS} and \cite{CFW}, which use the word ``closure'' of a
  Daniell space.

  \textbf{It also clashes with ``complete'' in other books.}
  \cite{CWS} uses ``complete'' for the sandwich condition ($f\le h\le g$
  together with $\ell(f)=\ell(g)$ implies that $h$ belongs as well), which is
  none of the three above but rather our
  Definition~\ref{def:mod-closed}.
  That book further uses ``complete'' for order-complete lattices and for the
  completion of a measure space.
  The correspondence of terms is collected in
  Remark~\ref{rem:extended-real-collapse}.

  In \textbf{the title}, however, \textbf{we follow custom and say
  ``completeness''}.
  Outside this paper that is the standard word for $L^1=L$, and the coexistence
  of the three senses is an internal matter of this paper.
  In the body and the abstract we use ``Daniell closed'' consistently.
\end{remark}

\begin{remark}\label{rem:complete-implies-mct}
  A Daniell closed integration system automatically satisfies the monotone
  convergence theorem.
  Indeed, the Daniell extension $(X,L^1,I^1)$ satisfies the monotone convergence
  theorem by Lemma~\ref{lem:daniell-ext-satisfies-mct}, and since $L^1=L$ and
  $I^1=I$ (Lemma~\ref{lem:daniell-ext-contains-source}), the system $(X,L,I)$
  itself satisfies the monotone convergence theorem.

  We may therefore define $\mathbf{DInt}$ as a full subcategory of
  $\mathbf{MInt}$.
  Where we list the monotone convergence theorem and Daniell closedness side by
  side below, it is in order to make the relation to $\mathbf{MInt}$ easy to see.
\end{remark}

\begin{definition}
  The full subcategory of the category $\mathbf{MInt}$ of integration systems
  satisfying the monotone convergence theorem consisting of the Daniell closed
  integration systems is written
  \[
    \mathbf{DInt}
  \]
  and called the \textbf{category of Daniell closed integration systems}.
\end{definition}

\begin{definition}\label{def:inclusion-functors}
  We write the inclusion functors among the three categories above as
  \begin{align*}
    \iota&:\mathbf{MInt}\to\mathbf{Int},\\
    \iota_{D}&:\mathbf{DInt}\to\mathbf{Int},\\
    \iota_M&:\mathbf{DInt}\to\mathbf{MInt}
  \end{align*}
  respectively. All of them are inclusions of full subcategories, and
  $\iota_{D}=\iota\circ\iota_M$.
\end{definition}

\begin{lemma}\label{lem:daniell-ext-is-complete}
  Let $(X,L,I)$ be an integration system on a set $X$ and let
  $(X,L^1,I^1)$ be its Daniell extension.
  Then the integration system $(X,L^1,I^1)$ is Daniell closed.
\end{lemma}

\noindent
This statement itself is stated as a corollary in \cite{Yamagami} and is not new.
However, since this paper uses it essentially in the main theorem of
Section~\ref{chap:functoriality}, in the form ``the codomain of the
factorization is Daniell closed'', we give a proof in this subsection only.
\begin{proof}
  We have already checked that the monotone convergence theorem holds, so it
  suffices to show $(L^1)^1=L^1$.
  Since $L^1\subset(L^1)^1$ has already been shown
  (Lemma~\ref{lem:daniell-ext-contains-source}), it suffices to show
  $(L^1)^1\subset L^1$.

  Let $f\in (L^1)^1$. That is, there are sequences of functions
  $\{\varphi_n\},\{f_n\}\subset L^1$ such that
  \begin{enumerate}
    \item $|f_n|\le\varphi_n$;
    \item $\sum_{n=1}^\infty I^1(\varphi_n)<+\infty$;
    \item if $\sum_{n=1}^\infty \varphi_n(x)<+\infty$ then
      $f(x)=\sum_{n=1}^\infty f_n(x)$.
  \end{enumerate}
  For each $n$ put
  \[
    f_n^+:=f_n\lor0,
    \qquad
    f_n^-:=(-f_n)\lor0 .
  \]
  Since $L^1$ is a Riesz space we have $f_n^+,f_n^-\in L^1$, and
  \[
    f_n=f_n^+-f_n^-,
    \qquad
    |f_n|=f_n^++f_n^- .
  \]
  Moreover, by positivity of $I^1$ and $|f_n|\le\varphi_n$,
  \[
    \begin{aligned}
      \sum_{n=1}^\infty I^1(f_n^+)
      +\sum_{n=1}^\infty I^1(f_n^-)
      &=\sum_{n=1}^\infty I^1(|f_n|)\\
      &\le\sum_{n=1}^\infty I^1(\varphi_n)
      <+\infty.
    \end{aligned}
  \]

  Now define functions $p,q:X\to\mathbb{R}$ by
  \[
    p(x):=
    \begin{cases}
      \sum_n f_n^+(x) & \text{if $\sum_n f_n^+(x)$ converges},\\
      0 & \text{otherwise}
    \end{cases}
  \]
  and
  \[
    q(x):=
    \begin{cases}
      \sum_n f_n^-(x) & \text{if $\sum_n f_n^-(x)$ converges},\\
      0 & \text{otherwise.}
    \end{cases}
  \]
  By Lemma~\ref{lem:positive-series-flattening} we have $p,q\in L^1$.

  Put
  \[
    A:=\left\{x\in X\mathrel{}\middle|\mathrel{}\sum_{n=1}^\infty\varphi_n(x)<+\infty\right\} .
  \]
  If $x\in A$ then
  \[
    \sum_{n=1}^\infty|f_n(x)|
    \le
    \sum_{n=1}^\infty\varphi_n(x)
    <+\infty ,
  \]
  so the two series defining $p(x)$ and $q(x)$ both converge. Hence, by the
  condition on the original series representation,
  \[
    f(x)=\sum_{n=1}^\infty f_n(x)=p(x)-q(x)
    \qquad(x\in A).
  \]
  Putting
  \[
    r:=f-(p-q) ,
  \]
  the function $r$ vanishes on $A$.

  Choose $\varepsilon_n>0$ for each $n$ so that
  \[
    \sum_{n=1}^\infty\varepsilon_n<+\infty .
  \]
  Since $\varphi_n\ge0$, by Lemma~\ref{lem:positive-good-representation} we may
  choose for each $n$ a series representation
  \[
    \varphi_n\overset{\psi_{n,k}}{\simeq}\sum_k g_{n,k}
  \]
  of $\varphi_n$ satisfying
  \[
    \sum_{k=1}^\infty I(\psi_{n,k})
    \le I^1(\varphi_n)+\varepsilon_n .
  \]
  Hence
  \[
    \sum_{n,k}I(\psi_{n,k})
    \le
    \sum_n I^1(\varphi_n)+\sum_n\varepsilon_n
    <+\infty.
  \]

  For each $n$,
  \[
    0\overset{\psi_{n,k}}{\simeq}\sum_k0
  \]
  is a series representation of the function constantly equal to $0$. Put
  \[
    B:=
    \left\{
      x\in X
      \mathrel{}\middle|\mathrel{}
      \sum_{n,k}\psi_{n,k}(x)<+\infty
    \right\} .
  \]
  If $x\in B$ then $\sum_k\psi_{n,k}(x)<+\infty$ for each $n$, so
  \[
    \varphi_n(x)=\sum_k g_{n,k}(x)
  \]
  and
  \[
    0\le\varphi_n(x)
    \le\sum_k|g_{n,k}(x)|
    \le\sum_k\psi_{n,k}(x).
  \]
  Hence
  \[
    \sum_n\varphi_n(x)
    \le\sum_{n,k}\psi_{n,k}(x)
    <+\infty.
  \]
  Therefore $x\in A$, and
  \[
    r(x)=0=\sum_n0
    \qquad(x\in B).
  \]
  Applying Lemma~\ref{lem:flatten-series-representations} with $u_n=0$ and
  $u=r$, we get $r\in L^1$.

  Since $p,q,r\in L^1$ and $L^1$ is a vector space,
  \[
    f=p-q+r\in L^1 .
  \]
  Therefore $(L^1)^1\subset L^1$, and $(L^1)^1=L^1$ is proved.
\end{proof}

\subsubsection{Splitting Daniell closedness into two closure conditions}

Daniell closedness $L^1=L$ is a single closure condition passing through the
extension by series representations.
In this subsubsection we show that this condition splits into two independent
conditions, \textbf{closure under the monotone convergence theorem} and
\textbf{closure under modification on null sets}
(Theorem~\ref{thm:complete-intrinsic}).
That neither condition can be dropped is seen in
Remark~\ref{rem:two-conditions-independent}.

The two conditions on the right are conditions that occur independently in
measure theory.
The monotone convergence theorem corresponds directly, and closure under
modification on null sets corresponds to completeness of the measure.
This decomposition is important for gauging the content of the main theorem of
Section~\ref{chap:functoriality}
(Remark~\ref{rem:main-theorem-not-tautology}).

\begin{definition}\label{def:null-set}
  Let $(X,L,I)$ be an integration system on a set $X$.
  We say that a subset $Z\subset X$ is an \textbf{$I$-null set} if there exists a
  sequence of functions $\{g_n\}$ in $L$ satisfying the following three
  conditions.
  \begin{enumerate}
    \item $g_n\ge0$ for every $n$;
    \item $\sum_{n=1}^{\infty}I(g_n)<+\infty$;
    \item $\sum_{n=1}^{\infty}g_n(x)=+\infty$ for every $x\in Z$.
  \end{enumerate}
\end{definition}

\noindent
This definition does not mention $L^1$.
Note that it is written using only the data given by $(X,L,I)$.
The set of all points at which condition (3) of a series representation
(Definition~\ref{def:series-representation}) does not constrain the value is
always an $I$-null set, as one sees by taking $\{g_n\}$ to be the dominating
functions $\{\varphi_n\}$.

\begin{lemma}\label{lem:null-sets-countable-union}
  A countable union of $I$-null sets is an $I$-null set.
\end{lemma}
\begin{proof}
  Let $Z_1,Z_2,\dots$ be $I$-null sets, choose for each $Z_k$ a sequence
  $\{g^{(k)}_n\}_n$ as in Definition~\ref{def:null-set}, and put
  $C_k:=\sum_nI(g^{(k)}_n)<+\infty$.
  Put $a_k:=2^{-k}/(1+C_k)>0$ and $\tilde{g}^{(k)}_n:=a_kg^{(k)}_n\in L$.
  Arranging the family $\{\tilde{g}^{(k)}_n\}_{k,n}$ into a single sequence,
  every term is nonnegative and
  \[
    \sum_{k,n}I(\tilde{g}^{(k)}_n)=\sum_ka_kC_k\le\sum_k2^{-k}=1<+\infty .
  \]
  Moreover, if $x\in\bigcup_kZ_k$ then $x\in Z_k$ for some $k$ and
  \[
    \sum_n\tilde{g}^{(k)}_n(x)=a_k\sum_ng^{(k)}_n(x)=+\infty .
  \]
  The remaining terms are nonnegative, so the sum of the rearranged sequence is
  $+\infty$ as well.
\end{proof}

\begin{lemma}\label{lem:vanish-off-null-set}
  Let $(X,L,I)$ be an integration system and let $Z\subset X$ be an $I$-null set.
  If $g\in L$ vanishes on $X\setminus Z$ then
  \[
    I(|g|)=0,
    \qquad\text{and in particular}\qquad
    I(g)=0 .
  \]
\end{lemma}
\begin{proof}
  Choose $\{g_n\}$ as in Definition~\ref{def:null-set} and put
  $C:=\sum_nI(g_n)<+\infty$.
  Take an arbitrary $\varepsilon>0$ and put
  \[
    u_N:=|g|\land\Bigl(\varepsilon\sum_{n\le N}g_n\Bigr) .
  \]
  Since $L$ is a Riesz space we have $u_N\in L$.

  We check that $u_N$ increases monotonically to $|g|$ pointwise.
  If $x\notin Z$ then $|g|(x)=0$, and the second term is nonnegative, so
  $u_N(x)=0=|g|(x)$.
  If $x\in Z$ then $\varepsilon\sum_{n\le N}g_n(x)\uparrow+\infty$, so
  $u_N(x)=|g|(x)$ for $N$ large enough, while before that
  $u_N(x)=\varepsilon\sum_{n\le N}g_n(x)$, which is nondecreasing.
  In either case $u_N(x)\uparrow|g|(x)$.

  Since $|g|\in L$, Lemma~\ref{lem:integration-is-continuous} gives
  $I(u_N)\uparrow I(|g|)$.
  On the other hand $u_N\le\varepsilon\sum_{n\le N}g_n$, so
  Lemma~\ref{lem:integration-preserve-order} gives
  \[
    I(u_N)\le\varepsilon\sum_{n\le N}I(g_n)\le\varepsilon C.
  \]
  Hence $I(|g|)\le\varepsilon C$. As $\varepsilon>0$ was arbitrary, $I(|g|)=0$,
  and $|I(g)|\le I(|g|)=0$ by
  Lemma~\ref{lem:triangle-ineq-for-integration}.
\end{proof}

\begin{definition}\label{def:mod-closed}
  We say that an integration system $(X,L,I)$ on a set $X$ is
  \textbf{closed under modification on null sets} if
  \[
    h|_{X\setminus Z}=f|_{X\setminus Z}
    \quad\Longrightarrow\quad
    h\in L
  \]
  holds for every $f\in L$, every $I$-null set $Z\subset X$ and every function
  $h:X\to\mathbb{R}$.
\end{definition}

\noindent
On $Z$ the function $h$ is subject to no constraint at all; only that it is
real-valued.

\begin{theorem}[decomposition of Daniell closedness]\label{thm:complete-intrinsic}
  Let $(X,L,I)$ be an integration system on a set $X$. The following are
  equivalent.
  \begin{enumerate}
    \item[(i)] $(X,L,I)$ is Daniell closed, that is, $L^1=L$.
    \item[(ii)] $(X,L,I)$ satisfies the monotone convergence theorem and is
      closed under modification on null sets.
  \end{enumerate}
\end{theorem}
\begin{proof}
  \textbf{(i)$\Rightarrow$(ii)}.
  The monotone convergence theorem was seen in
  Remark~\ref{rem:complete-implies-mct}.
  We show closure under modification on null sets.
  Let $f\in L$, let $Z$ be an $I$-null set, and let $h:X\to\mathbb{R}$ agree with
  $f$ on $X\setminus Z$.
  Choose $\{g_n\}$ as in Definition~\ref{def:null-set} and put
  \[
    f_1:=f,\quad f_n:=0\ (n\ge2),
    \qquad
    \varphi_1:=|f|,\quad \varphi_{n+1}:=g_n\ (n\ge1) .
  \]
  Then $|f_n|\le\varphi_n$ is clear, and
  \[
    \sum_{n}I(\varphi_n)=I(|f|)+\sum_nI(g_n)<+\infty .
  \]
  Moreover, if $\sum_n\varphi_n(x)<+\infty$ then $\sum_ng_n(x)<+\infty$, so
  $x\notin Z$ and therefore $h(x)=f(x)=\sum_nf_n(x)$.
  Hence $\{f_n\},\{\varphi_n\}$ is a series representation of $h$, and
  $h\in L^1=L$.

  \textbf{(ii)$\Rightarrow$(i)}.
  Since $L\subset L^1$ by
  Lemma~\ref{lem:daniell-ext-contains-source}, it suffices to show
  $L^1\subset L$.

  Let $f\in L^1$ and choose a series representation
  $f\overset{\varphi_n}{\simeq}\sum f_n$.
  Put
  \[
    Z:=\Bigl\{x\in X\;\Bigm|\;\sum_n\varphi_n(x)=+\infty\Bigr\} .
  \]
  From $|f_n|\le\varphi_n$ we get $\varphi_n\ge0$, and
  $\sum_nI(\varphi_n)<+\infty$, so $\{\varphi_n\}$ witnesses that $Z$ is an
  $I$-null set.

  Put $p_n:=f_n\lor0$ and $m_n:=(-f_n)\lor0$, so that $p_n,m_n\in L$ and
  \[
    f_n=p_n-m_n,
    \qquad
    0\le p_n,m_n\le|f_n|\le\varphi_n .
  \]
  Put
  \[
    P_N:=\sum_{n\le N}p_n\in L,
    \qquad
    M_N:=\sum_{n\le N}m_n\in L .
  \]
  Both are nondecreasing in $N$, and
  Lemma~\ref{lem:integration-preserve-order} gives
  \[
    I(P_N)\le\sum_{n\le N}I(\varphi_n)\le\sum_nI(\varphi_n)<+\infty ,
  \]
  and likewise for $M_N$.

  If $x\notin Z$ then $\sum_n\varphi_n(x)<+\infty$, so $\sum_np_n(x)$ and
  $\sum_nm_n(x)$ both converge. We therefore define
  \[
    P'(x):=
    \begin{cases}
      \sum_np_n(x) & (x\notin Z)\\
      0 & (x\in Z)
    \end{cases}
  \]
  and define $M'$ similarly. Both are real-valued functions on $X$.
  Let $P_N'$ be the function agreeing with $P_N$ on $X\setminus Z$ and vanishing
  on $Z$, and define $M_N'$ similarly.
  By hypothesis $(X,L,I)$ is closed under modification on null sets, so
  $P_N',M_N'\in L$.

  Since $P_N-P_N'$ vanishes on $X\setminus Z$,
  Lemma~\ref{lem:vanish-off-null-set} gives $I(P_N-P_N')=0$, that is,
  \[
    I(P_N')=I(P_N)\le\sum_nI(\varphi_n)<+\infty .
  \]
  Moreover $P_N'(x)\uparrow P'(x)$ for every $x\in X$ (for $x\notin Z$ by
  definition, and for $x\in Z$ both sides are $0$).
  By hypothesis $(X,L,I)$ satisfies the monotone convergence theorem, so
  $P'\in L$. Similarly $M'\in L$.

  Hence $P'-M'\in L$.
  If $x\notin Z$ then $\sum_nf_n(x)$ converges absolutely by
  Lemma~\ref{lem:series-exp-is-abs-conv}, so
  \[
    (P'-M')(x)=\sum_np_n(x)-\sum_nm_n(x)=\sum_nf_n(x)=f(x) .
  \]
  That is, $f$ agrees with the element $P'-M'$ of $L$ on $X\setminus Z$.
  Using closure under modification on null sets again, $f\in L$.
\end{proof}

\medskip
\noindent
\textbf{On the construction in (ii)$\Rightarrow$(i).}
The functions $P'$ and $M'$ built in the proof above are both monotonically
nondecreasing limits of elements of $L$ with bounded integrals, that is,
\textbf{upper functions} in the classical sense.
What the proof exhibits is therefore the representation of an $f\in L^1$ as
agreeing, off a null set, with the difference $P'-M'$ of two upper functions.
Constructing the Lebesgue integral as differences of upper functions is a
classical method, and a representation of this shape is well known in that
context.
We do not use this correspondence below, and so cite no particular source for
it.

\begin{corollary}\label{cor:mintbar-intrinsic}
  $\mathbf{DInt}$ is the full subcategory of $\mathbf{Int}$ consisting of the
  \textbf{integration systems that satisfy the monotone convergence theorem and
  are closed under modification on null sets}.
\end{corollary}
\begin{proof}
  By Theorem~\ref{thm:complete-intrinsic}.
\end{proof}

\begin{remark}\label{rem:not-about-mentioning-L1}
  The description above does not pass through the construction of the Daniell
  extension, but \textbf{that in itself is not where the point lies}.
  The condition $L^1=L$, when unfolded, also reads
  \begin{quote}
    every function admitting a series representation by elements of $L$ belongs
    to $L$,
  \end{quote}
  which uses nothing but $L$, $I$ and pointwise conditions.
  ``Whether $L^1$ is mentioned'' does not distinguish the two.

  What matters is that one unfamiliar closure condition splits into two
  conditions that occur independently in measure theory.
  We return to this point in
  Remark~\ref{rem:main-theorem-not-tautology}.
\end{remark}

\begin{example}
  Let us check this in the example of
  Remark~\ref{rem:D-I-not-left-adjoint} ($X=\{0,1,2\}$, $L=\{(a,b,b)\}$,
  $I(a,b,b)=a$).
  The function $g_n:=(0,1,1)\in L$ satisfies $g_n\ge0$ and $I(g_n)=0$, and
  $\sum_ng_n(x)=+\infty$ at $x=1,2$, so $Z=\{1,2\}$ is an $I$-null set.
  But $h=(0,1,0)$, obtained from $f=(0,0,0)\in L$ by modifying it on $Z$, does
  not belong to $L$.
  Hence this integration system is not closed under modification on null sets
  and, by Theorem~\ref{thm:complete-intrinsic}, is not Daniell closed.
  This agrees with the direct computation $L^1=\mathbb{R}^X\ne L$ in that
  example.

  On the other hand, in Example~\ref{ex:stone-fails} ($X=\{1,2\}$,
  $L=\{(a,2a)\}$, $I(a,2a)=a$), if $g_n=(c_n,2c_n)\ge0$ satisfies
  $\sum_nI(g_n)=\sum_nc_n<+\infty$ then $\sum_ng_n(x)<+\infty$ holds at both
  $x=1$ and $x=2$.
  Hence the only $I$-null set is the empty set, and the condition on
  modification on null sets holds vacuously.
  The monotone convergence theorem holds as well, so this integration system is
  Daniell closed by Theorem~\ref{thm:complete-intrinsic}, in agreement with
  $L^1=L$ as stated in that example.
\end{example}

\begin{remark}[neither condition can be dropped]\label{rem:two-conditions-independent}
  The two conditions in (ii) of Theorem~\ref{thm:complete-intrinsic} are
  independent.

  \textbf{That the monotone convergence theorem cannot be dropped} is shown by
  Example~\ref{ex:counting} ($X=\mathbb{N}$, the functions vanishing outside
  finitely many $n$, and summation).
  In this integration system, $g\in L$ with $g\ge0$ satisfies $g(x)\le I(g)$ for
  every $x$, so for a nonnegative sequence $\{g_n\}\subset L$ with
  $\sum_nI(g_n)<+\infty$ we have $\sum_ng_n(x)<+\infty$ at every $x$.
  Hence the only $I$-null set is the empty set, and the condition on
  modification on null sets holds vacuously.
  On the other hand
  \[
    f_N:=\sum_{k\le N}2^{-k}1_{\{k\}}\in L
  \]
  satisfies $I(f_N)\le1$ while increasing monotonically pointwise to
  $f(k)=2^{-k}$, and $f\notin L$ because $f$ does not have finite support.
  That is, the monotone convergence theorem fails.
  Indeed this integration system is not Daniell closed, and $L^1=\ell^1$.

  \textbf{That modification on null sets cannot be dropped} is shown by the
  example of Remark~\ref{rem:D-I-not-left-adjoint}.
  The integration system of that example satisfies the monotone convergence
  theorem but, as seen in the example above, is not closed under modification on
  null sets, hence not Daniell closed.
\end{remark}

\begin{remark}
  The two conditions in (ii) of Theorem~\ref{thm:complete-intrinsic} correspond
  to two properties that $L^1(\mu)$ satisfies in measure theory.
  The monotone convergence theorem corresponds directly, and closure under
  modification on null sets corresponds to the fact that $L^1(\mu)$ contains the
  functions measurable with respect to the completion of $\mu$.
  In Proposition~\ref{prop:l1-decomposition} of
  Subsection~\ref{sec:why-complete} we will see that the gap between the two
  reflections amounts exactly to the functions supported on null sets; the
  present theorem explains the reason for that from the side of Daniell
  closedness.
  That is, \textbf{to be an object of $\mathbf{DInt}$ one must be closed under
  modification on null sets}.
\end{remark}

\begin{remark}[the decomposition collapses in the extended-real formalism]\label{rem:extended-real-collapse}
  That the two conditions of Theorem~\ref{thm:complete-intrinsic} are
  independent \textbf{depends on this paper's restriction of functions to real
  values}.
  If the range is widened to
  $\overline{\mathbb{R}}=[-\infty,+\infty]$, then condition (b) follows from
  condition (a).

  Indeed \cite{CWS} adopts that formalism.
  It defines the counterpart of a Riesz space inside
  $\overline{\mathbb{R}}^X$ under the name Riesz lattice, and calls the triple
  corresponding to our integration system a Daniell space.
  It then gives the two conditions separate names.
  \begin{itemize}
    \item Its Definition 1.4.1 calls a triple \textbf{closed} if, for every
      monotone sequence whose integrals are bounded, the pointwise limit again
      belongs to the space and the integral commutes with the limit.
      This is our Definition~\ref{def:satisfies-mct} (satisfying the monotone
      convergence theorem).
    \item Its Definition 1.4.7 calls a triple \textbf{complete} if
      $f\le h\le g$ together with $\ell(f)=\ell(g)$ implies that $h$ belongs to
      the space as well.
      This sandwich condition is in substance our
      Definition~\ref{def:mod-closed}.
  \end{itemize}
  Its Corollary 1.4.8, however, states that
  \begin{quote}
    every closed Daniell space is complete.
  \end{quote}
  That is, \textbf{in the extended-real formalism (a)$\Rightarrow$(b), and the
  two conditions are not independent.}

  The reason is the following. Let $Z$ be an $I$-null set and take the sequence
  $\{g_n\}$ of Definition~\ref{def:null-set}.
  The partial sums $\sum_{n\le N}g_n$ are monotonically nondecreasing with
  bounded integrals, so under (a) their pointwise limit belongs to the space.
  If extended-real values are allowed, that limit is an element of the space as
  a \emph{function} taking the value $+\infty$ on $Z$.
  Null sets can then be named from inside the space, as the sets where an
  element diverges to $\infty$, and (b) follows from the axioms of a Riesz
  lattice, which originate in the fact that $\infty-\infty$ is undefined.

  \textbf{This argument does not go through here.}
  Since $L\subset\mathbb{R}^X$, the limit above cannot be an element of $L$.
  A null set cannot be captured from inside $L$ as divergence to $\infty$, and
  it is precisely for this reason that
  Definition~\ref{def:null-set} takes the form of a series and that condition
  (b) remains as an independent condition.
  The counterexample of
  Remark~\ref{rem:two-conditions-independent} ($X=\{0,1,2\}$,
  $L=\{(a,b,b)\}$, $I(a,b,b)=a$) was an integration system satisfying the
  monotone convergence theorem but not closed under modification on null sets.
  That is, it is also a \textbf{counterexample to the real-valued form of
  Corollary 1.4.8}.

  This is not a question of which formalism is better.
  The correspondence is that \textbf{widening the range makes null sets nameable
  inside the space, at the cost of collapsing the two closure conditions into
  one}.
  Theorem~\ref{thm:complete-intrinsic} is a decomposition visible only in a
  formalism restricted to real values.
\end{remark}

\subsection{The $L^1$ seminorm and approximation}

The content of this subsection is also known, and corresponds to \S4 and \S6 of
\cite{BM}.
We give proofs because it is used essentially in
Section~\ref{chap:operator}.

\begin{definition}\label{def:l1-seminorm}
  Let $(X,L,I)$ be an integration system on a set $X$.
  For each $f\in L$ we put
  \[
    \|f\|_I:=I(|f|) .
  \]
  Since $L$ is a Riesz space we have $|f|\in L$, so this value is defined.
\end{definition}

\begin{lemma}\label{lem:l1-seminorm}
  $\|\cdot\|_I$ is a seminorm on $L$. That is, for all $f,g\in L$ and
  $a\in\mathbb{R}$,
  \[
    \|af\|_I=|a|\,\|f\|_I,
    \qquad
    \|f+g\|_I\le\|f\|_I+\|g\|_I.
  \]
\end{lemma}
\begin{proof}
  The first identity follows from $|af|=|a||f|$ and linearity of $I$.
  Moreover $|f+g|\le|f|+|g|$ and
  Lemma~\ref{lem:integration-preserve-order} give
  $I(|f+g|)\le I(|f|)+I(|g|)$.
\end{proof}

\begin{remark}
  In general $\|f\|_I=0$ does not imply $f=0$.
  The operation that collapses this discrepancy is treated in
  Subsection~\ref{sec:ae}.
\end{remark}

\begin{lemma}[approximation of integrable functions by $L$]\label{lem:l1-approx-by-source}
  Let $(X,L,I)$ be an integration system on a set $X$ and let
  $(X,L^1,I^1)$ be its Daniell extension.
  Then for every $f\in L^1$ there is a sequence of functions
  $\{s_N\}\subset L$ satisfying
  \[
    \|f-s_N\|_{I^1}=I^1(|f-s_N|)\longrightarrow0
    \qquad(N\to\infty) .
  \]
  That is, $L$ is dense in $L^1$ with respect to $\|\cdot\|_{I^1}$.
\end{lemma}
\begin{proof}
  Choose a series representation of $f$, say
  \[
    f\overset{\varphi_n}{\simeq}\sum_{n=1}^{\infty}f_n ,
  \]
  and put
  \[
    s_N:=\sum_{n=1}^Nf_n\in L,
    \qquad
    r_N:=f-s_N .
  \]
  Since finitely many values $\varphi_n(x)$ are always finite,
  \[
    r_N\overset{\varphi_{N+m}}{\simeq}
    \sum_{m=1}^{\infty}f_{N+m}
  \]
  is a series representation of $r_N$.

  For each $m\ge1$ put
  \[
    g_{N,m}:=
    \left|\sum_{k=1}^{m}f_{N+k}\right|
    -\left|\sum_{k=1}^{m-1}f_{N+k}\right| ,
  \]
  where the empty sum is $0$. Since $L$ is a Riesz space we have
  $g_{N,m}\in L$, and the triangle inequality for absolute values gives
  \[
    |g_{N,m}|\le |f_{N+m}|\le\varphi_{N+m} .
  \]
  Moreover, if $\sum_m\varphi_{N+m}(x)<+\infty$ then, by telescoping the sum and
  continuity of the absolute value,
  \[
    \sum_{m=1}^{\infty}g_{N,m}(x)
    =\left|\sum_{m=1}^{\infty}f_{N+m}(x)\right|
    =|r_N(x)| .
  \]
  Hence
  \[
    |r_N|\overset{\varphi_{N+m}}{\simeq}
    \sum_{m=1}^{\infty}g_{N,m}
  \]
  is a series representation of $|r_N|$. By positivity of the integral,
  \[
    \begin{aligned}
      0\le I^1(|f-s_N|)
      &=\sum_{m=1}^{\infty}I(g_{N,m})\\
      &\le\sum_{m=1}^{\infty}I(|g_{N,m}|)\\
      &\le\sum_{m=1}^{\infty}I(\varphi_{N+m})
      \longrightarrow0 ,
    \end{aligned}
  \]
  which is the assertion.
\end{proof}

\begin{corollary}\label{cor:tail-estimate}
  Let $(X,L,I)$ be an integration system and suppose $\{h_j\}\subset L$
  satisfies $\sum_j I(|h_j|)<+\infty$.
  Then the function
  \[
    h(x):=\begin{cases}
      \sum_j h_j(x) & \text{if $\sum_j|h_j(x)|<+\infty$}\\
      0 & \text{otherwise}
    \end{cases}
  \]
  belongs to $L^1$, and for each $J\ge1$
  \[
    \left\|h-\sum_{j=1}^{J}h_j\right\|_{I^1}\le\sum_{j>J}I(|h_j|)
  \]
  holds.
\end{corollary}
\begin{proof}
  Putting $\varphi_j:=|h_j|\in L$, the pair $\{h_j\},\{\varphi_j\}$ is a series
  representation of $h$, so $h\in L^1$.
  Moreover $h-\sum_{j\le J}h_j$ has $\{h_{J+m}\}_{m\ge1}$ and
  $\{\varphi_{J+m}\}_{m\ge1}$ as a series representation, so the same estimate
  as at the end of the proof of Lemma~\ref{lem:l1-approx-by-source} gives the
  stated inequality.
\end{proof}

\begin{theorem}[$\|\cdot\|$ completeness]\label{thm:l1-norm-complete}
  Let $(X,L,I)$ be a Daniell closed integration system.
  Then $(L,\|\cdot\|_I)$ is complete.
  That is, for every Cauchy sequence $\{g_k\}\subset L$ with respect to
  $\|\cdot\|_I$ there is a $g\in L$ with $\|g-g_k\|_I\to0$.
\end{theorem}
\begin{proof}
  Let $\{g_k\}$ be a Cauchy sequence. Choose a subsequence $\{g_{k_j}\}$ so that
  \[
    \|g_{k_{j+1}}-g_{k_j}\|_I<2^{-j} ,
  \]
  and put
  \[
    h_1:=g_{k_1},\qquad h_{j+1}:=g_{k_{j+1}}-g_{k_j} .
  \]
  Then $\sum_j I(|h_j|)<+\infty$.
  By Corollary~\ref{cor:tail-estimate} there is a corresponding $h\in L^1$ with
  \[
    \left\|h-g_{k_J}\right\|_{I^1}
    =\left\|h-\sum_{j=1}^{J}h_j\right\|_{I^1}
    \le\sum_{j>J}I(|h_j|)
    \longrightarrow0 .
  \]
  Since $(X,L,I)$ is Daniell closed we have $L^1=L$, and
  $I^1|_L=I$ by Lemma~\ref{lem:daniell-ext-contains-source}.
  Hence $h\in L$ and $\|h-g_{k_J}\|_I\to0$.
  A Cauchy sequence with a convergent subsequence converges to the same limit as
  a whole, so $\|h-g_k\|_I\to0$.
\end{proof}

\subsection{Null functions and equality almost everywhere}\label{sec:ae}

The seminorm $\|\cdot\|_I$ is in general not a norm.
We set up here the operation that collapses its kernel.
This quotient is used essentially in Section~\ref{chap:operator}.

\subsubsection{Null functions and the quotient}

\begin{definition}\label{def:null-function}
  Let $(X,L,I)$ be an integration system on a set $X$. We put
  \[
    \mathcal{N}_I:=\{f\in L\mid I(|f|)=0\}
  \]
  and call the functions belonging to $\mathcal{N}_I$ the \textbf{null functions}
  of $(X,L,I)$.
\end{definition}

\begin{lemma}\label{lem:null-is-ideal}
  $\mathcal{N}_I$ is an ideal of $L$.
  That is, $\mathcal{N}_I$ is a linear subspace, and if $f\in\mathcal{N}_I$ and
  $f'\in L$ satisfies $|f'|\le|f|$ then $f'\in\mathcal{N}_I$.
\end{lemma}
\begin{proof}
  $\mathcal{N}_I$ is the kernel of the seminorm $\|\cdot\|_I$, hence a linear
  subspace (Lemma~\ref{lem:l1-seminorm}).
  Moreover, if $|f'|\le|f|$ then
  Lemma~\ref{lem:integration-preserve-order} gives $I(|f'|)\le I(|f|)=0$, and
  positivity of $I$ gives $I(|f'|)=0$.
\end{proof}

\begin{definition}\label{def:null-quotient}
  We write the quotient space as
  \[
    \widetilde{L}:=L/\mathcal{N}_I
  \]
  and denote the quotient map by $\pi_I:L\to\widetilde{L}$.
\end{definition}

\begin{lemma}\label{lem:quotient-riesz}
  $\widetilde{L}$ carries the structure of a Riesz space for the order defined by
  declaring $\pi_I(f)\ge0$ to mean ``there is an $f'\in\mathcal{N}_I$ with
  $f+f'\ge0$'', and $\pi_I$ is a Riesz homomorphism.
  That is,
  \[
    \pi_I(f\vee g)=\pi_I(f)\vee\pi_I(g),
    \qquad
    \pi_I(|f|)=|\pi_I(f)| .
  \]
  Furthermore
  \[
    \widetilde{I}:\widetilde{L}\to\mathbb{R};\quad \pi_I(f)\mapsto I(f)
  \]
  is a well-defined positive linear functional, and
  $\|\pi_I(f)\|:=\widetilde{I}(|\pi_I(f)|)=\|f\|_I$ is a norm.
\end{lemma}
\begin{proof}
  That the quotient of a Riesz space by an ideal is a Riesz space and that the
  quotient map is a Riesz homomorphism is general theory (see for instance
  \cite{LuxemburgZaanen}).
  By Lemma~\ref{lem:null-is-ideal} this applies here.

  $\widetilde{I}$ is well defined because $f\in\mathcal{N}_I$ implies
  $|I(f)|\le I(|f|)=0$
  (Lemma~\ref{lem:triangle-ineq-for-integration}).
  Positivity and linearity follow from those of $I$.
  Finally
  $\widetilde{I}(|\pi_I(f)|)=\widetilde{I}(\pi_I(|f|))=I(|f|)=\|f\|_I$, and this
  vanishes only when $f\in\mathcal{N}_I$, that is, when $\pi_I(f)=0$.
\end{proof}

\begin{lemma}\label{lem:quotient-is-banach}
  If $(X,L,I)$ is a Daniell closed integration system then
  $(\widetilde{L},\|\cdot\|)$ is a Banach space.
\end{lemma}
\begin{proof}
  If $\{\pi_I(f_k)\}$ is Cauchy then $\{f_k\}$ is a Cauchy sequence with respect
  to $\|\cdot\|_I$.
  By Theorem~\ref{thm:l1-norm-complete} its limit $f\in L$ exists, and
  $\|\pi_I(f)-\pi_I(f_k)\|=\|f-f_k\|_I\to0$.
\end{proof}

\begin{lemma}[compatibility of the Daniell extension with the quotient]\label{lem:quotient-embedding}
  Let $(X,L,I)$ be an integration system on a set $X$ and let
  $(X,L^1,I^1)$ be its Daniell extension.
  Then
  \[
    \mathcal{N}_{I^1}\cap L=\mathcal{N}_I .
  \]
  Hence the inclusion $L\subset L^1$ induces an injective Riesz homomorphism
  \[
    j_I:\widetilde{L}\longrightarrow\widetilde{L^1}
  \]
  satisfying
  \[
    \pi_{I^1}|_L=j_I\circ\pi_I .
  \]
  Furthermore $j_I$ preserves the norm.
\end{lemma}
\begin{proof}
  By Lemma~\ref{lem:daniell-ext-contains-source} we have $|f|\in L$ and
  $I^1(|f|)=I(|f|)$ for every $f\in L$.
  Hence $I^1(|f|)=0$ and $I(|f|)=0$ are equivalent, and we obtain
  $\mathcal{N}_{I^1}\cap L=\mathcal{N}_I$.

  It follows that the kernel of $\pi_{I^1}|_L:L\to\widetilde{L^1}$ is
  $\mathcal{N}_I$, so an injective linear map
  $j_I:\widetilde{L}\to\widetilde{L^1}$ is induced, satisfying
  $\pi_{I^1}|_L=j_I\circ\pi_I$.
  Both $\pi_I$ and $\pi_{I^1}$ are Riesz homomorphisms
  (Lemma~\ref{lem:quotient-riesz}) and $\pi_I$ is surjective, so $j_I$ is a
  Riesz homomorphism as well.
  Finally, for $f\in L$,
  \[
    \|j_I(\pi_I(f))\|=\|\pi_{I^1}(f)\|=I^1(|f|)=I(|f|)=\|\pi_I(f)\| .
  \]
\end{proof}

\subsubsection{ae-equivalence of morphisms}

\begin{definition}\label{def:ae-equal}
  Let $(X,L,I)$ and $(Y,M,J)$ be integration systems on sets $X$ and
  $Y$ respectively.
  We say that two morphisms of integration systems
  \[
    F,G:(X,L,I)\to(Y,M,J)
  \]
  are \textbf{equal almost everywhere} if
  \[
    J\bigl(|F^*f-G^*f|\bigr)=0
  \]
  holds for every $f\in L$, that is, if $F^*f-G^*f\in\mathcal{N}_J$ holds for
  every $f\in L$.
  In this case we write
  \[
    F=_{\mathrm{ae}}G .
  \]
\end{definition}

\begin{lemma}\label{lem:ae-is-congruence}
  For the morphisms of integration systems between each pair of objects, the
  relation $=_{\mathrm{ae}}$ is an equivalence relation and is compatible with
  composition of morphisms.
  That is, if morphisms of integration systems
  \[
    F,F':(X,L,I)\to(Y,M,J),
    \qquad
    G,G':(Y,M,J)\to(Z,N,K)
  \]
  satisfy
  \[
    F=_{\mathrm{ae}}F',
    \qquad
    G=_{\mathrm{ae}}G' ,
  \]
  then
  \[
    G\circ F=_{\mathrm{ae}}G'\circ F' .
  \]
\end{lemma}
\begin{proof}
  Reflexivity and symmetry are immediate from the definition. As for
  transitivity, suppose morphisms of integration systems
  $F,F',F'':(X,L,I)\to(Y,M,J)$ satisfy $F=_{\mathrm{ae}}F'$ and
  $F'=_{\mathrm{ae}}F''$. For every $f\in L$,
  \[
    \begin{aligned}
      0
      &\le J\bigl(|F^*f-(F'')^*f|\bigr)\\
      &\le J\bigl(|F^*f-(F')^*f|\bigr)
        +J\bigl(|(F')^*f-(F'')^*f|\bigr)
      =0 ,
    \end{aligned}
  \]
  so $F=_{\mathrm{ae}}F''$.

  Next we show compatibility with composition. For every $f\in L$, since
  morphisms of Riesz systems preserve absolute values and morphisms of
  integration systems preserve the integral,
  \[
    \begin{aligned}
      &K\bigl(|G^*F^*f-(G')^*(F')^*f|\bigr)\\
      &\quad\le
      K\bigl(|G^*F^*f-G^*(F')^*f|\bigr)
      +K\bigl(|G^*(F')^*f-(G')^*(F')^*f|\bigr)\\
      &\quad=
      J\bigl(|F^*f-(F')^*f|\bigr)
      +K\bigl(|G^*(F')^*f-(G')^*(F')^*f|\bigr)
      =0.
    \end{aligned}
  \]
  Therefore $G\circ F=_{\mathrm{ae}}G'\circ F'$.
\end{proof}

\begin{definition}\label{def:int-ae}
  The category having the same objects as the category $\mathbf{Int}$ of
  integration systems and whose morphisms are the equivalence classes under
  $=_{\mathrm{ae}}$ is written
  \[
    \mathbf{Int}_{\mathrm{ae}}
  \]
  and called the \textbf{ae-category of integration systems}.
  By Lemma~\ref{lem:ae-is-congruence}, composition of morphisms is defined
  independently of the choice of representatives.

  Among the objects of $\mathbf{Int}_{\mathrm{ae}}$, we write the full
  subcategory consisting of the integration systems satisfying the monotone
  convergence theorem as
  \[
    \mathbf{MInt}_{\mathrm{ae}},
  \]
  and the full subcategory consisting of the Daniell closed integration systems
  as
  \[
    \mathbf{DInt}_{\mathrm{ae}} .
  \]
\end{definition}

\subsection{Relation to measures}\label{sec:measure}

Up to here this paper has treated integration without using measures.
In this subsection we state the correspondence between integration systems and
measure theory.
The main results of this paper do not use this correspondence, but we record that
the notion of an integration system does recover measure theory.
Everything below is known, and we refer to the literature for proofs.

\subsubsection{From measures to integration systems}

This direction is the easy one.

\begin{proposition}\label{prop:measure-gives-integral-system}
  Let $(X,\Sigma,\mu)$ be a measure space, let $L$ be the set of all
  $\Sigma$-measurable simple functions that are integrable with respect to
  $\mu$, and let $I(f):=\int_Xf\,d\mu$.
  Then $(X,L,I)$ is an integration system.
\end{proposition}

\noindent
Positivity is immediate, and verifying the continuity axiom uses
$\sigma$-additivity of $\mu$.
For a proof see \cite{BM}, \S7.
The set of step functions in Example~\ref{ex:step-functions} is a subsystem of
the set of all simple functions with respect to Lebesgue measure.
Finite unions of half-open intervals do not form a $\sigma$-algebra, so it is not
itself a special case of this proposition.

\subsubsection{The Stone condition}

The converse direction, that of extracting a measure from an integration system,
needs one condition.

\begin{definition}\label{def:stone-condition}
  We say that an integration system $(X,L,I)$ on a set $X$ satisfies
  the \textbf{Stone condition} if
  \[
    f\wedge1\in L
  \]
  holds for every $f\in L$. Here $1$ is the constant function $1$ on $X$.
\end{definition}

\noindent
Note that $L$ itself need not contain $1$.
What is required is that the function obtained by truncating $f$ at $1$ belong to
$L$ again.
This condition is the standard one, and a Riesz space satisfying it is called a
\textbf{Stone vector lattice}.

This condition is needed because we want to define the measure of a set $A$ with
$1_A\in L^1$ by
\[
  \mu(A):=I^1(1_A) .
\]
If $L^1$ contains no indicator function at all, this definition cannot even be
started.

\begin{remark}\label{rem:mu-only-on-finite}
  The formula above is meaningful only for sets with $1_A\in L^1$.
  It cannot be used for general measurable sets, and in particular not for sets
  of infinite measure.

  Looking at Example~\ref{ex:counting} ($X=\mathbb{N}$, the functions of finite
  support, and summation), $L^1$ is the set $\ell^1$ of all absolutely summable
  sequences and the measure one should obtain is the counting measure.
  But $1_{\mathbb{N}}\notin\ell^1$, so $I^1(1_{\mathbb{N}})$ is undefined, while
  $\mu(\mathbb{N})=\infty$.

  To construct a general measure one must therefore start from the sets with
  $1_A\in L^1$ and pass through an upper integral or an outer measure.
  Theorem~\ref{thm:daniell-stone} below states only the result of that
  construction.
\end{remark}

\begin{example}\label{ex:stone-holds}
  Both Example~\ref{ex:counting} and Example~\ref{ex:step-functions} satisfy the
  Stone condition.
  Indeed, truncating a function of finite support at $1$ again has finite
  support, and truncating a step function at $1$ is again a step function.
\end{example}

\begin{example}[an integration system not satisfying the Stone condition]\label{ex:stone-fails}
  Let $X=\{1,2\}$ and put
  \[
    L:=\{(a,2a)\mid a\in\mathbb{R}\}\subset\mathbb{R}^X,
    \qquad
    I(a,2a):=a .
  \]
  Then $L$ is the one-dimensional subspace spanned by $(1,2)$, and
  $|(a,2a)|=|a|\cdot(1,2)\in L$, so $L$ is a Riesz space on $X$.
  Moreover $I$ is clearly positive linear, and continuity holds because $X$ is
  finite.
  Hence $(X,L,I)$ is an integration system.

  But $f=(1,2)$ gives $f\wedge1=(1,1)\notin L$, so the Stone condition fails.

  In this integration system $L^1=L$.
  Indeed, the dominating functions occurring in a series representation have the
  form $\varphi_n=(c_n,2c_n)$ with $c_n\ge0$, and if
  $\sum_nI(\varphi_n)=\sum_nc_n<+\infty$ then $\sum_n\varphi_n(x)<+\infty$ holds
  at both $x=1$ and $x=2$.
  Hence condition (3) of a series representation is imposed on all of $X$, and
  $f$ is confined to being the sum of an absolutely convergent series of
  elements of $L$, that is, an element of $L$.

  Consequently $1_{\{1\}}=(1,0)\notin L^1$, and $\mu(\{1\}):=I^1(1_{\{1\}})$
  cannot be defined.
\end{example}

\subsubsection{The Daniell--Stone representation theorem}

\begin{theorem}[Daniell--Stone]\label{thm:daniell-stone}
  Let $(X,L,I)$ be an integration system satisfying the Stone condition.
  Then there exist a $\sigma$-algebra $\Sigma$ on $X$ and a measure $\mu$ such
  that every $f\in L$ is $\mu$-integrable and
  \[
    I(f)=\int_Xf\,d\mu .
  \]
\end{theorem}

\noindent
For a proof see \cite{Stone}, \cite{Royden} Chapter 16, or \cite{CFW}.
Chapter 16 of \cite{Royden} is entitled ``The Daniell Integral''; its \S2 is the
extension theorem, \S3 uniqueness, and \S4 measurability and measure.
Here $\Sigma$ is chosen so that the functions belonging to $L^1$ are measurable,
and $\mu$ is constructed through an outer measure, as in
Remark~\ref{rem:mu-only-on-finite}, so as to satisfy $\mu(A)=I^1(1_A)$ for the
sets $A$ with $1_A\in L^1$.

\begin{example}
  Applying this theorem to Example~\ref{ex:step-functions}, the measure obtained
  is Lebesgue measure on $[0,1)$ and the Daniell extension $L^1$ coincides with
  the set of all Lebesgue integrable functions.
  This is how the fact stated in Example~\ref{ex:lebesgue} looks from the side
  of measure theory.
\end{example}

\begin{remark}
  This paper does not impose the Stone condition.
  The results of Section~\ref{chap:functoriality} and
  Section~\ref{chap:operator} hold without it.
  Indeed the integration system of Example~\ref{ex:stone-fails} is an object of
  $\mathbf{Int}$ as well, and is moreover Daniell closed.

  One cannot, however, conclude from this that ``the framework of this paper is
  wider than measure theory''.
  The integral of Example~\ref{ex:stone-fails} can be written as
  $I(a,2a)=\int_Xf\,d\mu$ by taking $\mu(\{1\})=1$ and $\mu(\{2\})=0$
  (and such a $\mu$ is moreover not unique: any $\mu$ with
  $\mu(\{1\})+2\mu(\{2\})=1$ will do).
  What that example shows is that without the Stone condition one
  \textbf{cannot canonically recover a measure from inside the Daniell extension
  using indicator functions}, and that \textbf{the representing measure is not
  uniquely determined}.

  Put precisely, this paper does not assume the Stone condition and does not take
  representation by a measure as a premise of the theory.
\end{remark}

\begin{problem}\label{prob:stone-measure-spaces}
  Is every Daniell closed integration system satisfying the Stone condition
  realized as the $L^1$ of a measure space, and does this correspondence lift to
  an equivalence of categories?
  That is, is the full subcategory of $\mathbf{DInt}$ consisting of the objects
  that satisfy the Stone condition equivalent to a suitable category of measure
  spaces and measure-preserving maps?
\end{problem}

\noindent
This question was raised by T.~Fritz in personal communication (2026).
It does not hold, however, in the form of an equivalence with the category of
measurable spaces.
A morphism of integration systems preserves the integral, $J(F^*f)=I(f)$, so the
integral itself is data that distinguishes objects up to isomorphism.
Indeed, let $X=\{*\}$, $L=\mathbb{R}$ and $I_c(a):=ca$ for $c>0$.
Since $L^1\subset\mathbb{R}^X=\mathbb{R}=L$, each of these is a Daniell closed
integration system satisfying the Stone condition.
But $(X,\mathbb{R},I_c)\not\cong(X,\mathbb{R},I_{c'})$ for $c\ne c'$, while the
underlying measurable space is a one-point space in every case.
What should be asked about is therefore measure spaces, not measurable spaces.
Moreover, by Theorem~\ref{thm:complete-intrinsic} a Daniell closed integration
system is closed under modification on null sets, so the measures to be matched
are expected to be the complete ones.

\section{Functoriality and universality of the Daniell extension}\label{chap:functoriality}

\subsection{Functoriality of the Daniell extension}

The inclusion functors $\iota,\iota_{D},\iota_M$ are as in
Definition~\ref{def:inclusion-functors}.

\begin{lemma}\label{lem:daniell-ext-is-functor-1}
  Let $(X,L,I)$ and $(Y,M,J)$ be integration systems on sets $X$ and
  $Y$ respectively, and let $F:(X,L,I)\to(Y,M,J)$ be a morphism of integration
  systems.
  Then the map of base spaces $\underline{F}:Y\to X$ determines a morphism
  \[
    F^1:(X,L^1,I^1)\to(Y,M^1,J^1)
  \]
  in the category $\mathbf{DInt}$ of Daniell closed integration systems.
\end{lemma}
\begin{proof}
  Let $f\in L^1$ and let
  \[
    f\overset{\varphi_n}{\simeq}\sum f_n
  \]
  be a series representation of it.
  Since $F$ is a morphism of integration systems, for each $n$
  \[
    F^*f_n,\ F^*\varphi_n\in M .
  \]
  Moreover, for every $y\in Y$,
  \[
    |F^*f_n(y)|
    =|f_n(\underline{F}(y))|
    \le\varphi_n(\underline{F}(y))
    =F^*\varphi_n(y)
  \]
  and
  \[
    \sum_{n=1}^\infty J(F^*\varphi_n)
    =\sum_{n=1}^\infty I(\varphi_n)
    <+\infty.
  \]

  Furthermore, if $y\in Y$ satisfies
  \[
    \sum_{n=1}^\infty F^*\varphi_n(y)<+\infty
  \]
  then, putting $x:=\underline{F}(y)$,
  \[
    \sum_{n=1}^\infty\varphi_n(x)<+\infty .
  \]
  By the condition on the original series representation,
  \[
    \begin{aligned}
      F^*f(y)
      &=f(\underline{F}(y))\\
      &=\sum_{n=1}^\infty f_n(\underline{F}(y))\\
      &=\sum_{n=1}^\infty F^*f_n(y) .
    \end{aligned}
  \]
  Hence
  \[
    F^*f\overset{F^*\varphi_n}{\simeq}\sum F^*f_n
  \]
  is a series representation of $F^*f$, and $F^*f\in M^1$ follows. Therefore the
  map of base spaces $\underline{F}:Y\to X$ determines a morphism of Riesz
  systems
  \[
    F^1:(X,L^1)\to(Y,M^1) .
  \]

  Finally, by the definition of the extended integral via series
  representations and the fact that $F$ preserves the integral,
  \[
    \begin{aligned}
      J^1(F^*f)
      &=\sum_{n=1}^\infty J(F^*f_n)\\
      &=\sum_{n=1}^\infty I(f_n)\\
      &=I^1(f) .
    \end{aligned}
  \]
  Hence $F^1$ is a morphism of integration systems. Since the Daniell extension
  is a Daniell closed integration system by
  Lemma~\ref{lem:daniell-ext-is-complete}, the morphism $F^1$ is a morphism in
  $\mathbf{DInt}$.
\end{proof}

\begin{lemma}\label{lem:daniell-ext-is-functor-2}
  Assigning to an integration system $(X,L,I)$ on a set $X$ its Daniell
  extension $(X,L^1,I^1)$ determines a functor
  \[
    \mathcal{D}:\mathbf{Int}\to\mathbf{DInt} .
  \]
\end{lemma}
\begin{proof}
  It suffices to check the following.
  \begin{itemize}
    \item For an integration system $(X,L,I)$ on a set $X$, the
      identity morphism $\operatorname{id}_{(X,L,I)}$ satisfies
      $\mathcal{D}(\operatorname{id}_{(X,L,I)})=\operatorname{id}_{(X,L^1,I^1)}$.
    \item For integration systems $(X,L,I),(Y,M,J),(Z,N,K)$ on sets
      $X,Y,Z$ respectively and morphisms of integration systems
      $F:(X,L,I)\to(Y,M,J)$ and $G:(Y,M,J)\to(Z,N,K)$, we have
      $\mathcal{D}(G\circ F)=\mathcal{D}(G)\circ\mathcal{D}(F)$.
  \end{itemize}

  By Lemma~\ref{lem:daniell-ext-is-functor-1}, a morphism of integration systems
  \[
    F:(X,L,I)\to(Y,M,J)
  \]
  determines the morphism
  \[
    F^1:(X,L^1,I^1)\to(Y,M^1,J^1)
  \]
  induced by the same map of base spaces $\underline{F}:Y\to X$. We therefore
  set
  \[
    \mathcal{D}(F):=F^1 .
  \]

  First, the map of base spaces underlying the identity morphism
  $\operatorname{id}_{(X,L,I)}$ is $\operatorname{id}_X$. The morphism
  $\mathcal{D}(\operatorname{id}_{(X,L,I)})$ after the Daniell extension is
  induced by the same map $\operatorname{id}_X$, so
  \[
    \mathcal{D}(\operatorname{id}_{(X,L,I)})
    =\operatorname{id}_{(X,L^1,I^1)} .
  \]

  Next, take morphisms of integration systems
  \[
    F:(X,L,I)\to(Y,M,J),
    \qquad
    G:(Y,M,J)\to(Z,N,K) .
  \]
  The map of base spaces underlying the composite $G\circ F$ is
  \[
    \underline{G\circ F}=\underline{F}\circ\underline{G}:Z\to X .
  \]
  On the other hand, the maps of base spaces underlying
  $\mathcal{D}(F)=F^1$ and $\mathcal{D}(G)=G^1$ are $\underline{F}$ and
  $\underline{G}$ respectively, so the map of base spaces underlying
  $\mathcal{D}(G)\circ\mathcal{D}(F)$ is also
  $\underline{F}\circ\underline{G}$. Hence
  \[
    \mathcal{D}(G\circ F)
    =\mathcal{D}(G)\circ\mathcal{D}(F) .
  \]
  Therefore $\mathcal{D}$ is a functor.
\end{proof}

\subsection{Proof of the main theorem}\label{sec:main-theorem}

\begin{theorem}[main theorem]\label{thm:main-theorem}
  The functor
  \[
    \mathcal{D}:\mathbf{Int}\to\mathbf{DInt}
  \]
  is left adjoint to the inclusion functor
  \[
    \iota_{D}:\mathbf{DInt}\to\mathbf{Int} .
  \]
\end{theorem}
\begin{proof}
  We argue for the functor $\mathcal{D}$ defined in
  Lemma~\ref{lem:daniell-ext-is-functor-2}. Let $(X,L,I)$ be an arbitrary
  integration system and $(X,L^1,I^1)$ its Daniell extension.
  It suffices to show that the morphism of integration systems
  \[
    i_{(X,L,I)}:(X,L,I)\to(X,L^1,I^1)
  \]
  induced by the identity map $\operatorname{id}_X$ is a universal morphism from
  $(X,L,I)$ to $\iota_{D}$.

  Let $(Y,M,J)$ be a Daniell closed integration system on a set $Y$ and
  let
  \[
    F:(X,L,I)\to\iota_{D}(Y,M,J)
  \]
  be a morphism of integration systems. By
  Lemma~\ref{lem:daniell-ext-is-functor-1}, the map of base spaces
  $\underline{F}:Y\to X$ induces a morphism of integration systems
  \[
    F^1:(X,L^1,I^1)\to(Y,M^1,J^1) .
  \]

  Since $(Y,M,J)$ is Daniell closed we have $M^1=M$. Moreover
  $J^1(g)=J(g)$ for every $g\in M$ by
  Lemma~\ref{lem:daniell-ext-contains-source}. Hence
  \[
    (Y,M^1,J^1)=(Y,M,J) .
  \]
  Under this identification we set
  \[
    F^\sharp:=F^1:(X,L^1,I^1)\to(Y,M,J) .
  \]

  The inclusion functor does not change morphisms, and the map of base spaces
  underlying $i_{(X,L,I)}$ is $\operatorname{id}_X$, so the map of base spaces
  underlying the composite $\iota_{D}(F^\sharp)\circ i_{(X,L,I)}$ is
  \[
    \underline{i_{(X,L,I)}}\circ\underline{F^\sharp}
    =\operatorname{id}_X\circ\underline{F}
    =\underline{F} .
  \]
  A morphism of integration systems is precisely its map of base spaces, so
  \[
    F=\iota_{D}(F^\sharp)\circ i_{(X,L,I)} .
  \]

  Next, suppose a morphism
  \[
    H:(X,L^1,I^1)\to(Y,M,J)
  \]
  satisfies
  \[
    F=\iota_{D}(H)\circ i_{(X,L,I)} .
  \]
  Comparing the maps of base spaces on both sides,
  \[
    \underline{F}
    =\underline{i_{(X,L,I)}}\circ\underline{H}
    =\operatorname{id}_X\circ\underline{H}
    =\underline{H} .
  \]
  On the other hand $\underline{F^\sharp}=\underline{F}$, so $H=F^\sharp$.
  Hence $F^\sharp$ is the unique morphism satisfying the factorization above,
  and $i_{(X,L,I)}$ is a universal morphism from $(X,L,I)$ to $\iota_{D}$.

  Since $(X,L,I)$ was arbitrary,
  Theorem~\ref{thm:left-adjunction-existence-from-universality} yields a left
  adjoint of $\iota_{D}$, and its action on objects is the Daniell extension.
  Moreover, for every morphism of integration systems
  \[
    T:(X,L,I)\to(X',L',I')
  \]
  the morphism $\mathcal{D}(T)=T^1$ satisfies
  \[
    \iota_{D}(\mathcal{D}(T))\circ i_{(X,L,I)}
    =i_{(X',L',I')}\circ T ,
  \]
  since the maps of base spaces on both sides are $\underline{T}$.
  By the uniqueness coming from the universal property, the action on morphisms
  of the left adjoint constructed in that theorem agrees with the action of
  $\mathcal{D}$ on morphisms. Therefore $\mathcal{D}$ is left adjoint to
  $\iota_{D}$.
\end{proof}

\begin{remark}[on the content of the main theorem]\label{rem:main-theorem-not-tautology}
  The skeleton of the proof of the main theorem is formal.
  The pair consisting of $T:=\iota_{D}\circ\mathcal{D}$ and the morphisms $i_A$
  is an idempotent pointed endofunctor
  (one has $T(i_A)=i_{TA}=\operatorname{id}$), and it is a standard fact of
  category theory that the full subcategory of fixed points of an idempotent
  pointed endofunctor is reflective.
  Indeed, all that the proof above uses is the following two points:
  \begin{itemize}
    \item idempotence $(L^1)^1=L^1$
      (Lemma~\ref{lem:daniell-ext-is-complete});
    \item that pullback preserves integrability and the value of the integral
      (Lemma~\ref{lem:daniell-ext-is-functor-1}).
  \end{itemize}
  The former is stated as a corollary in \cite{Yamagami} and is known.

  \textbf{Hence, as long as $\mathbf{DInt}$ is defined as ``the category of
  integration systems satisfying $L^1=L$'', the main theorem says no more than
  that the subcategory of fixed points of $\mathcal{D}$ has $\mathcal{D}$ as its
  reflector, and carries no new content.}

  Let us say in advance that the following line of argument does not work: that
  the main theorem has content because
  Theorem~\ref{thm:complete-intrinsic} makes it possible to write Daniell
  closedness ``without mentioning $L^1$''.
  That is invalid. The condition $L^1=L$, when unfolded, also reads ``every
  function admitting a series representation by elements of $L$ belongs to
  $L$'', which uses nothing but $L$, $I$ and pointwise conditions
  (Remark~\ref{rem:not-about-mentioning-L1}).
  Being a reflective subcategory, too, follows formally under either
  description.
  If the rule of closure is preserved under pullback then the subcategory of
  closed objects is reflective, and the proof of
  Theorem~\ref{thm:mint-inclusion-has-left-adjoint} (the reflection into
  $L^{\mathrm{M}}$) has exactly that shape.

  The significance of Theorem~\ref{thm:complete-intrinsic} lies elsewhere.
  It is that \textbf{one closure operation, the Daniell extension, splits into
  two independent closure operations}:
  \[
    L^1=L
    \quad\Longleftrightarrow\quad
    \text{closure under the monotone convergence theorem}
    \quad\textbf{and}\quad
    \text{closure under modification on null sets} .
  \]
  Both conditions on the right occur independently in measure theory.
  The monotone convergence theorem corresponds directly, and closure under
  modification on null sets corresponds to completeness of the measure.
  That neither can be dropped was checked in
  Remark~\ref{rem:two-conditions-independent}.

  Under this decomposition, the main theorem asserts the following.

  \begin{quote}
    The full subcategory of $\mathbf{Int}$ cut out by the two conditions of the
    monotone convergence theorem and modification on null sets is reflective,
    and its reflector is given by the Daniell extension via series
    representations, a construction that can be written down
    \emph{explicitly}.
  \end{quote}

  \noindent
  The content thus moves away from the proof of the adjunction and over to the
  fact that the subcategory cut out by the two conditions coincides with the
  whole set of fixed points of $\mathcal{D}$
  (Theorem~\ref{thm:complete-intrinsic}).
  The two places in this paper with technical content are this one and
  Theorem~\ref{thm:LM-equals-Lsigma}; the proof of the adjunction itself is
  formal.

  This account concerns $\mathbf{Int}$ in
  Section~\ref{chap:functoriality} and says nothing about the category of
  operator morphisms in Section~\ref{chap:operator}.
\end{remark}

\subsection{Why Daniell closedness is imposed on the codomain}\label{sec:why-complete}

The main theorem (Theorem~\ref{thm:main-theorem}) is stated after restricting the
codomain of the factorization to the category $\mathbf{DInt}$ of Daniell closed
integration systems.
In this subsection we show that this restriction cannot be dispensed with.
That is, the functor $\iota_M\circ\mathcal{D}:\mathbf{Int}\to\mathbf{MInt}$
defined by the Daniell extension is \emph{not} left adjoint to the inclusion
functor $\iota:\mathbf{MInt}\to\mathbf{Int}$.
We then show that $\iota$ does nevertheless have a left adjoint, but that its
reflector is \emph{not} the Daniell extension.

\subsubsection{The Daniell extension is not left adjoint to $\iota$}

\begin{remark}\label{rem:D-I-not-left-adjoint}
  By Lemma~\ref{lem:daniell-ext-is-functor-2}, the Daniell extension defines a
  functor
  \[
    \iota_M\circ\mathcal{D}:\mathbf{Int}\to\mathbf{MInt}
  \]
  from $\mathbf{Int}$ to $\mathbf{MInt}$, but it is not left adjoint to the
  inclusion functor
  \[
    \iota:\mathbf{MInt}\to\mathbf{Int} .
  \]

  Indeed, let $X=\{0,1,2\}$, write $\mathbb{R}^X$ in coordinates, and put
  \[
    L:=\{(a,b,b)\mid a,b\in\mathbb{R}\},
    \qquad
    I(a,b,b):=a .
  \]
  Then $L$ is a Riesz space on $X$ and $I$ is an integral. Moreover, if a
  sequence $\{g_n\}$ of functions in $L$ converges pointwise and monotonically to
  a function $g$, then $g$ too takes equal values at the points $1$ and $2$, and
  $I(g_n)=g_n(0)\to g(0)=I(g)$. Hence $(X,L,I)$ satisfies the monotone
  convergence theorem.

  On the other hand, let $u:=(0,1,1)\in L$. For an arbitrary
  $h=(a,b,c)\in\mathbb{R}^X$, putting
  \[
    f_1=(a,a,a),\quad f_n=0\ (n\ge2),
    \qquad
    \varphi_1=|a|(1,1,1),\quad \varphi_n=u\ (n\ge2) ,
  \]
  we get
  \[
    |f_n|\le\varphi_n,
    \qquad
    \sum_{n=1}^{\infty}I(\varphi_n)=|a|<+\infty.
  \]
  Furthermore $\sum_n\varphi_n(x)<+\infty$ holds only at $x=0$, and at that
  point $h(0)=a=\sum_n f_n(0)$. Hence $h\in L^1$ and
  \[
    L^1=\mathbb{R}^X,
    \qquad
    I^1(h)=h(0).
  \]

  Note here that a morphism of integration systems is determined by its map of
  base spaces. Putting $A:=(X,L,I)$, the condition on a map $\tau:X\to X$ for it
  to define a morphism $A\to A$ is
  \[
    \tau(0)=0,
    \qquad
    \tau(1)=\tau(2)\ \text{or}\ \tau(1),\tau(2)\in\{1,2\} ,
  \]
  so that
  \[
    \left|\operatorname{Hom}_{\mathbf{Int}}(A,\iota A)\right|=5.
  \]
  By contrast, a map $\tau:X\to X$ defining a morphism
  $(\iota_M\circ\mathcal{D})(A)\to A$ must satisfy $h\circ\tau\in L$ for every
  $h\in\mathbb{R}^X$. Hence
  \[
    \tau(0)=0,
    \qquad
    \tau(1)=\tau(2),
  \]
  and
  \[
    \left|\operatorname{Hom}_{\mathbf{MInt}}((\iota_M\circ\mathcal{D})(A),A)\right|=3.
  \]
  Therefore there is no bijection between these two Hom-sets, and
  $\iota_M\circ\mathcal{D}$ is not left adjoint to $\iota$. In this example the
  Daniell extension adjoins all functions on the points $1,2$ that $I$ does not
  distinguish, while the codomain $A$ does not contain them. Hence, in order to
  obtain the universal property of a left adjoint from the Daniell extension
  itself, one must restrict the codomain to the Daniell closed integration
  systems.
\end{remark}

\subsubsection{Another reflection possessed by $\iota$}

That the Daniell extension is not left adjoint to $\iota$ does not mean that
$\iota$ has no left adjoint.
In fact $\iota$ does have a left adjoint.
Its reflector is not the Daniell extension, however, but a different
construction, which takes inside $L^1$ the smallest Riesz subspace satisfying
the monotone convergence theorem.

\begin{theorem}\label{thm:mint-inclusion-has-left-adjoint}
  The inclusion functor
  \[
    \iota:\mathbf{MInt}\to\mathbf{Int}
  \]
  has a left adjoint.
\end{theorem}
\begin{proof}
  Let $A=(X,L,I)$ be an object of $\mathbf{Int}$ and let $(X,L^1,I^1)$ be its
  Daniell extension. Consider the family of subspaces of $L^1$
  \[
    \mathfrak{N}_A
    :=\left\{
      N\subset L^1
      \mathrel{}\middle|\mathrel{}
      \begin{array}{l}
        L\subset N,\quad N\text{ is a Riesz space on $X$},\\
        (X,N,I^1|_N)\in\mathbf{MInt}
      \end{array}
    \right\}
  \]
  and put
  \[
    L^{\mathrm{M}}:=\bigcap_{N\in\mathfrak{N}_A}N,
    \qquad
    I^{\mathrm{M}}:=I^1|_{L^{\mathrm{M}}} .
  \]
  By Lemma~\ref{lem:daniell-ext-contains-source} and
  Lemma~\ref{lem:daniell-ext-is-complete} we have $L\subset L^1$ and
  $(X,L^1,I^1)\in\mathbf{MInt}$, so $L^1\in\mathfrak{N}_A$ and
  $\mathfrak{N}_A$ is nonempty.

  $L^{\mathrm{M}}$ is an intersection of Riesz spaces, hence a Riesz space, and
  $L\subset L^{\mathrm{M}}$. Moreover $I^{\mathrm{M}}$ is a restriction of the
  integral $I^1$, so $(X,L^{\mathrm{M}},I^{\mathrm{M}})$ is an integration
  system.

  We show that this integration system satisfies the monotone convergence
  theorem. Suppose $\{f_n\}\subset L^{\mathrm{M}}$ and $f:X\to\mathbb{R}$ satisfy
  \[
    f_n(x)\uparrow f(x)\quad(\forall x\in X),
    \qquad
    \sup_n I^{\mathrm{M}}(f_n)<+\infty .
  \]
  For every $N\in\mathfrak{N}_A$ we have $f_n\in N$, and $(X,N,I^1|_N)$
  satisfies the monotone convergence theorem, so
  \[
    f\in N,
    \qquad
    I^1(f)=\lim_{n\to\infty}I^1(f_n).
  \]
  As $N$ was arbitrary, $f\in L^{\mathrm{M}}$ and
  \[
    I^{\mathrm{M}}(f)=\lim_{n\to\infty}I^{\mathrm{M}}(f_n).
  \]
  Hence
  \[
    A^{\mathrm{M}}:=(X,L^{\mathrm{M}},I^{\mathrm{M}})
  \]
  is an object of $\mathbf{MInt}$. Since $I^1|_L=I$ by
  Lemma~\ref{lem:daniell-ext-contains-source}, the identity map
  $\operatorname{id}_X$ induces a morphism of integration systems
  \[
    i_A:A\to\iota(A^{\mathrm{M}}) .
  \]

  Let $B=(Y,M,J)$ be an object of $\mathbf{MInt}$ and take a morphism of
  integration systems
  \[
    F:A\to\iota(B) .
  \]
  By Lemma~\ref{lem:daniell-ext-is-functor-1}, the same map of base spaces
  $\underline{F}:Y\to X$ induces a morphism of integration systems
  \[
    F^1:(X,L^1,I^1)\to(Y,M^1,J^1) .
  \]
  Put
  \[
    N_F:=\{f\in L^1\mid F^*f\in M\} .
  \]
  Since $F$ is a morphism of integration systems we have $L\subset N_F$, and
  since $M$ is a Riesz space so is $N_F$. Moreover, for every $f\in N_F$, since
  $F^1$ preserves the integral and $I^1,J^1$ extend the original integrals,
  \[
    J(F^*f)=J^1(F^*f)=I^1(f).
  \]

  We show that $(X,N_F,I^1|_{N_F})$ satisfies the monotone convergence theorem.
  Suppose $\{f_n\}\subset N_F$ and $f:X\to\mathbb{R}$ satisfy
  \[
    f_n(x)\uparrow f(x)\quad(\forall x\in X),
    \qquad
    \sup_n I^1(f_n)<+\infty .
  \]
  Since $(X,L^1,I^1)$ satisfies the monotone convergence theorem,
  \[
    f\in L^1,
    \qquad
    I^1(f)=\lim_{n\to\infty}I^1(f_n).
  \]
  On the other hand $F^*f_n\uparrow F^*f$ and
  \[
    \sup_nJ(F^*f_n)=\sup_nI^1(f_n)<+\infty.
  \]
  Since $B$ satisfies the monotone convergence theorem we get $F^*f\in M$, hence
  $f\in N_F$. Therefore $N_F\in\mathfrak{N}_A$, and we obtain
  $L^{\mathrm{M}}\subset N_F$.

  Hence $\underline{F}:Y\to X$ induces a morphism of integration systems
  \[
    F^\sharp:A^{\mathrm{M}}\to B .
  \]
  Since the map of base spaces underlying $i_A$ is $\operatorname{id}_X$,
  \[
    F=\iota(F^\sharp)\circ i_A.
  \]
  Moreover, if a morphism $H:A^{\mathrm{M}}\to B$ satisfies this equality, then
  comparing the maps of base spaces gives $\underline{H}=\underline{F}$. A
  morphism of integration systems is precisely its map of base spaces, so
  $H=F^\sharp$. Hence $i_A$ is a universal morphism from $A$ to $\iota$.

  Since $A$ was arbitrary,
  Theorem~\ref{thm:left-adjunction-existence-from-universality} shows that
  $\iota$ has a left adjoint.
\end{proof}

\subsubsection{When do the two reflections agree?}

We have thus obtained two reflections for an object $A=(X,L,I)$ of
$\mathbf{Int}$.
\begin{itemize}
  \item the reflection $\mathcal{D}A=(X,L^1,I^1)$ into $\mathbf{DInt}$
    (Theorem~\ref{thm:main-theorem});
  \item the reflection $A^{\mathrm{M}}=(X,L^{\mathrm{M}},I^{\mathrm{M}})$ into
    $\mathbf{MInt}$ (Theorem~\ref{thm:mint-inclusion-has-left-adjoint}).
\end{itemize}
In both cases the base space remains $X$, and by construction
\[
  L\subset L^{\mathrm{M}}\subset L^1 .
\]
When do these two agree?

\begin{lemma}\label{lem:two-reflections-same-daniell}
  Let $A$ be an object of $\mathbf{Int}$. Then
  \[
    \mathcal{D}(A^{\mathrm{M}})=\mathcal{D}(A),
  \]
  that is, $(L^{\mathrm{M}})^1=L^1$, and the extended integrals agree on them.
\end{lemma}
\begin{proof}
  Since $L\subset L^{\mathrm{M}}$ and $I^{\mathrm{M}}|_L=I$, a series
  representation of $f\in L^1$ with respect to $(X,L,I)$ is also, as it stands, a
  series representation with respect to
  $(X,L^{\mathrm{M}},I^{\mathrm{M}})$.
  Hence $L^1\subset(L^{\mathrm{M}})^1$.

  Conversely, since $L^{\mathrm{M}}\subset L^1$ and
  $I^1|_{L^{\mathrm{M}}}=I^{\mathrm{M}}$, the same reasoning gives
  $(L^{\mathrm{M}})^1\subset(L^1)^1$, and
  $(L^1)^1=L^1$ by Lemma~\ref{lem:daniell-ext-is-complete}.
  Together these give $(L^{\mathrm{M}})^1=L^1$.
  That the extended integrals agree is immediate from the definition via series
  representations.
\end{proof}

\begin{proposition}\label{prop:two-reflections-coincide}
  Let $A$ be an object of $\mathbf{Int}$. The following are equivalent.
  \begin{enumerate}
    \item[(i)] $L^{\mathrm{M}}=L^1$, that is, $A^{\mathrm{M}}=\mathcal{D}A$.
    \item[(ii)] $A^{\mathrm{M}}$ is Daniell closed.
    \item[(iii)] The morphism $A^{\mathrm{M}}\to\mathcal{D}A$ induced by the
      identity map $\operatorname{id}_X$ is an isomorphism.
    \item[(iv)] There is no proper Riesz subspace of $L^1$ containing $L$ that
      defines an integration system satisfying the monotone convergence theorem.
  \end{enumerate}
\end{proposition}
\begin{proof}
  (i)$\Rightarrow$(ii) holds because $\mathcal{D}A$ is Daniell closed
  (Lemma~\ref{lem:daniell-ext-is-complete}).
  (ii)$\Rightarrow$(i) holds because, if $A^{\mathrm{M}}$ is Daniell closed, then
  $(L^{\mathrm{M}})^1=L^{\mathrm{M}}$, and the left-hand side equals $L^1$ by
  Lemma~\ref{lem:two-reflections-same-daniell}.
  (i)$\Leftrightarrow$(iii) is immediate, since the objects share the base space
  and a morphism is precisely its map of base spaces.
  (i)$\Leftrightarrow$(iv) holds because $L^{\mathrm{M}}$ is the intersection of
  $\mathfrak{N}_A$, that is, the smallest of the subspaces satisfying the
  conditions.
\end{proof}

\begin{example}
  In the example of Remark~\ref{rem:D-I-not-left-adjoint}, $A=(X,L,I)$ was
  already an object of $\mathbf{MInt}$, so $L^{\mathrm{M}}=L$, while
  $L^1=\mathbb{R}^X$.
  Hence $L^{\mathrm{M}}\subsetneq L^1$ and the two reflections do not agree.

  It is easy to see what is going on in this example.
  The integral $I$ does not distinguish the point $1$ from the point $2$, and the
  elements of $L$ take the same value at these two points.
  The Daniell extension adjoins all functions on these two points, but they are
  not needed in order to satisfy the monotone convergence theorem.
\end{example}

The example above shows that failure of $L$ to separate the points of the base
space can produce a gap between the two reflections.
But separation is not a sufficient condition.
Below we determine $L^{\mathrm{M}}$ in terms of the $\sigma$-algebra generated by
$L$ and derive that fact as a corollary.

\begin{definition}\label{def:sigma-L}
  Let $(X,L,I)$ be an integration system. We call the $\sigma$-algebra on $X$
  \[
    \sigma(L)
    :=\sigma\bigl(\{f^{-1}(B)\mid f\in L,\ B\in\mathcal{B}(\mathbb{R})\}\bigr)
  \]
  the \textbf{$\sigma$-algebra generated by $L$}.
  That is, $\sigma(L)$ is the smallest $\sigma$-algebra making every element of
  $L$ measurable.
  We also put
  \[
    L^{\sigma}:=\{f\in L^1\mid f\text{ is }\sigma(L)\text{-measurable}\},
    \qquad
    I^{\sigma}:=I^1|_{L^{\sigma}} .
  \]
\end{definition}

\begin{proposition}\label{prop:LM-in-Lsigma}
  For every integration system $A=(X,L,I)$,
  \[
    L\subset L^{\mathrm{M}}\subset L^{\sigma}\subset L^1 .
  \]
\end{proposition}
\begin{proof}
  It suffices to show $(X,L^{\sigma},I^{\sigma})\in\mathfrak{N}_A$.

  Every element of $L$ is $\sigma(L)$-measurable (by definition) and
  $L\subset L^1$, so $L\subset L^{\sigma}$.
  The set of all $\sigma(L)$-measurable real-valued functions is a Riesz space on
  $X$, and $L^1$ is a Riesz space as well
  (Theorem~\ref{thm:l1-is-riesz-space}), so their intersection $L^{\sigma}$ is a
  Riesz space on $X$.
  $I^{\sigma}$ is the restriction of the integral $I^1$ to a Riesz subspace, so
  $(X,L^{\sigma},I^{\sigma})$ is an integration system.

  We check the monotone convergence theorem. Suppose
  $\{f_n\}\subset L^{\sigma}$ and $f:X\to\mathbb{R}$ satisfy
  $f_n(x)\uparrow f(x)$ for all $x\in X$ and
  $\sup_nI^{\sigma}(f_n)<+\infty$.
  Since $(X,L^1,I^1)$ satisfies the monotone convergence theorem
  (Lemma~\ref{lem:daniell-ext-is-complete}), we get $f\in L^1$ and
  $I^1(f)=\lim_nI^1(f_n)$.
  Moreover $f$ is a pointwise limit of $\sigma(L)$-measurable functions, hence
  $\sigma(L)$-measurable.
  Therefore $f\in L^{\sigma}$ and $I^{\sigma}(f)=\lim_nI^{\sigma}(f_n)$.

  Hence $L^{\sigma}\in\mathfrak{N}_A$, and since $L^{\mathrm{M}}$ is the
  intersection of $\mathfrak{N}_A$ we obtain
  $L^{\mathrm{M}}\subset L^{\sigma}$.
\end{proof}

\begin{definition}\label{def:supported-on-null}
  Let $(X,L,I)$ be an integration system. We put
  \[
    \mathcal{Z}_I:=\bigl\{h:X\to\mathbb{R}\;\bigm|\;
      \text{there is an $I$-null set $Z$ with $h|_{X\setminus Z}=0$}\bigr\}
  \]
  and call the elements of $\mathcal{Z}_I$ \textbf{functions supported on a null
  set}.
\end{definition}

\noindent
By Lemma~\ref{lem:null-sets-countable-union}, $\mathcal{Z}_I$ is a Riesz space on
$X$.
Indeed, if $h_1$ and $h_2$ vanish outside $Z_1$ and $Z_2$ respectively, then
both $h_1+h_2$ and $|h_1|$ vanish outside $Z_1\cup Z_2$.

\begin{proposition}\label{prop:l1-decomposition}
  For every integration system $A=(X,L,I)$,
  \[
    L^1=L^{\sigma}+\mathcal{Z}_I .
  \]
  Moreover $h\in\mathcal{Z}_I$ implies $I^1(|h|)=0$, so for a decomposition
  $f=g+h$ with $g\in L^{\sigma}$ and $h\in\mathcal{Z}_I$ we have
  $I^1(f)=I^1(g)$.
\end{proposition}
\begin{proof}
  First we check $\mathcal{Z}_I\subset L^1$.
  Let $h\in\mathcal{Z}_I$ and take an $I$-null set $Z$ and a corresponding
  sequence $\{g_n\}$.
  Putting $f_n:=0$ and $\varphi_n:=g_n$ we have $|f_n|\le\varphi_n$ and
  $\sum_nI(\varphi_n)<+\infty$, and any $x$ with $\sum_n\varphi_n(x)<+\infty$
  does not belong to $Z$, so $h(x)=0=\sum_nf_n(x)$.
  Hence $\{f_n\},\{\varphi_n\}$ is a series representation of $h$ and
  $h\in L^1$.
  Moreover $\{g_n\}\subset L\subset L^1$ and $I^1|_L=I$
  (Lemma~\ref{lem:daniell-ext-contains-source}), so $Z$ is a null set for the
  integration system $(X,L^1,I^1)$ as well.
  Since $|h|$ vanishes outside $Z$, applying
  Lemma~\ref{lem:vanish-off-null-set} to $(X,L^1,I^1)$ gives $I^1(|h|)=0$.

  $L^{\sigma}\subset L^1$ holds by definition. Since $L^1$ is a Riesz space
  (Theorem~\ref{thm:l1-is-riesz-space}) we get
  $L^{\sigma}+\mathcal{Z}_I\subset L^1$.

  We show the converse. Let $f\in L^1$, choose a series representation
  $f\overset{\varphi_n}{\simeq}\sum f_n$, and put
  \[
    \Phi_N:=\sum_{n\le N}\varphi_n,
    \qquad
    Z:=\Bigl\{x\in X\;\Bigm|\;\sum_n\varphi_n(x)=+\infty\Bigr\} .
  \]
  Then $\Phi_N\in L$ is $\sigma(L)$-measurable and
  \[
    Z=\bigcap_{M\ge1}\bigcup_{N\ge1}\{\Phi_N>M\} ,
  \]
  so $Z\in\sigma(L)$.
  Moreover $|f_n|\le\varphi_n$ gives $\varphi_n\ge0$, and
  $\sum_nI(\varphi_n)<+\infty$, so $\{\varphi_n\}$ shows that $Z$ is an $I$-null
  set.

  Define a function $g:X\to\mathbb{R}$ by
  \[
    g(x):=
    \begin{cases}
      \sum_nf_n(x) & (x\notin Z)\\
      0 & (x\in Z).
    \end{cases}
  \]
  For $x\notin Z$ the series $\sum_nf_n(x)$ converges absolutely
  (Lemma~\ref{lem:series-exp-is-abs-conv}), so $g$ is real-valued.
  Moreover $g$ is a pointwise limit of the $\sigma(L)$-measurable functions
  $1_{X\setminus Z}\sum_{n\le N}f_n$, hence $\sigma(L)$-measurable.

  The pair $\{f_n\},\{\varphi_n\}$ is also a series representation of $g$.
  Indeed, conditions (1) and (2) are unchanged, and condition (3) is precisely
  the statement that $g(x)=\sum_nf_n(x)$ when
  $\sum_n\varphi_n(x)<+\infty$, that is, when $x\notin Z$.
  Hence $g\in L^1$, and so $g\in L^{\sigma}$.

  Finally, if $x\notin Z$ then $f(x)=\sum_nf_n(x)=g(x)$, so $f-g$ vanishes
  outside $Z$. Therefore $f-g\in\mathcal{Z}_I$ and
  $f=g+(f-g)\in L^{\sigma}+\mathcal{Z}_I$.
\end{proof}

\noindent
The reverse inclusion holds under the Stone condition
(Definition~\ref{def:stone-condition}).
The Stone condition is weaker than requiring that $L$ contain the constant
function $1$, and it includes integration systems of infinite total mass
(Example~\ref{ex:step-functions-R}).

\begin{theorem}\label{thm:LM-equals-Lsigma}
  Let $A=(X,L,I)$ be an integration system satisfying the Stone condition. Then
  \[
    L^{\mathrm{M}}=L^{\sigma}.
  \]
\end{theorem}
\begin{proof}
  Since $L^{\mathrm{M}}\subset L^{\sigma}$ by
  Proposition~\ref{prop:LM-in-Lsigma}, it suffices to show that
  $L^{\sigma}\subset N$ holds for every $N\in\mathfrak{N}_A$.
  Fix $N\in\mathfrak{N}_A$ in what follows.

  \textbf{Step 1: if $g\in L$, $g\ge0$ and $a>0$ then
  $1_{\{g>a\}}\in N\cap L^1$ and $I^1(1_{\{g>a\}})\le I(g)/a$.}
  By the Stone condition $g\wedge a=a\bigl((g/a)\land1\bigr)\in L$, so
  \[
    h_n:=\bigl(n(g-g\land a)\bigr)\land1\in L .
  \]
  Putting $t:=g(x)$ at a point $x$: if $t\le a$ then $t-(t\land a)=0$ and so
  $h_n(x)=0$, and if $t>a$ then $n(t-a)\uparrow+\infty$ and so
  $h_n(x)\uparrow1$.
  Hence $h_n\uparrow1_{\{g>a\}}$ pointwise and monotonically nondecreasing.
  Moreover $0\le h_n\le1$ and $h_n$ vanishes outside $\{g>a\}$, while on
  $\{g>a\}$ we have $g/a>1$, so $h_n\le g/a$ and
  \[
    I(h_n)\le I(g)/a<+\infty.
  \]
  Since $N$ contains $L$ and satisfies the monotone convergence theorem,
  $1_{\{g>a\}}\in N$.
  As $N\subset L^1$ and $I^1$ extends $I$, the estimate on the integral follows
  as well.

  \textbf{Step 2: if $G\in\sigma(L)$ satisfies $1_G\in N$, then $1_E\in N$ for
  every $E\in\sigma(L)$ with $E\subset G$.}
  Put
  \[
    \mathcal{A}_G:=\{E\in\sigma(L)\mid E\subset G,\ 1_E\in N\} .
  \]
  We check that this is a $\sigma$-algebra on $G$.
  We have $G\in\mathcal{A}_G$. If $E\in\mathcal{A}_G$ then
  $1_{G\setminus E}=1_G-1_E\in N$.
  If $E_1,E_2\in\mathcal{A}_G$ then
  $1_{E_1\cup E_2}=1_{E_1}\lor1_{E_2}\in N$.
  Finally, let $E_k\in\mathcal{A}_G$, put $E:=\bigcup_kE_k$ and
  $E'_m:=E_1\cup\dots\cup E_m$; then $E'_m\in\mathcal{A}_G$ and
  $1_{E'_m}\uparrow1_E$, and $E'_m\subset G$ gives
  \[
    I^1(1_{E'_m})\le I^1(1_G)<+\infty .
  \]
  By the monotone convergence theorem for $N$ we get $1_E\in N$.

  Next we show that the trace $\{E\cap G\mid E\in\sigma(L)\}$ of $\sigma(L)$ on
  $G$ is contained in $\mathcal{A}_G$.
  The trace $\sigma$-algebra is generated by the traces of a generating family of
  $\sigma(L)$, so it suffices to check $\{f>a\}\cap G\in\mathcal{A}_G$ for
  $f\in L$ and $a\in\mathbb{R}$.
  \begin{itemize}
    \item Case $a>0$. We have $\{f>a\}=\{f^+>a\}$ with $f^+\in L$ and
      $f^+\ge0$, so Step 1 gives $1_{\{f>a\}}\in N$. Hence
      $1_{\{f>a\}\cap G}=1_{\{f>a\}}\land1_G\in N$.
    \item Case $a=0$. Here $\{f>0\}=\bigcup_{k\ge1}\{f>1/k\}$ is a countable
      union of sets from the previous case.
    \item Case $a<0$. Then $-a>0$, and for $k>1/(-a)$ we have $-a-1/k>0$, so
      \[
        \{f\le a\}=\{-f\ge-a\}=\bigcap_{k>1/(-a)}\{-f>-a-\tfrac1k\}
      \]
      is a countable intersection of sets from the case $a>0$. The claim follows
      from $\{f>a\}\cap G=G\setminus(\{f\le a\}\cap G)$.
  \end{itemize}
  Since $\mathcal{A}_G$ is a $\sigma$-algebra it contains the whole trace.
  An $E\in\sigma(L)$ with $E\subset G$ belongs to the trace as $E=E\cap G$, so
  $1_E\in N$.

  \textbf{Step 3: $L^{\sigma}\subset N$.}
  Both $L^{\sigma}$ and $N$ are Riesz spaces and $f$ decomposes as
  $f=f^+-f^-$, so it suffices to treat the case $f\in L^{\sigma}$ with $f\ge0$.
  Choose a series representation $f\overset{\varphi_n}{\simeq}\sum f_n$.

  At a point $x$ where $\varphi_n(x)=0$ for all $n$ we have
  $\sum_n\varphi_n(x)=0<+\infty$, so $f(x)=\sum_nf_n(x)$, and
  $|f_n(x)|\le\varphi_n(x)=0$ gives $f(x)=0$. Hence
  \[
    \{f\ne0\}\subset\bigcup_{n,k\ge1}G_{n,k},
    \qquad
    G_{n,k}:=\{\varphi_n>1/k\} .
  \]
  Since $\varphi_n\in L$ and $\varphi_n\ge0$, Step 1 gives
  $1_{G_{n,k}}\in N\cap L^1$, and $G_{n,k}\in\sigma(L)$.
  Arranging the family $\{G_{n,k}\}_{n,k}$ into a single sequence
  $G_1,G_2,\dots$ and putting $H_m:=G_1\cup\dots\cup G_m$, we get
  $1_{H_m}=1_{G_1}\lor\dots\lor1_{G_m}\in N$ and $H_m\in\sigma(L)$.

  Take $c>0$. Since $f$ is $\sigma(L)$-measurable we have
  $\{f\ge c\}\in\sigma(L)$, and $\{f\ge c\}\cap H_m$ belongs to $\sigma(L)$ and
  is contained in $H_m$, so Step 2 gives
  $1_{\{f\ge c\}\cap H_m}\in N$.
  From $c>0$ and $f\ge0$ we get
  $\{f\ge c\}\subset\{f\ne0\}\subset\bigcup_mH_m$, so
  $\{f\ge c\}\cap H_m\uparrow\{f\ge c\}$, that is,
  $1_{\{f\ge c\}\cap H_m}\uparrow1_{\{f\ge c\}}$.
  Moreover $f/c\ge1$ on $\{f\ge c\}$ and $f/c\ge0$ off it, so
  $1_{\{f\ge c\}}\le f/c$ and
  \[
    I^1(1_{\{f\ge c\}\cap H_m})\le I^1(f)/c<+\infty.
  \]
  By the monotone convergence theorem for $N$ we get $1_{\{f\ge c\}}\in N$.

  Finally, putting
  \[
    f_m:=\sum_{k=1}^{m2^m}2^{-m}1_{\{f\ge k2^{-m}\}}
    \qquad\Bigl(=\bigl(2^{-m}\lfloor2^mf\rfloor\bigr)\land m\Bigr) ,
  \]
  the function $f_m$ belongs to $N$ as a finite linear combination of elements of
  $N$.
  It is nondecreasing in $m$ and converges pointwise to $f$, and
  $0\le f_m\le f$ gives $I^1(f_m)\le I^1(f)<+\infty$.
  By the monotone convergence theorem for $N$ we get $f\in N$.

  Since $N\in\mathfrak{N}_A$ was arbitrary we have
  $L^{\sigma}\subset L^{\mathrm{M}}$, and together with
  Proposition~\ref{prop:LM-in-Lsigma} we obtain
  $L^{\mathrm{M}}=L^{\sigma}$.
\end{proof}

\begin{corollary}\label{cor:LM-equals-Lsigma-one}
  If the constant function $1$ belongs to $L$ then
  $L^{\mathrm{M}}=L^{\sigma}$.
\end{corollary}
\begin{proof}
  Since $1\in L$ and $L$ is a Riesz space, $f\land1\in L$ for every $f\in L$,
  that is, the Stone condition holds.
\end{proof}

\noindent
When $1\in L$ one may take $G:=X$, so the localization to $H_m$ in Step 2 of the
proof above is unnecessary.

\medskip
\noindent
\textbf{Relation to classical Daniell--Stone theory.}
The content of Theorem~\ref{thm:LM-equals-Lsigma} lies in the neighbourhood of
the classical theory.
The smallest family containing a Riesz space satisfying the Stone condition (a
Stone vector lattice) and closed under monotone limits is a standard object,
called the $\sigma$-ring of real functions generated by $L$, and that it
corresponds to $\sigma(L)$-measurability is classically known.
\S4 of Chapter 16 of \cite{Royden}, ``Measurability and measure'', is this part.

What distinguishes the present theorem from the classical assertion is that the
monotone limits are taken inside $L^1$, and only under the boundedness condition
$\sup_nI(f_n)<+\infty$ on the integrals
(Definition~\ref{def:satisfies-mct}).
This is why the generated monotone family cannot be used as it stands, and it is
also the source of the truncation in Step 1 and the localization to $H_m$ in the
proof.
The conclusion obtained is the identification of \textbf{that restricted closure
with the second reflection $L^{\mathrm{M}}$}
(Theorem~\ref{thm:mint-inclusion-has-left-adjoint}).

\begin{corollary}\label{cor:l1-decomposition-LM}
  If $A=(X,L,I)$ satisfies the Stone condition then
  \[
    L^1=L^{\mathrm{M}}+\mathcal{Z}_I .
  \]
  That is, \textbf{the difference between $L^1$ and the second reflection
  $L^{\mathrm{M}}$ amounts exactly to the functions supported on null sets}.
\end{corollary}
\begin{proof}
  By Theorem~\ref{thm:LM-equals-Lsigma} we have
  $L^{\mathrm{M}}=L^{\sigma}$, so the claim follows from
  Proposition~\ref{prop:l1-decomposition}.
\end{proof}

\begin{remark}[two-stage factorization of the Daniell extension]\label{rem:D-factors}
  Under the Stone condition, the results so far split the Daniell extension into
  the following two stages.
  \[
    L
    \quad\subset\quad
    \underbrace{L^{\mathrm{M}}=L^{\sigma}}_{\text{first stage}}
    \quad\subset\quad
    \underbrace{L^1=L^{\mathrm{M}}+\mathcal{Z}_I}_{\text{second stage}}
  \]
  The first stage is the operation of closing $L$ under monotone limits, and it
  is obtained as the reflection into $\mathbf{MInt}$
  (Theorem~\ref{thm:mint-inclusion-has-left-adjoint}).
  Its destination is the set of all $\sigma(L)$-measurable integrable functions
  (Theorem~\ref{thm:LM-equals-Lsigma}).
  The second stage is the operation of adjoining the functions supported on null
  sets (Corollary~\ref{cor:l1-decomposition-LM}).

  Read in the language of measures, the first stage is the operation of taking
  the integrable functions measurable with respect to the $\sigma$-algebra
  generated by $L$, and the second is the operation of completing the measure
  (passing through Theorem~\ref{thm:daniell-stone}).
  That is,
  \begin{quote}
    \textbf{the Daniell extension is the composite of Daniell--Stone style
    measurability and completion of the measure.}
  \end{quote}

  \noindent
  Both stages are classical operations.
  What this paper supplies is the identification of this factorization
  (Theorem~\ref{thm:LM-equals-Lsigma} and
  Proposition~\ref{prop:l1-decomposition}) and the categorical placement of the
  composite as a reflector.
  As stated in Subsection~\ref{sec:prior-work},
  \textbf{the position of this paper, in a phrase, is ``the Daniell counterpart
  of the completion of a measure space''.}
\end{remark}

\begin{corollary}\label{cor:when-coincide-sigma}
  Let $A=(X,L,I)$ be an integration system satisfying the Stone condition.
  Then the equivalent conditions of
  Proposition~\ref{prop:two-reflections-coincide} are equivalent to the
  following.
  \begin{enumerate}
    \item[(v)] Every element of $L^1$ is $\sigma(L)$-measurable.
  \end{enumerate}
\end{corollary}
\begin{proof}
  By Theorem~\ref{thm:LM-equals-Lsigma} we have
  $L^{\mathrm{M}}=L^{\sigma}$, and $L^{\sigma}=L^1$ is a restatement of (v).
\end{proof}

\noindent
This is a characterization in terms of properties of $(X,L,I)$ itself, and makes
no reference to the construction of $L^{\mathrm{M}}$.

\begin{corollary}\label{cor:lebesgue-reflections-differ}
  Let $A$ be the integration system of Example~\ref{ex:step-functions}, that is,
  the integration system consisting of the step functions on $[0,1)$ and the
  Riemann integral. Then
  \[
    L^{\mathrm{M}}=\{f:[0,1)\to\mathbb{R}\text{ Borel measurable and Lebesgue integrable}\}
    \subsetneq L^1
    =\{f:[0,1)\to\mathbb{R}\text{ Lebesgue integrable}\} ,
  \]
  and the two reflections do not agree.
\end{corollary}
\begin{proof}
  Truncating a step function at $1$ is again a step function, so this
  integration system satisfies the Stone condition
  (Example~\ref{ex:stone-holds}).
  Moreover $\sigma(L)$ is the $\sigma$-algebra generated by the indicator
  functions of the half-open intervals $[a,b)\subset[0,1)$, that is, the Borel
  $\sigma$-algebra of $[0,1)$.
  The equalities follow from the fact that $L^1$ is the set of all Lebesgue
  integrable functions (Example~\ref{ex:lebesgue}) together with
  Theorem~\ref{thm:LM-equals-Lsigma}.

  We check that the inclusion is proper. The Cantor set $C\subset[0,1]$ is a
  Lebesgue null set and has the cardinality of the continuum $\mathfrak{c}$.
  Hence $C\setminus\{1\}$ has $2^{\mathfrak{c}}$ subsets, while $[0,1)$ has only
  $\mathfrak{c}$ Borel sets.
  There is therefore a non-Borel set $E\subset C\setminus\{1\}$.
  As a subset of a Lebesgue null set, $E$ is Lebesgue measurable, and $1_E$ is
  Lebesgue integrable with integral $0$, so $1_E\in L^1$.
  On the other hand $E\notin\mathcal{B}([0,1))$, so $1_E$ is not Borel
  measurable and $1_E\notin L^{\mathrm{M}}$.
\end{proof}

\begin{remark}\label{rem:separation-not-sufficient}
  For the integration system of
  Corollary~\ref{cor:lebesgue-reflections-differ}, $L$ separates the points of
  $[0,1)$ (for $x\ne y$ one can take a half-open interval containing $x$ and not
  $y$).
  It satisfies the Stone condition as well
  (Example~\ref{ex:stone-holds}).
  Nevertheless the two reflections do not agree.
  Hence \textbf{neither separation of points nor the Stone condition is a
  sufficient condition for the two reflections to agree}.

  What produces the gap can be identified in general, not only in this example.
  By Proposition~\ref{prop:l1-decomposition} we have
  $L^1=L^{\sigma}+\mathcal{Z}_I$ for every integration system, and under the
  Stone condition $L^1=L^{\mathrm{M}}+\mathcal{Z}_I$ by
  Corollary~\ref{cor:l1-decomposition-LM}.
  That is, \textbf{the gap between the two reflections amounts exactly to the
  functions supported on null sets}.
  In this example it appears in the form that $L^{\mathrm{M}}$ contains only
  Borel measurable functions whereas $L^1$ contains functions behaving
  arbitrarily on null sets.
  That is, it is the difference between Borel and Lebesgue, the amount by which
  the measure is completed.

  It is the fact that the Daniell extension imposes condition (3) of a series
  representation ``only at the points where the sum of the dominating functions
  converges'' (see also Remark~\ref{rem:why-half-open}) that performs this
  completion automatically.
\end{remark}

\begin{example}
  Let us verify Theorem~\ref{thm:LM-equals-Lsigma} in the example of
  Remark~\ref{rem:D-I-not-left-adjoint} ($X=\{0,1,2\}$, $L=\{(a,b,b)\}$,
  $I(a,b,b)=a$).
  Here $1=(1,1,1)\in L$.
  The elements of $L$ take the same value at the points $1$ and $2$, so
  \[
    \sigma(L)=\{\emptyset,\{0\},\{1,2\},X\} ,
  \]
  and the $\sigma(L)$-measurable functions are precisely the functions taking the
  same value at the points $1$ and $2$, that is, the elements of $L$.
  Hence $L^{\sigma}=L^1\cap L=L$, and $L^{\mathrm{M}}=L$ as the theorem asserts.
  This agrees with the value computed directly in the example above.
\end{example}

\begin{example}\label{ex:step-functions-R}
  Theorem~\ref{thm:LM-equals-Lsigma} applies to integration systems of infinite
  total mass as well.
  Let $X=\mathbb{R}$, let $L$ be the set of all step functions with bounded
  support, and let $I$ be their integral.
  Then $1_{\mathbb{R}}\notin L$, but truncating a step function at $1$ is again
  a step function, so the Stone condition is satisfied.
  Moreover $\sigma(L)$ is the $\sigma$-algebra generated by the bounded
  half-open intervals, that is, the Borel $\sigma$-algebra of $\mathbb{R}$.

  Hence, by Theorem~\ref{thm:LM-equals-Lsigma}, $L^{\mathrm{M}}$ is the set of
  all Borel measurable elements of $L^1$.
  The two reflections do not agree for this integration system either.
  Indeed, the Cantor set $C\subset[0,1]$ is an $I$-null set (take step functions
  $g_n$ covering $C$ by finitely many intervals of total length at most
  $2^{-n}$), and by the cardinality argument (the same as in the proof of
  Corollary~\ref{cor:lebesgue-reflections-differ}) there is a non-Borel set
  $E\subset C$.
  Since $1_E$ vanishes outside $C$ we have
  $1_E\in\mathcal{Z}_I\subset L^1$
  (Proposition~\ref{prop:l1-decomposition}), but $1_E$ is not Borel measurable,
  so it does not belong to $L^{\mathrm{M}}$.
\end{example}

\begin{problem}\label{prob:when-reflections-coincide}
  Is the hypothesis ``the Stone condition'' in
  Theorem~\ref{thm:LM-equals-Lsigma} necessary?
  That is, is there an integration system that does not satisfy the Stone
  condition and has $L^{\mathrm{M}}\ne L^{\sigma}$?
\end{problem}

\noindent
That the Stone condition is not necessary is shown by
Example~\ref{ex:stone-fails}.
That example ($X=\{1,2\}$, $L=\{(a,2a)\}$, $I(a,2a)=a$) does not satisfy the
Stone condition, but $L^1=L$ as seen there.
Moreover $\sigma(L)=2^X$ from the preimages under $f=(1,2)$, so
$L^{\sigma}=L^1=L$ and therefore $L^{\mathrm{M}}=L^{\sigma}$ holds.
The question is whether the equality can fail once the Stone condition is
dropped.

The Stone condition was used only in Step 1 of the proof, that is, at the point
where the indicator functions of the upper sets of elements of $L$ are extracted
inside $N$.
Without that operation there is no guarantee that the indicator functions of
elements of $\sigma(L)$ belong to $N$, and the footing of Steps 2 and 3 is lost.

\subsection{The ae version of the main theorem}\label{sec:main-theorem-ae}

An adjunction of the same shape as the main theorem holds in the category
obtained by identifying morphisms up to null functions
(Definition~\ref{def:int-ae}).

\begin{lemma}\label{lem:daniell-extension-preserves-ae}
  Let $(X,L,I)$ and $(Y,M,J)$ be integration systems and let
  \[
    F,G:(X,L,I)\to(Y,M,J)
  \]
  be morphisms of integration systems. Then
  \[
    F=_{\mathrm{ae}}G
    \quad\Longrightarrow\quad
    F^1=_{\mathrm{ae}}G^1 ,
  \]
  where $F^1,G^1$ are the morphisms
  \[
    (X,L^1,I^1)\to(Y,M^1,J^1)
  \]
  induced by the Daniell extension.
\end{lemma}
\begin{proof}
  Take $f\in L^1$ and take the sequence of functions $\{s_N\}\subset L$ of
  Lemma~\ref{lem:l1-approx-by-source}.
  Since $F^1,G^1$ are morphisms of Riesz systems and of integration systems and
  $J^1$ extends $J$,
  \[
    \begin{aligned}
      &J^1\bigl(|(F^1)^*f-(G^1)^*f|\bigr)\\
      &\quad\le
      J^1\bigl(|(F^1)^*(f-s_N)|\bigr)
      +J^1\bigl(|F^*s_N-G^*s_N|\bigr)
      +J^1\bigl(|(G^1)^*(s_N-f)|\bigr)\\
      &\quad=
      I^1(|f-s_N|)
      +J\bigl(|F^*s_N-G^*s_N|\bigr)
      +I^1(|s_N-f|)\\
      &\quad=2I^1(|f-s_N|).
    \end{aligned}
  \]
  The last equality used $F=_{\mathrm{ae}}G$ and $s_N\in L$.
  Letting $N\to\infty$, the left-hand side is $0$. As $f\in L^1$ was arbitrary,
  $F^1=_{\mathrm{ae}}G^1$.
\end{proof}

\begin{corollary}\label{cor:daniell-extension-ae-functor}
  The assignment sending an object to its Daniell extension and a morphism to
  \[
    [F]\longmapsto[F^1]
  \]
  defines a functor
  \[
    \mathcal{D}_{\mathrm{ae}}:
    \mathbf{Int}_{\mathrm{ae}}
    \to
    \mathbf{DInt}_{\mathrm{ae}} .
  \]
\end{corollary}
\begin{proof}
  By Lemma~\ref{lem:daniell-extension-preserves-ae} the assignment on morphisms
  does not depend on the choice of representative.
  Moreover, that the Daniell extension is a Daniell closed integration system
  and that identities and composition of morphisms are preserved were shown up
  to Subsection~\ref{sec:main-theorem}.
\end{proof}

\begin{definition}
  We denote the inclusion functor of the full subcategory by
  \[
    \iota_{D,\mathrm{ae}}:
    \mathbf{DInt}_{\mathrm{ae}}
    \to
    \mathbf{Int}_{\mathrm{ae}} .
  \]
\end{definition}

\begin{theorem}[ae version of the main theorem]\label{thm:main-ae}
  The functor
  \[
    \mathcal{D}_{\mathrm{ae}}:
    \mathbf{Int}_{\mathrm{ae}}
    \to
    \mathbf{DInt}_{\mathrm{ae}}
  \]
  defined by the Daniell extension is left adjoint to the inclusion functor
  $\iota_{D,\mathrm{ae}}$.
\end{theorem}
\begin{proof}
  Let $A=(X,L,I)$ be an integration system and write its Daniell extension as
  $\mathcal{D}_{\mathrm{ae}}A=(X,L^1,I^1)$.
  We show that the ae-equivalence class
  \[
    [i_A]:A\to\iota_{D,\mathrm{ae}}\mathcal{D}_{\mathrm{ae}}A
  \]
  of the morphism of integration systems
  \[
    i_A:A\to\iota_{D}\mathcal{D}A
  \]
  constructed in Subsection~\ref{sec:main-theorem} is a universal morphism from
  $A$ to $\iota_{D,\mathrm{ae}}$.

  Let $B=(Y,M,J)$ be a Daniell closed integration system and let
  \[
    [F]:A\to\iota_{D,\mathrm{ae}}B
  \]
  be a morphism. Choose a representative $F$. By the proof of
  Theorem~\ref{thm:main-theorem} and $M^1=M$, there is a morphism
  \[
    F^1:\mathcal{D}_{\mathrm{ae}}A\to B
  \]
  induced by the Daniell extension such that, at the level of representatives,
  \[
    F=\iota_{D}(F^1)\circ i_A .
  \]
  Hence
  \[
    [F]=\iota_{D,\mathrm{ae}}([F^1])\circ[i_A].
  \]

  Next, suppose a morphism $[H]:\mathcal{D}_{\mathrm{ae}}A\to B$ satisfies
  \[
    [F]=\iota_{D,\mathrm{ae}}([H])\circ[i_A] .
  \]
  This means that, for representatives,
  \[
    F=_{\mathrm{ae}}H\circ i_A .
  \]
  Applying Lemma~\ref{lem:daniell-extension-preserves-ae} to these two
  morphisms,
  \[
    F^1=_{\mathrm{ae}}(H\circ i_A)^1 .
  \]
  The Daniell extension does not change the map of base spaces underlying a
  morphism, and $B$ is Daniell closed, so $(H\circ i_A)^1=H$. Hence
  \[
    [F^1]=[H].
  \]
  Therefore the ae-equivalence class of the factorizing morphism is unique and
  $[i_A]$ is a universal morphism.

  Since $A$ was arbitrary, $\iota_{D,\mathrm{ae}}$ has a left adjoint.
  Moreover, for a morphism of integration systems $T:A\to A'$ the class $[T^1]$
  satisfies
  \[
    \iota_{D,\mathrm{ae}}([T^1])\circ[i_A]
    =[i_{A'}]\circ[T] ,
  \]
  so the action on morphisms of the left adjoint obtained from the universal
  property agrees with $\mathcal{D}_{\mathrm{ae}}$.
\end{proof}

\begin{corollary}
  The functor
  \[
    \mathcal{D}_{M,\mathrm{ae}}:
    \mathbf{MInt}_{\mathrm{ae}}
    \to
    \mathbf{DInt}_{\mathrm{ae}}
  \]
  obtained by restricting $\mathcal{D}_{\mathrm{ae}}$ to
  $\mathbf{MInt}_{\mathrm{ae}}$ is left adjoint to the inclusion functor
  \[
    \iota_{M,\mathrm{ae}}:
    \mathbf{DInt}_{\mathrm{ae}}
    \to
    \mathbf{MInt}_{\mathrm{ae}} .
  \]
\end{corollary}
\begin{proof}
  $\mathbf{MInt}_{\mathrm{ae}}$ is a full subcategory of
  $\mathbf{Int}_{\mathrm{ae}}$, and $\mathbf{DInt}_{\mathrm{ae}}$ is a full
  subcategory of it.
  It therefore suffices to restrict the adjunction isomorphism of
  Theorem~\ref{thm:main-ae} to the objects of
  $\mathbf{MInt}_{\mathrm{ae}}$.
\end{proof}

\begin{remark}\label{rem:ae-mct-counterexample-fails}
  Dividing morphisms by ae-equivalence does not weaken the membership condition
  for functions in the target of the factorization.
  Accordingly, in the theorem above as in the main theorem, Daniell closedness is
  assumed of the target of the factorization.

  However, \textbf{the finite counterexample of
  Remark~\ref{rem:D-I-not-left-adjoint} cannot be used as it stands as grounds
  for the claim that Daniell closedness cannot be dropped.}
  In that example only $\tau(0)=0$ has any effect, so all $5$ morphisms and all
  $3$ morphisms are ae-equivalent to one another, and in the ae-category both
  $\operatorname{Hom}$ sets collapse to a single point.
  A counterexample for the ae version is given in the next subsubsection
  (Corollary~\ref{cor:ae-completeness-needed}).

  This point contrasts with Section~\ref{chap:operator}.
  There, once morphisms are replaced by operators, the extension is determined
  only up to null functions, and dividing the codomain by null functions becomes
  the natural thing to do.
  Even there, however, we do not prove that ``it is undetermined unless one
  divides'' (Problem~\ref{prob:no-lift}).
\end{remark}

\subsubsection{Daniell closedness cannot be dropped in the ae version either}

As stated in Remark~\ref{rem:ae-mct-counterexample-fails}, the finite
counterexample of Remark~\ref{rem:D-I-not-left-adjoint} has no effect in the
ae-category.
A counterexample for the ae version is given by the pair consisting of the
integration system of Borel measurable functions and the integration system of
Lebesgue integrable functions.

\begin{proposition}\label{prop:no-morphism-lebesgue-to-borel}
  Let $X=[0,1)$ and put
  \begin{align*}
    P&:=\bigl(X,\ \{f:X\to\mathbb{R}\text{ Borel measurable and Lebesgue integrable}\},\ \text{the Lebesgue integral}\bigr),\\
    Q&:=\bigl(X,\ \{f:X\to\mathbb{R}\text{ Lebesgue integrable}\},\ \text{the Lebesgue integral}\bigr).
  \end{align*}
  Then
  \[
    \operatorname{Hom}_{\mathbf{Int}}(Q,P)=\emptyset.
  \]
\end{proposition}
\begin{proof}
  Suppose a morphism $F:Q\to P$ exists and let $\tau:X\to X$ be its map of base
  spaces.
  By definition, for every Lebesgue integrable $h:X\to\mathbb{R}$ the function
  $h\circ\tau$ is Borel measurable and Lebesgue integrable, and
  \[
    \int_X h\circ\tau\,d\lambda=\int_X h\,d\lambda
  \]
  holds, where $\lambda$ is Lebesgue measure.

  \textbf{$\tau$ is Borel measurable and measure-preserving.}
  For a Borel set $B\subset X$ the function $1_B$ is Lebesgue integrable, and
  since $1_B\circ\tau=1_{\tau^{-1}(B)}$ is Borel measurable, the set
  $\tau^{-1}(B)$ is Borel.
  Moreover preservation of the integral gives
  $\lambda(\tau^{-1}(B))=\lambda(B)$.

  \textbf{The image of $\tau$ is uncountable.}
  If $\tau(X)$ were contained in a Borel null set $N_0$, then
  $\tau^{-1}(X\setminus N_0)=\emptyset$ while
  $\lambda(X\setminus N_0)=1$, contradicting preservation of measure.
  A countable set is contained in a Borel null set, so $\tau(X)$ is
  uncountable.

  \textbf{The image of $\tau$ contains a Lebesgue null set of cardinality
  $\mathfrak{c}$.}
  Since $\tau$ is Borel measurable, its image $\tau(X)$ is an analytic set.
  An uncountable analytic set contains a nonempty perfect set (the perfect set
  property of analytic sets; see \cite{Kechris}).
  Carrying out a Cantor-type construction inside such a perfect set $K$, with
  all $2^n$ closed intervals at stage $n$ of diameter at most $4^{-n}$, the
  resulting intersection $P_0\subset K$ is a set of cardinality $\mathfrak{c}$,
  and $\lambda(P_0)\le2^n\cdot4^{-n}=2^{-n}$ for every $n$, so
  $\lambda(P_0)=0$.

  \textbf{Contradiction.}
  Take an arbitrary $S\subset P_0$. As a subset of a Lebesgue null set, $S$ is
  Lebesgue measurable, and $1_S$ is Lebesgue integrable.
  Hence $1_S\circ\tau=1_{\tau^{-1}(S)}$ is Borel measurable, that is,
  $\tau^{-1}(S)$ is a Borel set.

  The assignment $S\mapsto\tau^{-1}(S)$ is injective on the subsets of $P_0$.
  Indeed, if $S\ne S'$ then, without loss of generality, we may take
  $p\in S\setminus S'$.
  Since $p\in P_0\subset\tau(X)$ there is an $x\in X$ with $\tau(x)=p$, and
  $x\in\tau^{-1}(S)\setminus\tau^{-1}(S')$.

  We thus obtain an injection from $\mathcal{P}(P_0)$ into the set of all Borel
  subsets of $X$.
  But $|\mathcal{P}(P_0)|=2^{\mathfrak{c}}$ while $X$ has only $\mathfrak{c}$
  Borel sets (the same fact as used in the proof of
  Corollary~\ref{cor:lebesgue-reflections-differ}).
  This is a contradiction.
\end{proof}

\begin{corollary}\label{cor:ae-completeness-needed}
  The functor
  \[
    \iota_{M,\mathrm{ae}}\circ\mathcal{D}_{\mathrm{ae}}:
    \mathbf{Int}_{\mathrm{ae}}\longrightarrow\mathbf{MInt}_{\mathrm{ae}}
  \]
  defined by the Daniell extension is not left adjoint to the inclusion functor
  $\iota_{\mathrm{ae}}:\mathbf{MInt}_{\mathrm{ae}}\to\mathbf{Int}_{\mathrm{ae}}$
  of the full subcategory.
  That is, \textbf{imposing Daniell closedness on the codomain of the
  factorization cannot be dispensed with in the ae version either}.
\end{corollary}
\begin{proof}
  Taking $A$ to be the integration system of
  Example~\ref{ex:step-functions}, by
  Corollary~\ref{cor:lebesgue-reflections-differ} the object $P$ of
  Proposition~\ref{prop:no-morphism-lebesgue-to-borel} is $A^{\mathrm{M}}$ and
  $Q$ is $\mathcal{D}A$.
  Here $P$ is an object of $\mathbf{MInt}$
  (Theorem~\ref{thm:mint-inclusion-has-left-adjoint}) and is not Daniell closed
  (the same corollary).
  Moreover $\mathcal{D}P=\mathcal{D}(A^{\mathrm{M}})=\mathcal{D}A=Q$ by
  Lemma~\ref{lem:two-reflections-same-daniell}.

  The morphisms of $\mathbf{Int}_{\mathrm{ae}}$ are the ae-equivalence classes of
  morphisms of $\mathbf{Int}$, so
  Proposition~\ref{prop:no-morphism-lebesgue-to-borel} gives
  \[
    \operatorname{Hom}_{\mathbf{MInt}_{\mathrm{ae}}}(\mathcal{D}_{\mathrm{ae}}P,\,P)=\emptyset .
  \]
  On the other hand
  $\operatorname{Hom}_{\mathbf{Int}_{\mathrm{ae}}}(P,\iota_{\mathrm{ae}}P)$
  contains the ae-equivalence class of the identity morphism and is therefore
  nonempty.
  There is no bijection between the empty set and a nonempty set, so no
  adjunction isomorphism exists.
\end{proof}

\begin{remark}
  This proof is the only place in this paper that uses a result of descriptive
  set theory (the perfect set property of analytic sets).
  This contrasts with the non-ae counterexample
  (Remark~\ref{rem:D-I-not-left-adjoint}), for which a three-point set
  sufficed.
  Once morphisms are divided by ae-equivalence, finite objects can no longer
  tell them apart, and a counterexample appears only in the form of the
  difference between Borel and Lebesgue.
\end{remark}

\section{Morphisms as operators}\label{chap:operator}

The morphisms used in Section~\ref{chap:functoriality} were the maps of base
spaces $\underline{F}:Y\to X$ themselves.
Under that definition a $\operatorname{Hom}$ set is no more than a subset of
$\operatorname{Map}(Y,X)$, so the uniqueness of the factorization in
Theorem~\ref{thm:main-theorem} followed automatically from the agreement of the
maps of base spaces.
That is, for the non-ae version of the main theorem the uniqueness part has no
content (this is not so for the ae version; see
Subsection~\ref{sec:main-theorem-ae}).

In Section~\ref{chap:operator} we detach morphisms from maps of base spaces and
redefine them as linear maps between Riesz spaces.

\begin{quote}
  \textbf{A note on attribution.}
  This formulation is based on the one introduced in \cite{Yamagami} in the
  context of discussing iterated integrals, but it is \textbf{not the same
  thing.}
  That book imposes \emph{order continuity} in addition to linearity and
  positivity, whereas the operator morphisms of integration systems in
  Section~\ref{chap:operator} (Definition~\ref{def:operator-morphism}) do not
  impose it.
  That is, what we adopt is that book's definition with one condition removed.

  Order continuity is a condition that book imposes in order to
  \emph{construct} $J\circ\Lambda$ as an integral, and it becomes unnecessary in
  the setting of Section~\ref{chap:operator}, which requires preservation of the
  integral.
  For details see Remark~\ref{rem:order-continuity-redundant}.
\end{quote}

\noindent
Note that the definition of \cite{Yamagami}, with order continuity imposed
(Remark~\ref{rem:yamagami-operator}), includes Markov operators and transition
probabilities, and that book says so explicitly.

This replacement changes the following three things.

\begin{itemize}
  \item Morphisms no longer pass through the base space, and a
    $\operatorname{Hom}$ set is no longer a subset of a set of maps.
  \item Uniqueness of the factorization no longer follows automatically; it is
    obtained using the density of $L$ in $L^1$
    (Lemma~\ref{lem:l1-approx-by-source}).
    The same density is used in the ae version in
    Subsection~\ref{sec:main-theorem-ae}.
  \item The extension is determined only up to null functions.
    It therefore becomes natural to divide the codomain by null functions
    (Subsection~\ref{sec:ae}).
    We do not, however, prove that ``it is undetermined unless one divides''
    (Problem~\ref{prob:no-lift}).
\end{itemize}

In the end we show that the Daniell extension gives a reflection in the category
with operator morphisms taken in this way as well
(Theorem~\ref{thm:operator-main-theorem}).

\subsection{Operator morphisms}

\begin{definition}\label{def:operator-morphism}
  Let $(X,L,I)$ and $(Y,M,J)$ be integration systems on sets $X$ and
  $Y$ respectively.
  An \textbf{operator morphism} of integration systems
  \[
    \Lambda:(X,L,I)\to(Y,M,J)
  \]
  is a linear map $\Lambda:L\to M$ satisfying the following two conditions.
  \begin{itemize}
    \item (positivity) for every $f\in L$, if $f\ge0$ then $\Lambda f\ge0$;
    \item (preservation of the integral) $J(\Lambda f)=I(f)$ for every $f\in L$.
  \end{itemize}
\end{definition}

\begin{remark}[relation to the operators of \cite{Yamagami}]\label{rem:yamagami-operator}
  \cite{Yamagami} introduces linear maps $\Lambda$ from a Riesz space $L$ on a
  set $X$ to a Riesz space $M$ on a set $Y$ satisfying
  \begin{itemize}
    \item (positivity) for every $f\in L$, if $f\ge0$ then $\Lambda f\ge0$;
    \item (order continuity) for every sequence of functions
      $\{f_n\}\subset L$, if $f_n(x)\downarrow0$ for all $x\in X$ then
      $(\Lambda f_n)(y)\downarrow0$ for all $y\in Y$.
  \end{itemize}
  Definition~\ref{def:operator-morphism} adds preservation of the integral to
  this and \textbf{removes order continuity}.

  Anything that is an operator in the sense of that book and preserves the
  integral is an operator morphism in the sense of
  Definition~\ref{def:operator-morphism}. The converse does not hold in general.
\end{remark}

\begin{remark}\label{rem:order-continuity-redundant}
  Note that we do not impose the order continuity of
  Remark~\ref{rem:yamagami-operator}.
  \cite{Yamagami} imposes it in order to \emph{construct} an integral
  $J\circ\Lambda$ on $L$ out of $J$ and $\Lambda$. Indeed, if
  $f_n(x)\downarrow0$ then $(\Lambda f_n)(y)\downarrow0$, and since $J$ is an
  integral, $J(\Lambda f_n)\to0$ follows.

  In this paper, by contrast, the integral $I$ is given from the outset as part
  of the object $(X,L,I)$.
  In that situation preservation of the integral, $J\circ\Lambda=I$,
  \textbf{already implies that $J\circ\Lambda$ is an integral}.
  That is, order continuity is unnecessary for the purpose of making
  $J\circ\Lambda$ an integral.

  Order continuity is a genuinely strong condition on the pointwise behaviour of
  $\Lambda f_n$, but it is never used in the discussion below.
  And in fact, for the morphisms of this category, pointwise conditions cannot
  be stated as propositions at all.
  The values of the extension given by
  Theorem~\ref{thm:operator-extension} are elements of the quotient
  $\widetilde{M}$, so a pointwise condition would depend on the choice of
  representative.
  That is, a pointwise condition is not false but meaningless.
  Imposing pointwise conditions at all is inappropriate in this category.

  Note also that the fact that the condition in \cite{Yamagami} is stated in
  terms of pointwise convergence is essential in that book's setting.
  If it were replaced by convergence in the order of $M$ then, since order
  convergence is strictly weaker than pointwise convergence (see the example in
  Appendix~\ref{sec:forget}), $J(\Lambda f_n)\to0$ would no longer follow.
\end{remark}

\begin{lemma}
  \begin{enumerate}
    \item For every integration system $(X,L,I)$, the identity map
      $\operatorname{id}_L:L\to L$ is an operator morphism
      $(X,L,I)\to(X,L,I)$.
    \item For operator morphisms $\Lambda:(X,L,I)\to(Y,M,J)$ and
      $\Xi:(Y,M,J)\to(Z,N,K)$, the composite $\Xi\circ\Lambda:L\to N$ is an
      operator morphism $(X,L,I)\to(Z,N,K)$.
  \end{enumerate}
\end{lemma}
\begin{proof}
  (1) is immediate. As for (2), a composite of linear maps is linear and a
  composite of positive maps is positive.
  Preservation of the integral holds because
  $K(\Xi\Lambda f)=J(\Lambda f)=I(f)$.
\end{proof}

\begin{definition}
  The category whose objects are integration systems and whose morphisms are
  operator morphisms of integration systems is written
  \[
    \mathbf{Int}_{\mathrm{op}} .
  \]
\end{definition}

\subsubsection{Relation to morphisms given by maps of base spaces}

\begin{lemma}\label{lem:base-map-gives-operator}
  Let $F:(X,L,I)\to(Y,M,J)$ be a morphism of integration systems in the sense of
  Section~\ref{chap:preliminaries}.
  Then
  \[
    F^*:L\to M;\quad f\mapsto f\circ\underline{F}
  \]
  is an operator morphism of integration systems.
\end{lemma}
\begin{proof}
  Linearity is immediate from the definition of composition.
  If $f\ge0$ then $(F^*f)(y)=f(\underline{F}(y))\ge0$ for every $y\in Y$, so
  $F^*$ is positive.
  Preservation of the integral is precisely the definition of a morphism of
  integration systems.
  Note that $F^*$ is order continuous as well. Indeed, if $f_n(x)\downarrow0$
  holds for every $x\in X$ then, putting $x:=\underline{F}(y)$ for each
  $y\in Y$, we get $(F^*f_n)(y)=f_n(x)\downarrow0$.
\end{proof}

\begin{proposition}\label{prop:base-to-operator-functor}
  The assignment that leaves objects unchanged and sends a morphism $F$ to
  $F^*$ defines a functor
  \[
    \mathcal{E}:\mathbf{Int}\to\mathbf{Int}_{\mathrm{op}} .
  \]
\end{proposition}
\begin{proof}
  By Lemma~\ref{lem:base-map-gives-operator} morphisms are sent to operator
  morphisms.
  For identity morphisms, $(\operatorname{id}_X)^*=\operatorname{id}_L$.
  As for composition, the map of base spaces underlying $G\circ F$ was
  $\underline{F}\circ\underline{G}$, so for every $f\in L$
  \[
    (G\circ F)^*f
    =f\circ(\underline{F}\circ\underline{G})
    =(f\circ\underline{F})\circ\underline{G}
    =G^*(F^*f) ,
  \]
  and $(G\circ F)^*=G^*\circ F^*$ holds.
\end{proof}

\begin{remark}
  The functor $\mathcal{E}$ need not be faithful.
  Indeed, when $L$ does not separate the points of $X$, distinct maps of base
  spaces $\underline{F}\ne\underline{G}$ can give the same pullback
  $F^*=G^*$.
  In this sense $\mathbf{Int}_{\mathrm{op}}$ does not carry from the outset the
  information about the base space that the integral does not distinguish.
\end{remark}

\begin{remark}
  $\mathcal{E}$ is not full either.
  Among operator morphisms, those coming from maps of base spaces are confined
  to the ones preserving the lattice operations $\wedge,\vee$ (Riesz
  homomorphisms).
  A general operator morphism need not be a lattice homomorphism, and, as
  \cite{Yamagami} points out, this difference is what decides whether
  probability kernels and Markov operators are included.
  The next example exhibits this difference concretely.
\end{remark}

\begin{example}[the case of finite sets: column-stochastic matrices]\label{ex:stochastic-matrix}
  Let $X,Y$ be finite sets and let $L=\mathbb{R}^X$ and $M=\mathbb{R}^Y$.
  Take positive weights $p:X\to(0,\infty)$ and $q:Y\to(0,\infty)$ and put
  \[
    I(f):=\sum_{x\in X}p(x)f(x),
    \qquad
    J(g):=\sum_{y\in Y}q(y)g(y) ;
  \]
  then $(X,L,I)$ and $(Y,M,J)$ are integration systems.

  A linear map $\Lambda:L\to M$ is represented by a matrix $(\Lambda_{yx})$.
  Then
  \begin{itemize}
    \item positivity is equivalent to all entries satisfying
      $\Lambda_{yx}\ge0$;
    \item preservation of the integral, $J(\Lambda f)=I(f)$, says that
      $\sum_{x}\bigl(\sum_{y}q(y)\Lambda_{yx}\bigr)f(x)=\sum_x p(x)f(x)$ holds
      for every $f$, hence is equivalent to
      \[
        \sum_{y\in Y}q(y)\Lambda_{yx}=p(x)
        \qquad(\forall x\in X) .
      \]
  \end{itemize}

  In particular, in the case $p\equiv1$ and $q\equiv1$ (counting, the finite
  version of Example~\ref{ex:counting}), the latter condition says
  \[
    \sum_{y\in Y}\Lambda_{yx}=1
    \qquad(\forall x\in X) ,
  \]
  that is, that \textbf{each column sums to $1$}.
  Hence in this case the operator morphisms of integration systems are precisely
  the \textbf{column-stochastic matrices}.
  Below, ``stochastic matrix'' in this paper means a column-stochastic matrix.

  On the other hand, for a morphism coming from a map of base spaces
  $\underline{F}:Y\to X$ we have $\Lambda_{yx}=1$ when
  $\underline{F}(y)=x$ and $\Lambda_{yx}=0$ otherwise, so each row contains
  exactly one $1$.
  In the counting case, moreover, the column-sum condition forces
  $|\underline{F}^{-1}(x)|=1$, that is, that $\underline{F}$ be a bijection.

  That is, between the same objects,
  \begin{center}
    with maps of base spaces as morphisms only the ``permutation matrices'',\\
    with operators as morphisms all ``column-stochastic matrices''
  \end{center}
  are morphisms. This is the simplest manifestation of the fact that
  $\mathcal{E}$ is not full.
\end{example}

\begin{remark}
  Imposing $\Lambda1_X=1_Y$ (row sums equal to $1$) in the setting of
  Example~\ref{ex:stochastic-matrix} gives what \cite{Yamagami} calls a
  \textbf{Markov operator}, the corresponding family of measures being a
  \textbf{transition probability}.

  Note, however, that in the counting case imposing preservation of the integral
  (column sums equal to $1$) and row sums equal to $1$ simultaneously yields a
  \textbf{doubly stochastic matrix}.
  Hence one cannot say that ``$\mathbf{Int}_{\mathrm{op}}$ contains all
  transition matrices of Markov chains''.
  For general weights $p,q$, the transition matrices that are morphisms are
  confined to those with $\sum_{y}q(y)\Lambda_{yx}=p(x)$, that is, those sending
  $q$ to $p$.
  What is included is exactly those satisfying this condition.
\end{remark}

\subsection{The Daniell extension of operator morphisms}

The main result of this subsection is
Theorem~\ref{thm:operator-extension}.
Unlike in Section~\ref{chap:functoriality},
\textbf{uniqueness does not follow from the agreement of maps of base spaces};
it is obtained by way of contractivity and the density of $L$ in $L^1$.

In what follows, for $(Y,M,J)$ the symbols $\mathcal{N}_J$,
$\widetilde{M}=M/\mathcal{N}_J$, $\pi_J$ and $\widetilde{J}$ are those of
Subsection~\ref{sec:ae}.

\subsubsection{Contractivity}

Below, we say that $\Lambda:L\to\widetilde{M}$ is an \textbf{operator morphism
into the quotient} if $\Lambda$ is linear, positive
($f\ge0\Rightarrow\Lambda f\ge0$) and satisfies
$\widetilde{J}(\Lambda f)=I(f)$.
For an ordinary operator morphism $\Lambda:L\to M$, the composite
$\pi_J\circ\Lambda$ is an operator morphism into the quotient.

\begin{lemma}\label{lem:operator-dominates-abs}
  Let $\Lambda:L\to\widetilde{M}$ be linear and positive.
  Then $|\Lambda f|\le\Lambda|f|$ holds for every $f\in L$.
\end{lemma}
\begin{proof}
  We have $f\le|f|$ and $-f\le|f|$. Since $\Lambda$ is linear and positive,
  $\Lambda f\le\Lambda|f|$ and $-\Lambda f\le\Lambda|f|$, and in the Riesz space
  $\widetilde{M}$ we obtain
  $|\Lambda f|=(\Lambda f)\vee(-\Lambda f)\le\Lambda|f|$.
\end{proof}

\begin{lemma}[contractivity]\label{lem:operator-is-contraction}
  Let $\Lambda:L\to\widetilde{M}$ be an operator morphism into the quotient.
  Then
  \[
    \|\Lambda f\|\le\|f\|_I
  \]
  holds for every $f\in L$. Moreover equality holds when $f\ge0$.
\end{lemma}
\begin{proof}
  By Lemma~\ref{lem:operator-dominates-abs} and positivity of $\widetilde{J}$,
  \[
    \|\Lambda f\|=\widetilde{J}(|\Lambda f|)\le\widetilde{J}(\Lambda|f|)=I(|f|)=\|f\|_I.
  \]
  Here the second equality is preservation of the integral for $|f|\in L$.
  When $f\ge0$ we have $|f|=f$ and $\Lambda f\ge0$, so
  $\|\Lambda f\|=\widetilde{J}(\Lambda f)=I(f)=\|f\|_I$.
\end{proof}

\begin{remark}
  That equality does not hold in general is the difference from a lattice
  homomorphism.
  If $\Lambda$ is a lattice homomorphism then $|\Lambda f|=\Lambda|f|$ and it is
  isometric, whereas a general operator morphism remains merely contractive.
  All that the discussion below needs is contractivity.
\end{remark}

\subsubsection{The extension theorem}

\begin{theorem}[Daniell extension of operator morphisms]\label{thm:operator-extension}
  Let $(X,L,I)$ be an integration system with Daniell extension
  $(X,L^1,I^1)$, and let $(Y,M,J)$ be an integration system
  \textbf{such that $\widetilde{M}$ is a Banach space}.
  Let $\Lambda:L\to\widetilde{M}$ be an operator morphism into the quotient.
  Then there \textbf{exists a unique} operator morphism into the quotient
  \[
    \Lambda^1:L^1\to\widetilde{M}
  \]
  satisfying $\Lambda^1|_L=\Lambda$.
  Moreover $\Lambda^1$ is contractive with respect to $\|\cdot\|$.
\end{theorem}
\begin{proof}
  \textbf{Existence.}
  Let $f\in L^1$. By Lemma~\ref{lem:l1-approx-by-source} we may take
  $\{s_N\}\subset L$ with $\|f-s_N\|_{I^1}\to0$.
  A convergent sequence is a Cauchy sequence, and
  $\|s_N-s_M\|_{I^1}=\|s_N-s_M\|_I$ for $s_N,s_M\in L$
  (Lemma~\ref{lem:daniell-ext-contains-source}).
  Hence Lemma~\ref{lem:operator-is-contraction} gives
  \[
    \|\Lambda s_N-\Lambda s_M\|
    =\|\Lambda(s_N-s_M)\|
    \le\|s_N-s_M\|_I
    \longrightarrow0 ,
  \]
  so $\{\Lambda s_N\}$ is a Cauchy sequence in $\widetilde{M}$.
  By hypothesis $\widetilde{M}$ is a Banach space, so this sequence converges.
  We define $\Lambda^1(f)$ to be its limit.

  \textbf{Independence of the approximating sequence.}
  If $\{t_N\}\subset L$ also satisfies $\|f-t_N\|_{I^1}\to0$ then
  $\|s_N-t_N\|_{I^1}\to0$, and contractivity gives
  $\|\Lambda s_N-\Lambda t_N\|\to0$.
  As $\widetilde{M}$ is a normed space and limits are unique, the two limits
  agree.

  \textbf{Linearity, contractivity, and agreement on $L$.}
  For $f,g\in L^1$ and $a\in\mathbb{R}$, taking approximating sequences
  $\{s_N\},\{u_N\}$ for each, the sequence $\{s_N+au_N\}$ is an approximating
  sequence for $f+ag$, whence linearity.
  Moreover $\|\Lambda s_N\|\le\|s_N\|_I=\|s_N\|_{I^1}\to\|f\|_{I^1}$, so
  continuity of the norm gives $\|\Lambda^1 f\|\le\|f\|_{I^1}$.
  For $f\in L$, taking the constant sequence $s_N:=f$ gives
  $\Lambda^1 f=\Lambda f$.

  \textbf{Positivity.}
  Let $f\in L^1$ with $f\ge0$.
  For an approximating sequence $\{s_N\}$ put $s_N^+:=s_N\vee0\in L$.
  Since $f\ge0$ we have $|f-s_N^+|\le|f-s_N|$ pointwise, hence
  $\|f-s_N^+\|_{I^1}\le\|f-s_N\|_{I^1}\to0$.
  So $\{s_N^+\}$ is also an approximating sequence for $f$, and
  $\Lambda^1 f=\lim_N\Lambda s_N^+$.
  For each $N$, $s_N^+\ge0$ gives $\Lambda s_N^+\ge0$.
  The positive cone of $\widetilde{M}$ is closed with respect to $\|\cdot\|$
  (indeed $\|(\Lambda^1 f)\wedge0\|\le\|\Lambda^1 f-\Lambda s_N^+\|\to0$), so
  $\Lambda^1 f\ge0$.

  \textbf{Preservation of the integral.}
  $\widetilde{J}$ satisfies $|\widetilde{J}(v)|\le\|v\|$ and is therefore
  continuous.
  Hence $\widetilde{J}(\Lambda^1 f)=\lim_N\widetilde{J}(\Lambda s_N)
  =\lim_N I(s_N)=\lim_N I^1(s_N)=I^1(f)$.
  The last equality holds because
  $|I^1(f)-I^1(s_N)|\le I^1(|f-s_N|)\to0$.

  \textbf{Uniqueness.}
  Suppose $\Psi:L^1\to\widetilde{M}$ is an operator morphism into the quotient
  with $\Psi|_L=\Lambda$.
  Applying Lemma~\ref{lem:operator-is-contraction} to $(X,L^1,I^1)$, the map
  $\Psi$ is contractive and hence continuous with respect to
  $\|\cdot\|_{I^1}$.
  $\Psi$ and $\Lambda^1$ agree on $L$, and $L$ is dense in $L^1$ by
  Lemma~\ref{lem:l1-approx-by-source}, so as two continuous maps they agree on
  all of $L^1$.
\end{proof}

\begin{remark}\label{rem:uniqueness-has-content}
  Note that the proof of uniqueness differs essentially from the case of the
  main theorem of Section~\ref{chap:functoriality}.
  There, morphisms were the maps of base spaces themselves, so uniqueness of the
  factorizing morphism followed automatically from the agreement of the maps of
  base spaces.
  Here \textbf{contractivity
  (Lemma~\ref{lem:operator-is-contraction}) and density
  (Lemma~\ref{lem:l1-approx-by-source}) actually do the work.}

  Section~\ref{chap:operator} is not, however, the first place where an analytic
  lemma from the preliminary section contributes to a proof of a universal
  property.
  In the ae version of the main theorem in
  Subsection~\ref{sec:main-theorem-ae}, the same density lemma is used by way of
  Lemma~\ref{lem:daniell-extension-preserves-ae}.
  Precisely, \textbf{in the non-ae version uniqueness follows automatically from
  the agreement of the maps of base spaces, while the ae version and the
  operator version of Section~\ref{chap:operator} require an analytic
  argument}.
\end{remark}

\begin{remark}\label{rem:ae-is-forced}
  What Theorem~\ref{thm:operator-extension} gives is a map into
  $\widetilde{M}$, not a map into $M$.
  Choosing representatives in $M$ for the elements of $\widetilde{M}$ yields a
  map $L^1\to M$, but this is determined only up to null functions.
  Hence, in order to discuss the Daniell extension in a category whose morphisms
  are operator morphisms, it is natural to divide the codomain by null
  functions.

  This paper cannot, however, assert that ``the extension is undetermined unless
  one divides''.
  That would require an \textbf{example of an operator morphism into the
  quotient admitting no linear positive lift $L^1\to M$}, and this paper does
  not supply one (Problem~\ref{prob:no-lift}).
  Indeed, for morphisms coming from maps of base spaces a lift
  $F^*:L^1\to M$ does exist, so there is no need to pass to ae-equivalence in
  the setting of Section~\ref{chap:functoriality}.

  Note that the situation in which the extension is determined only up to null
  sets occurs likewise in the theorem on iterated integrals in
  \cite{Yamagami}, and that book states the extension for integrable functions
  in that form.
\end{remark}

\begin{problem}\label{prob:no-lift}
  Is there an operator morphism into the quotient $L^1\to\widetilde{M}$
  admitting no linear positive lift $L^1\to M$?

  If one exists, it would show that passing to ae-equivalence is a requirement,
  and at the same time add one more reason to keep
  Section~\ref{chap:functoriality} as an independent section (see
  Subsection~\ref{sec:two-chapters}).
\end{problem}

\subsubsection{The ae-category and the main theorem}

\begin{lemma}[an operator morphism into the quotient is determined by the quotient]\label{lem:operator-descends}
  Let $\Lambda:L\to\widetilde{M}$ be an operator morphism into the quotient.
  Then $\Lambda$ maps $\mathcal{N}_I$ to $0$, and there is a unique linear map
  $\widetilde{\Lambda}:\widetilde{L}\to\widetilde{M}$ satisfying
  $\widetilde{\Lambda}\circ\pi_I=\Lambda$.
  Moreover $\widetilde{\Lambda}$ is positive and satisfies
  $\widetilde{J}(\widetilde{\Lambda}(\pi_I(f)))=\widetilde{I}(\pi_I(f))$ for
  every $f\in L$.
\end{lemma}
\begin{proof}
  If $f\in\mathcal{N}_I$ then $\|f\|_I=0$, and
  Lemma~\ref{lem:operator-is-contraction} gives
  $\|\Lambda f\|\le\|f\|_I=0$.
  As $\|\cdot\|$ is a norm on $\widetilde{M}$ we get $\Lambda f=0$.
  Hence $\Lambda$ maps $\mathcal{N}_I$ to $0$ and a unique linear map
  $\widetilde{\Lambda}$ is determined.

  We check positivity. If $\pi_I(f)\ge0$ then, by the way the order on
  $\widetilde{L}$ is defined (Lemma~\ref{lem:quotient-riesz}), we may take
  $f'\in\mathcal{N}_I$ with $f+f'\ge0$.
  Since $\Lambda f'=0$,
  \[
    \widetilde{\Lambda}(\pi_I(f))=\Lambda(f+f')\ge0 .
  \]
  Finally
  $\widetilde{J}(\widetilde{\Lambda}(\pi_I(f)))=\widetilde{J}(\Lambda f)=I(f)=\widetilde{I}(\pi_I(f))$.
\end{proof}

\begin{definition}\label{def:int-op-ae}
  The category whose objects are integration systems and whose morphisms from
  $(X,L,I)$ to $(Y,M,J)$ are the operator morphisms into the quotient
  $\Lambda:L\to\widetilde{M}$ is written
  \[
    \mathbf{Int}_{\mathrm{op},\mathrm{ae}} .
  \]
  Composition is defined as follows.
  For $\Lambda:L\to\widetilde{M}$ and $\Xi:M\to\widetilde{N}$,
  Lemma~\ref{lem:operator-descends} shows that $\Xi$ induces a linear positive
  map $\widetilde{\Xi}:\widetilde{M}\to\widetilde{N}$.
  We define the composite to be $\widetilde{\Xi}\circ\Lambda$.
  This is linear and positive and satisfies
  $\widetilde{K}(\widetilde{\Xi}(\Lambda f))=\widetilde{J}(\Lambda f)=I(f)$, so
  it is again an operator morphism into the quotient.
  The identity morphism is $\pi_I:L\to\widetilde{L}$.

  We write the full subcategory of objects for which $\widetilde{L}$ is a Banach
  space as
  \[
    \mathbf{BInt}_{\mathrm{op},\mathrm{ae}} .
  \]
  Daniell closed integration systems belong to this subcategory by
  Lemma~\ref{lem:quotient-is-banach}.
  The reason for taking this subcategory is given in
  Subsection~\ref{sec:what-op-ae-sees}.
\end{definition}

\begin{remark}
  For ordinary operator morphisms $\Lambda,\Lambda':L\to M$, defining the same
  morphism $\pi_J\circ\Lambda=\pi_J\circ\Lambda'$ is equivalent to
  $J(|\Lambda f-\Lambda' f|)=0$ holding for every $f\in L$.
  The latter is precisely the relation $=_{\mathrm{ae}}$ of
  Definition~\ref{def:ae-equal} transported to operator morphisms.

  It is not accurate, however, to call
  $\mathbf{Int}_{\mathrm{op},\mathrm{ae}}$ ``the category obtained by dividing
  morphisms by ae-equivalence''.
  Its morphisms are defined from the outset as \emph{all} operator morphisms
  into the quotient $L\to\widetilde{M}$, and not all of them need come from
  $M$-valued morphisms (Remark~\ref{rem:ae-is-forced}).
  Precisely, the relation is that \textbf{there is a faithful functor from the
  category obtained by dividing ordinary operator morphisms by ae-equivalence
  into $\mathbf{Int}_{\mathrm{op},\mathrm{ae}}$}.
\end{remark}

\begin{theorem}[operator-morphism version of the main theorem]\label{thm:operator-main-theorem}
  The inclusion functor
  \[
    \iota_{B,\mathrm{op}}:
    \mathbf{BInt}_{\mathrm{op},\mathrm{ae}}
    \to
    \mathbf{Int}_{\mathrm{op},\mathrm{ae}}
  \]
  has a left adjoint
  \[
    \mathcal{D}_{\mathrm{op}}:
    \mathbf{Int}_{\mathrm{op},\mathrm{ae}}
    \to
    \mathbf{BInt}_{\mathrm{op},\mathrm{ae}} .
  \]
  Its action on objects is the Daniell extension.
  Its action on morphisms sends a morphism $\Lambda:A\to B$ to the extension to
  $L^1\to\widetilde{M^1}$, given by
  Theorem~\ref{thm:operator-extension}, of the composite
  $i_B\circ\Lambda:L\to\widetilde{M^1}$
  (since $\mathcal{D}B$ is Daniell closed, $\widetilde{M^1}$ is a Banach space
  by Lemma~\ref{lem:quotient-is-banach} and that theorem applies).
\end{theorem}
\begin{proof}
  Let $A=(X,L,I)$ be an object of
  $\mathbf{Int}_{\mathrm{op},\mathrm{ae}}$ and let
  $\mathcal{D}A:=(X,L^1,I^1)$ be its Daniell extension.
  By Lemma~\ref{lem:daniell-ext-is-complete} the system $\mathcal{D}A$ is
  Daniell closed, and by Lemma~\ref{lem:quotient-is-banach} the space
  $\widetilde{L^1}$ is a Banach space, so $\mathcal{D}A$ is an object of
  $\mathbf{BInt}_{\mathrm{op},\mathrm{ae}}$.
  Moreover, from the inclusion $L\subset L^1$ and $I^1|_L=I$
  (Lemma~\ref{lem:daniell-ext-contains-source}), the map
  \[
    i_A:=\pi_{I^1}|_L:L\to\widetilde{L^1}
  \]
  is an operator morphism into the quotient and defines a morphism
  $i_A:A\to\iota_{B,\mathrm{op}}(\mathcal{D}A)$ of
  $\mathbf{Int}_{\mathrm{op},\mathrm{ae}}$.

  We show that $i_A$ is a universal morphism from $A$ to
  $\iota_{B,\mathrm{op}}$.
  Let $B=(Y,M,J)$ be an object of
  $\mathbf{BInt}_{\mathrm{op},\mathrm{ae}}$ and let
  $\Lambda:A\to\iota_{B,\mathrm{op}}(B)$ be a morphism, that is, an operator
  morphism into the quotient $\Lambda:L\to\widetilde{M}$.
  Since $\widetilde{M}$ is a Banach space,
  Theorem~\ref{thm:operator-extension} applies and there is a unique operator
  morphism into the quotient $\Lambda^1:L^1\to\widetilde{M}$ with
  $\Lambda^1|_L=\Lambda$.

  Write $\Lambda^\sharp$ for the morphism $\mathcal{D}A\to B$ determined by
  $\Lambda^1$.
  By the definition of composition in
  Definition~\ref{def:int-op-ae}, the composite
  $\iota_{B,\mathrm{op}}(\Lambda^\sharp)\circ i_A$ sends $f\in L$ to
  $\widetilde{\Lambda^1}(\pi_{I^1}(f))=\Lambda^1(f)=\Lambda(f)$.
  Hence
  \[
    \Lambda=\iota_{B,\mathrm{op}}(\Lambda^\sharp)\circ i_A.
  \]

  Conversely, if a morphism $H:\mathcal{D}A\to B$ satisfies
  $\Lambda=\iota_{B,\mathrm{op}}(H)\circ i_A$, then the operator morphism into
  the quotient $L^1\to\widetilde{M}$ representing $H$ agrees with $\Lambda$ on
  $L$.
  By the uniqueness in Theorem~\ref{thm:operator-extension} we get
  $H=\Lambda^\sharp$.

  Since $A$ was arbitrary,
  Theorem~\ref{thm:left-adjunction-existence-from-universality} shows that
  $\iota_{B,\mathrm{op}}$ has a left adjoint.
\end{proof}

\begin{remark}
  What Section~\ref{chap:operator} imposes on the codomain is only that
  $\widetilde{M}$ be a Banach space, not Daniell closedness itself.
  For Daniell closed integration systems,
  Lemma~\ref{lem:quotient-is-banach} bridges the two.
  That this choice of subcategory is the natural one is seen in the next
  subsubsection (Proposition~\ref{prop:BInt-characterization}).

  Note also that the definition of the monotone convergence theorem
  (Definition~\ref{def:satisfies-mct}) is used with pointwise convergence
  unchanged.
  That is a property of \emph{objects} only, and it works solely to yield
  completeness of $\widetilde{M}$ by way of
  Lemma~\ref{lem:quotient-is-banach}.
  As long as no pointwise condition is imposed on the side of morphisms
  (Remark~\ref{rem:order-continuity-redundant}), there is no need to rewrite the
  monotone convergence theorem in an almost-everywhere sense.
\end{remark}

\subsubsection{What this category sees}\label{sec:what-op-ae-sees}

There is a point about the formulation of the main theorem of the previous
subsubsection that deserves attention.
The morphisms of $\mathbf{Int}_{\mathrm{op},\mathrm{ae}}$ are operator morphisms
into the quotient $L\to\widetilde{M}$, and by contractivity
(Lemma~\ref{lem:operator-is-contraction}) they map $\mathcal{N}_I$ to $0$.
Hence morphisms see only $\widetilde{L}$ and $\widetilde{M}$.
In this subsubsection we formulate this and state what follows from it.

\begin{definition}\label{def:preAL}
  We say that a pair $(E,\phi)$ is a \textbf{pre-AL space} if $E$ is a Riesz
  space and $\phi:E\to\mathbb{R}$ is a positive linear functional such that
  \[
    \|u\|_\phi:=\phi(|u|)
  \]
  is a norm on $E$.
  A morphism from $(E,\phi)$ to $(F,\psi)$ is a linear map $T:E\to F$ that is
  positive and satisfies $\psi\circ T=\phi$.
  We write this category as $\mathbf{AL}_0$.
\end{definition}

\noindent
Since $\|u+v\|_\phi=\phi(u+v)=\|u\|_\phi+\|v\|_\phi$ for $u,v\ge0$, the norm is
additive on the positive cone.
If it is complete this is an AL-space, the object of Kakutani's representation
theorem touched on in Appendix~\ref{sec:forget}.

\begin{lemma}\label{lem:quotient-is-preAL}
  If $(X,L,I)$ is an integration system then
  $(\widetilde{L},\widetilde{I})$ is a pre-AL space.
\end{lemma}
\begin{proof}
  By Lemma~\ref{lem:quotient-riesz}.
\end{proof}

\begin{proposition}\label{prop:Q-full-faithful}
  The assignment sending an object $A=(X,L,I)$ to
  $(\widetilde{L},\widetilde{I})$ and a morphism $\Lambda:L\to\widetilde{M}$ to
  the induced map $\widetilde{\Lambda}:\widetilde{L}\to\widetilde{M}$ defines a
  \textbf{faithful and full} functor
  \[
    Q:\mathbf{Int}_{\mathrm{op},\mathrm{ae}}\longrightarrow\mathbf{AL}_0 .
  \]
\end{proposition}
\begin{proof}
  \textbf{That $\widetilde{\Lambda}$ is a morphism of $\mathbf{AL}_0$} is
  Lemma~\ref{lem:operator-descends}.

  \textbf{Functoriality.}
  The identity morphism is $\pi_I$ and
  $\widetilde{\pi_I}=\operatorname{id}_{\widetilde{L}}$.
  As for composition, Definition~\ref{def:int-op-ae} gives
  $\Xi\circ\Lambda=\widetilde{\Xi}\circ\Lambda
  =\widetilde{\Xi}\circ\widetilde{\Lambda}\circ\pi_I$, so
  $\widetilde{\Xi\circ\Lambda}=\widetilde{\Xi}\circ\widetilde{\Lambda}$.

  \textbf{Faithfulness.}
  If $\widetilde{\Lambda}=\widetilde{\Lambda'}$ then
  $\Lambda=\widetilde{\Lambda}\circ\pi_I=\widetilde{\Lambda'}\circ\pi_I=\Lambda'$.

  \textbf{Fullness.}
  Given a morphism
  $T:(\widetilde{L},\widetilde{I})\to(\widetilde{M},\widetilde{J})$ of
  $\mathbf{AL}_0$, the map $\Lambda:=T\circ\pi_I:L\to\widetilde{M}$ is linear
  and positive and satisfies
  $\widetilde{J}(\Lambda f)=\widetilde{I}(\pi_I(f))=I(f)$, so it is an operator
  morphism into the quotient.
  By construction $\widetilde{\Lambda}=T$.
\end{proof}

\begin{corollary}\label{cor:op-ae-iso-criterion}
  $A\cong B$ in $\mathbf{Int}_{\mathrm{op},\mathrm{ae}}$ is equivalent to
  $(\widetilde{L},\widetilde{I})\cong(\widetilde{M},\widetilde{J})$ in
  $\mathbf{AL}_0$.
  Moreover the isomorphism is then isometric.
\end{corollary}
\begin{proof}
  A faithful and full functor reflects isomorphisms, and also preserves them.
  As for isometry, since both $T$ and its inverse $T^{-1}$ are contractive (use
  Lemma~\ref{lem:operator-is-contraction} with the same computation for
  morphisms of $\mathbf{AL}_0$),
  \[
    \|u\|=\|T^{-1}Tu\|\le\|Tu\|\le\|u\|
  \]
  follows.
\end{proof}

\noindent
That is, \textbf{$\mathbf{Int}_{\mathrm{op},\mathrm{ae}}$ sees neither the base
space $X$ nor the realization of $\widetilde{L}$ as a function space.}
If $\widetilde{L}$ is the same as a pre-AL space, the objects are not
distinguished.

\begin{proposition}\label{prop:D-op-is-banach-completion}
  Let $A=(X,L,I)$ be an integration system. Then
  $Q(i_A)=j_I:\widetilde{L}\to\widetilde{L^1}$ is a dense isometric embedding
  and $\widetilde{L^1}$ is a Banach space.
  That is, $Q(\mathcal{D}A)$ is the Banach completion of $Q(A)$.
\end{proposition}
\begin{proof}
  That $j_I$ is an injective Riesz homomorphism and is isometric is
  Lemma~\ref{lem:quotient-embedding}.
  Density of the image: for $f\in L^1$, taking $\{s_N\}\subset L$ as in
  Lemma~\ref{lem:l1-approx-by-source},
  \[
    \|\pi_{I^1}(f)-j_I(\pi_I(s_N))\|
    =\|\pi_{I^1}(f-s_N)\|
    =I^1(|f-s_N|)
    \longrightarrow0.
  \]
  Since $\mathcal{D}A$ is a Daniell closed integration system
  (Lemma~\ref{lem:daniell-ext-is-complete}),
  Lemma~\ref{lem:quotient-is-banach} shows that $\widetilde{L^1}$ is a Banach
  space.
\end{proof}

\begin{proposition}[characterization of $\mathbf{BInt}_{\mathrm{op},\mathrm{ae}}$]\label{prop:BInt-characterization}
  Let $A=(X,L,I)$ be an object of
  $\mathbf{Int}_{\mathrm{op},\mathrm{ae}}$.
  Then $A$ belonging to $\mathbf{BInt}_{\mathrm{op},\mathrm{ae}}$ is equivalent
  to $i_A:A\to\mathcal{D}A$ being an isomorphism of
  $\mathbf{Int}_{\mathrm{op},\mathrm{ae}}$.
  Moreover $\mathbf{BInt}_{\mathrm{op},\mathrm{ae}}$ is isomorphism-closed, and
  is the isomorphism closure of the class of Daniell closed objects.
\end{proposition}
\begin{proof}
  $\Leftarrow$: if $i_A$ is an isomorphism then
  Corollary~\ref{cor:op-ae-iso-criterion} shows that $\widetilde{L}$ is
  isometrically isomorphic to $\widetilde{L^1}$, and the latter is a Banach
  space by Proposition~\ref{prop:D-op-is-banach-completion}.

  $\Rightarrow$: if $\widetilde{L}$ is a Banach space then the image of the
  isometric embedding $j_I$ of
  Proposition~\ref{prop:D-op-is-banach-completion} is complete, hence closed,
  and it is also dense, so $j_I$ is surjective.
  As $j_I$ is a Riesz isomorphism, its inverse is positive and preserves the
  integral, so $Q(i_A)$ is an isomorphism of $\mathbf{AL}_0$.
  By Proposition~\ref{prop:Q-full-faithful} the morphism $i_A$ is an
  isomorphism.

  Being isomorphism-closed follows from
  Corollary~\ref{cor:op-ae-iso-criterion}.
  Isomorphisms are isometric, so if $\widetilde{L}$ is a Banach space then any
  pre-AL space isomorphic to it is a Banach space too.

  Finally we check that it is the closure.
  Daniell closed objects belong to
  $\mathbf{BInt}_{\mathrm{op},\mathrm{ae}}$ by
  Lemma~\ref{lem:quotient-is-banach}, and as just seen that subcategory is
  isomorphism-closed, so the isomorphism closure of the Daniell closed objects
  is contained in $\mathbf{BInt}_{\mathrm{op},\mathrm{ae}}$.
  Conversely, if $A\in\mathbf{BInt}_{\mathrm{op},\mathrm{ae}}$ then
  $A\cong\mathcal{D}A$, and $\mathcal{D}A$ is Daniell closed.
\end{proof}

\begin{remark}[the class of Daniell closed objects is not isomorphism-closed]\label{rem:op-ae-not-replete}
  The proposition above means that the subcategory to be taken when stating the
  main theorem in this category is
  $\mathbf{BInt}_{\mathrm{op},\mathrm{ae}}$, not the class of all Daniell closed
  objects.
  The latter is not isomorphism-closed.

  Indeed, take the example of
  Remark~\ref{rem:D-I-not-left-adjoint}: $X=\{0,1,2\}$,
  $L=\{(a,b,b)\}$, $I(a,b,b)=a$.
  Since $I(|(a,b,b)|)=|a|$ we have
  $\mathcal{N}_I=\{(0,b,b)\mid b\in\mathbb{R}\}$ and
  \[
    \widetilde{L}\cong\mathbb{R},
    \qquad
    \widetilde{I}(t)=t
  \]
  (the isomorphism being $\pi_I(a,b,b)\mapsto a$).
  As $\mathbb{R}$ is a Banach space, the proposition above gives
  $A\cong\mathcal{D}A$.
  On the other hand, as seen in that example, $L^1=\mathbb{R}^X\ne L$, so $A$ is
  not Daniell closed.

  This is not a counterexample to the main theorem
  (Theorem~\ref{thm:operator-main-theorem}), since being a reflective
  subcategory does not require being isomorphism-closed.
  The point is that $\mathbf{BInt}_{\mathrm{op},\mathrm{ae}}$ is the more
  natural choice for this category.
\end{remark}

\begin{remark}[the position of this section]\label{rem:op-ae-is-banach-completion}
  To summarize, what the Daniell extension gives in
  $\mathbf{Int}_{\mathrm{op},\mathrm{ae}}$ is the \textbf{Banach completion of
  the pre-AL space $(\widetilde{L},\widetilde{I})$}
  (Proposition~\ref{prop:D-op-is-banach-completion}).
  This category sees neither the base space nor the realization as a function
  space, so \textbf{one cannot expect the main theorem of
  Section~\ref{chap:operator} to carry the content ``accompanied by a concrete
  realization''.}
  The distinction from an abstract completion described in
  Subsection~\ref{sec:prior-work} is lost in the category of
  Section~\ref{chap:operator}.

  That distinction is retained on the side of the category of
  Section~\ref{chap:functoriality}.
  There, morphisms are the maps of base spaces themselves, and objects too are
  distinguished together with their realization $L\subset\mathbb{R}^X$.
  For the relation between the two sections see
  Subsection~\ref{sec:two-chapters}.

  Section~\ref{chap:operator} does not thereby become devoid of content.
  That a Banach completion \emph{exists} is general theory, but what
  Section~\ref{chap:operator} states is that \textbf{that completion is given as
  $\widetilde{L^1}$, the quotient of the integrable functions on the same base
  space} (Proposition~\ref{prop:D-op-is-banach-completion}).
  This is an assertion at the level of objects and is invisible to the morphisms
  of $\mathbf{Int}_{\mathrm{op},\mathrm{ae}}$.
  In other words, the concrete realization survives \emph{within the statement
  of the theorem} but not \emph{as categorical structure}.

  There is one more reason to keep Section~\ref{chap:operator}.
  Operator morphisms include, as morphisms and without modification, operations
  occurring in probability theory.
  Conditional expectation $E[\cdot\mid\mathcal{G}]$ is linear and positive and
  preserves the integral, so it is exactly an operator morphism of
  Section~\ref{chap:operator} (and it satisfies order continuity in the sense of
  Remark~\ref{rem:yamagami-operator} as well).
  In the framework of Section~\ref{chap:functoriality}, whose morphisms are maps
  of base spaces, this operation cannot be expressed as a morphism.
  The entrance to treating the It\^o integral and martingales in the same
  framework (Appendix~\ref{sec:forget}) is on the side of
  Section~\ref{chap:operator}.
\end{remark}

\subsection{Relation to the case of maps of base spaces as morphisms}

Section~\ref{chap:functoriality} and Section~\ref{chap:operator} put different
morphisms on the same objects and obtained a reflection theorem of the same
shape for each.
In this subsection we clarify the relation between the two.

We state the conclusion first.
The functor $\mathcal{E}$ commutes with the two reflections
(Proposition~\ref{prop:E-commutes-with-D}), and as a consequence
\textbf{the uniqueness part of the main theorem of
Section~\ref{chap:functoriality} follows from the theorem of
Section~\ref{chap:operator}}
(Corollary~\ref{cor:ch3-uniqueness-from-ch4}).
The existence part, however, \textbf{does not follow}.
That is, one cannot say from the operator version alone that the extension comes
from a map of base spaces
(Remark~\ref{rem:two-chapters-relation}).

\subsubsection{A functor between the ae-categories}

\begin{lemma}\label{lem:E-descends-to-ae}
  For a morphism of integration systems
  $F:(X,L,I)\to(Y,M,J)$ put
  \[
    \mathcal{E}_{\mathrm{ae}}(F):=\pi_J\circ F^*:L\to\widetilde{M} .
  \]
  Then the following hold.
  \begin{enumerate}
    \item $\mathcal{E}_{\mathrm{ae}}(F)$ is an operator morphism into the
      quotient.
    \item $F=_{\mathrm{ae}}G$ and
      $\mathcal{E}_{\mathrm{ae}}(F)=\mathcal{E}_{\mathrm{ae}}(G)$ are
      equivalent.
    \item The assignment $[F]\mapsto\mathcal{E}_{\mathrm{ae}}(F)$ defines a
      \textbf{faithful} functor
      \[
        \mathcal{E}_{\mathrm{ae}}:
        \mathbf{Int}_{\mathrm{ae}}\to\mathbf{Int}_{\mathrm{op},\mathrm{ae}} ,
      \]
      which is the identity on objects.
  \end{enumerate}
\end{lemma}
\begin{proof}
  (1) By Lemma~\ref{lem:base-map-gives-operator} the map $F^*$ is linear and
  positive and preserves the integral, and $\pi_J$ is a linear Riesz
  homomorphism (Lemma~\ref{lem:quotient-riesz}), so the composite is linear and
  positive as well.
  Moreover $\widetilde{J}(\pi_J(F^*f))=J(F^*f)=I(f)$.

  (2) $\pi_J(F^*f)=\pi_J(G^*f)$ is equivalent to
  $F^*f-G^*f\in\mathcal{N}_J$, that is, to
  $J(|F^*f-G^*f|)=0$.
  That this holds for every $f\in L$ was the definition of
  $F=_{\mathrm{ae}}G$ (Definition~\ref{def:ae-equal}).

  (3) The ``only if'' direction of (2) shows that the assignment is independent
  of the choice of representative.
  The converse direction gives injectivity, that is, faithfulness.
  For identity morphisms,
  $\mathcal{E}_{\mathrm{ae}}(\operatorname{id})=\pi_I$, which is the identity
  morphism of $\mathbf{Int}_{\mathrm{op},\mathrm{ae}}$
  (Definition~\ref{def:int-op-ae}).
  As for composition, let $F:(X,L,I)\to(Y,M,J)$ and
  $G:(Y,M,J)\to(Z,N,K)$. By the definition of composition in
  Definition~\ref{def:int-op-ae}, the map
  $\widetilde{M}\to\widetilde{N}$ induced by
  $\mathcal{E}_{\mathrm{ae}}(G)$ is $\pi_J(h)\mapsto\pi_K(G^*h)$. Hence, for
  every $f\in L$,
  \[
    \bigl(\widetilde{\mathcal{E}_{\mathrm{ae}}(G)}\circ\mathcal{E}_{\mathrm{ae}}(F)\bigr)(f)
    =\pi_K\bigl(G^*(F^*f)\bigr)
    =\pi_K\bigl((G\circ F)^*f\bigr)
    =\mathcal{E}_{\mathrm{ae}}(G\circ F)(f) .
  \]
\end{proof}

\begin{remark}\label{rem:E-ae-faithful-not-full}
  As stated at the beginning of Section~\ref{chap:operator}, the functor
  $\mathcal{E}:\mathbf{Int}\to\mathbf{Int}_{\mathrm{op}}$ before passing to the
  ae-quotient need not be faithful, because distinct maps of base spaces can
  give the same pullback.
  Once the ae-quotient is taken, exactly that difference is collapsed by
  $=_{\mathrm{ae}}$, so $\mathcal{E}_{\mathrm{ae}}$ becomes faithful.

  On the other hand $\mathcal{E}_{\mathrm{ae}}$ is \textbf{not full}.
  Only the morphisms coming from lattice homomorphisms, that is, the
  deterministic ones, lie in the image, and a general operator morphism is not
  among them.
  In the case of finite sets, the image consists only of those
  column-stochastic matrices each of whose rows has exactly one $1$.
\end{remark}

\subsubsection{Commutativity of the two reflections}

We write $\mathcal{D}_{\mathrm{op}}$ for the left adjoint of
$\iota_{B,\mathrm{op}}$ given by
Theorem~\ref{thm:operator-main-theorem}.
Its action is as stated in Theorem~\ref{thm:operator-main-theorem}.

\begin{proposition}\label{prop:E-commutes-with-D}
  The following diagram commutes.
  \[
    \begin{array}{ccc}
      \mathbf{Int}_{\mathrm{ae}}
        & \xrightarrow{\ \mathcal{D}_{\mathrm{ae}}\ }
        & \mathbf{DInt}_{\mathrm{ae}}\\[2pt]
      \big\downarrow{\scriptstyle\mathcal{E}_{\mathrm{ae}}}
        & & \big\downarrow{\scriptstyle\mathcal{E}_{\mathrm{ae}}}\\[2pt]
      \mathbf{Int}_{\mathrm{op},\mathrm{ae}}
        & \xrightarrow{\ \mathcal{D}_{\mathrm{op}}\ }
        & \mathbf{BInt}_{\mathrm{op},\mathrm{ae}}
    \end{array}
  \]
  That is,
  $\mathcal{E}_{\mathrm{ae}}\circ\mathcal{D}_{\mathrm{ae}}
  =\mathcal{D}_{\mathrm{op}}\circ\mathcal{E}_{\mathrm{ae}}$.
\end{proposition}
\begin{proof}
  On objects both are the Daniell extension
  $(X,L,I)\mapsto(X,L^1,I^1)$, and $\mathcal{E}_{\mathrm{ae}}$ does not change
  objects, so they agree.

  We argue for morphisms. Let $[F]:(X,L,I)\to(X',L',I')$ be a morphism of
  $\mathbf{Int}_{\mathrm{ae}}$.
  Write the injective Riesz homomorphism determined by
  Lemma~\ref{lem:quotient-embedding} as
  \[
    j':\widetilde{L'}\longrightarrow\widetilde{L'^1} .
  \]
  The unit morphism
  $i_{A'}:(X',L',I')\to(X',L'^1,I'^1)$ in
  $\mathbf{Int}_{\mathrm{op},\mathrm{ae}}$ is
  $\pi_{I'^1}|_{L'}=j'\circ\pi_{I'}$.

  Since $\mathcal{D}_{\mathrm{op}}$ is a left adjoint, for an operator morphism
  into the quotient $\Lambda:L\to\widetilde{L'}$ the morphism
  $\mathcal{D}_{\mathrm{op}}(\Lambda)$ is the unique operator morphism into the
  quotient $L^1\to\widetilde{L'^1}$ satisfying
  \[
    \mathcal{D}_{\mathrm{op}}(\Lambda)\bigr|_{L}=j'\circ\Lambda .
  \]
  Applying this to
  $\Lambda=\mathcal{E}_{\mathrm{ae}}[F]=\pi_{I'}\circ F^*$ gives the composite
  going down then across in the diagram, which we write as $\Phi_2$.

  The composite going across then down, on the other hand, is
  $\mathcal{E}_{\mathrm{ae}}$ applied to
  $\mathcal{D}_{\mathrm{ae}}[F]=[F^1]$, that is,
  \[
    \Phi_1:=\pi_{I'^1}\circ(F^1)^*:L^1\to\widetilde{L'^1} .
  \]
  Applying Lemma~\ref{lem:E-descends-to-ae} to the morphism $F^1$, the map
  $\Phi_1$ is an operator morphism into the quotient.

  We check that $\Phi_1$ satisfies the same characterization as $\Phi_2$.
  For $f\in L$ we have $(F^1)^*f=F^*f\in L'$, so the equality
  $\pi_{I'^1}|_{L'}=j'\circ\pi_{I'}$ of
  Lemma~\ref{lem:quotient-embedding} gives
  \[
    \Phi_1(f)
    =\pi_{I'^1}(F^*f)
    =j'\bigl(\pi_{I'}(F^*f)\bigr)
    =\bigl(j'\circ\mathcal{E}_{\mathrm{ae}}[F]\bigr)(f) .
  \]
  That is, $\Phi_1|_L=j'\circ\mathcal{E}_{\mathrm{ae}}[F]=\Phi_2|_L$.

  The system $(X',L'^1,I'^1)$ is a Daniell closed integration system
  (Lemma~\ref{lem:daniell-ext-is-complete}), and
  Lemma~\ref{lem:quotient-is-banach} shows that $\widetilde{L'^1}$ is a Banach
  space, so the \textbf{uniqueness} in
  Theorem~\ref{thm:operator-extension} applies and we obtain
  $\Phi_1=\Phi_2$.
\end{proof}

\begin{remark}
  Note that this proof relies on the uniqueness of
  Section~\ref{chap:operator}.
  Within the framework of Section~\ref{chap:functoriality} alone there is no
  means of saying that the two extensions agree, because they are morphisms of
  different categories.
\end{remark}

\subsubsection{Relation between the two main theorems}\label{sec:two-chapters}

\begin{corollary}\label{cor:ch3-uniqueness-from-ch4}
  Let $A$ be an object of $\mathbf{Int}_{\mathrm{ae}}$, let $B$ be an object of
  $\mathbf{DInt}_{\mathrm{ae}}$, and let
  $[H],[H']:\mathcal{D}_{\mathrm{ae}}A\to B$ be morphisms.
  Then
  \[
    \iota_{D,\mathrm{ae}}([H])\circ[i_A]
    =\iota_{D,\mathrm{ae}}([H'])\circ[i_A]
    \quad\Longrightarrow\quad
    [H]=[H']
  \]
  holds. That is, the uniqueness part of
  Theorem~\ref{thm:main-ae} follows from
  Theorem~\ref{thm:operator-extension} of
  Section~\ref{chap:operator}.
\end{corollary}
\begin{proof}
  Apply $\mathcal{E}_{\mathrm{ae}}$ to both sides.
  By Lemma~\ref{lem:E-descends-to-ae} the assignment
  $\mathcal{E}_{\mathrm{ae}}$ is a functor, and
  $\mathcal{E}_{\mathrm{ae}}[i_A]$ is the morphism determined by the inclusion
  $L\subset L^1$, so $\mathcal{E}_{\mathrm{ae}}[H]$ and
  $\mathcal{E}_{\mathrm{ae}}[H']$ are operator morphisms into the quotient
  $L^1\to\widetilde{M}$ agreeing on $L$.
  Since $B$ is Daniell closed, $\widetilde{M}$ is a Banach space by
  Lemma~\ref{lem:quotient-is-banach}, and the uniqueness in
  Theorem~\ref{thm:operator-extension} gives
  $\mathcal{E}_{\mathrm{ae}}[H]=\mathcal{E}_{\mathrm{ae}}[H']$.
  As $\mathcal{E}_{\mathrm{ae}}$ is faithful
  (Lemma~\ref{lem:E-descends-to-ae}), we obtain $[H]=[H']$.
\end{proof}

\begin{remark}[relation between the two sections]\label{rem:two-chapters-relation}
  By contrast, the \textbf{existence} part of
  Theorem~\ref{thm:main-ae} does not follow from the theorem of
  Section~\ref{chap:operator}.
  Given $[F]:A\to\iota_{D,\mathrm{ae}}B$,
  Theorem~\ref{thm:operator-extension} supplies an extension of
  $\mathcal{E}_{\mathrm{ae}}[F]$ as a morphism of
  $\mathbf{Int}_{\mathrm{op},\mathrm{ae}}$, but it does not guarantee that this
  lies in the image of $\mathcal{E}_{\mathrm{ae}}$, that is, that it comes from
  a map of base spaces, because $\mathcal{E}_{\mathrm{ae}}$ is not full
  (Remark~\ref{rem:E-ae-faithful-not-full}).
  To obtain existence one needs the construction specific to
  Section~\ref{chap:functoriality}
  (Lemma~\ref{lem:daniell-ext-is-functor-1}), by which the extension is
  obtained without changing the map of base spaces.

  The contributions of the two theorems divide as follows.

  \begin{center}
  \begin{tabular}{@{}lll@{}}
    \toprule
     & Existence & Uniqueness \\
    \midrule
    Section~\ref{chap:functoriality}, non-ae & needs its own construction & follows automatically (no content) \\
    Section~\ref{chap:functoriality}, ae & as above & by density \\
    Section~\ref{chap:operator} (operators) & not obtainable from the operator version & by contractivity and density \\
    \bottomrule
  \end{tabular}
  \end{center}

  Care is needed in how this is phrased.
  Since $\mathcal{E}_{\mathrm{ae}}$ is faithful,
  $\mathbf{Int}_{\mathrm{ae}}$ is embedded (without changing objects) in
  $\mathbf{Int}_{\mathrm{op},\mathrm{ae}}$.
  In that sense \textbf{the category of
  Section~\ref{chap:functoriality} is a wide subcategory of the category of
  Section~\ref{chap:operator}}, and it is not accurate to say that ``neither is
  a special case of the other''.
  Nor is Section~\ref{chap:functoriality} in the relation of completing a single
  theorem together with Section~\ref{chap:operator}: it proves both existence
  and uniqueness on its own.
  The relations that can be stated precisely are the following.

  \begin{itemize}
    \item $\mathcal{E}_{\mathrm{ae}}$ is faithful but not full
      (Remark~\ref{rem:E-ae-faithful-not-full}).
    \item The two reflections are compatible under this functor
      (Proposition~\ref{prop:E-commutes-with-D}).
    \item The \emph{uniqueness} of the base-space-map version can be derived
      from the operator version
      (Corollary~\ref{cor:ch3-uniqueness-from-ch4}).
    \item But \textbf{that the extension stays in the image of
      $\mathcal{E}_{\mathrm{ae}}$, that is, that it comes from a map of base
      spaces, cannot be derived from the operator version alone.}
    \item Furthermore, \textbf{$\mathbf{Int}_{\mathrm{op},\mathrm{ae}}$ sees
      neither the base space nor the realization as a function space}
      (Proposition~\ref{prop:Q-full-faithful},
      Corollary~\ref{cor:op-ae-iso-criterion}).
      Objects are distinguished only up to the pre-AL space
      $(\widetilde{L},\widetilde{I})$, and what the Daniell extension gives is
      its Banach completion
      (Proposition~\ref{prop:D-op-is-banach-completion}).
  \end{itemize}

  The last two points are the reasons for keeping
  Section~\ref{chap:functoriality} as an independent section.
  The first is a reason at the level of theorems and the second at the level of
  categories, and the latter is the stronger of the two.
  Only the side of Section~\ref{chap:functoriality} retains the distinction
  between an abstract completion and a concrete realization
  (Subsection~\ref{sec:prior-work}) as categorical structure.
  For details see Subsection~\ref{sec:what-op-ae-sees}.
\end{remark}

\section{Appendix}

\subsection{Conditions for the existence of a left adjoint}

The content of this subsection is standard general category theory and contains
nothing specific to this paper.
We refer to \cite{MacLane} for proofs and state only the assertions, in the form
used in the body.

\begin{definition}
  Let $F$ be a functor from a category $C$ to a category $D$.
  Let $Y$ be an object of $D$, let $G_Y$ be an object of $C$, and let
  \[
    \eta_Y:Y\to F(G_Y)
  \]
  be a morphism of $D$.
  We say that $\eta_Y$ is a \textbf{universal morphism from $Y$ to $F$} if the
  following holds.

  For every object $X$ of $C$ and every morphism $u:Y\to F(X)$ of $D$, there is
  a unique morphism
  \[
    u^\sharp:G_Y\to X
  \]
  of $C$ satisfying
  \[
    u=F(u^\sharp)\circ\eta_Y .
  \]

  We call the above condition the \textbf{universal property} of
  $\eta_Y:Y\to F(G_Y)$.
\end{definition}

\begin{theorem}\label{thm:left-adjunction-existence-from-universality}
  Let $F$ be a functor from a category $C$ to a category $D$.
  If for every object $Y$ of $D$ one can choose a universal morphism
  \[
    \eta_Y:Y\to F(G_Y)
  \]
  from $Y$ to $F$, then a left adjoint
  \[
    G:D\to C
  \]
  of $F$ exists.
  Moreover the action of $G$ on morphisms is uniquely determined by the
  universal property of each $\eta_Y$.
\end{theorem}

\noindent
For a proof see \cite{MacLane}, Chapter IV, \S1, Theorem 2.
The argument runs as follows: assigning $G_Y$ to an object $Y$, and to a
morphism $f:Y\to Y'$ the morphism determined by the universal property of
$\eta_{Y'}\circ f$, makes $G$ a functor, and
$\Psi_{Y,X}(u):=F(u)\circ\eta_Y$ gives a natural bijection
$\operatorname{Hom}_C(G(Y),X)\cong\operatorname{Hom}_D(Y,F(X))$.

\subsection{What is the base space doing?}\label{sec:forget}

A Daniell-style integral is nothing more than a positive linear functional,
satisfying a certain continuity, on a Riesz space of functions on a base space.
From this appearance one might think that the base space is not essential and
could be removed.
In this subsection we see that the role the base space plays is exactly one
thing: \textbf{to supply a place where infima are computed outside $L$}.
We then explain that the direction of ``forgetting the base space'' does not
work as it stands.

\subsubsection{The base space appears in only one place}

In the definition of an integration system $(X,L,I)$, the only place where the
base space $X$ is essentially at work is the hypothesis of the continuity axiom,
that is, the left-hand side of
\[
  f_n\downarrow0\ \text{(pointwise)}\implies I(f_n)\to0 .
\]
If this is replaced by decrease in the order of $L$ itself, then even the
ordinary integral ceases to be an integral.

\begin{example}\label{ex:fat-cantor}
  Let $L$ be the set of all step functions on the closed interval $[0,1]$, that
  is, the set of all functions for which there is a finite partition
  $0=x_0<x_1<\cdots<x_m=1$ such that the function is constant on each
  $[x_{i-1},x_i)$ with $i<m$ and on $[x_{m-1},1]$.
  The partition points are strictly increasing, and singletons are not allowed.
  Taking $I$ to be the ordinary Riemann integral, $(L,I)$ is an integration
  system on $[0,1]$.

  Let $A\subset[0,1]$ be a fat Cantor set of measure $1/2$ and write
  $[0,1]\setminus A=\bigsqcup_kI_k$ with $I_k=(a_k,b_k)$.
  Putting
  \[
    f_n:=1_{[0,1]}-\sum_{k=1}^n1_{[a_k,b_k)}\in L ,
  \]
  the sequence $\{f_n\}$ is monotonically decreasing and $I(f_n)\to1/2$.

  On the other hand, if $g\in L$ satisfies $g\le f_n$ for every $n$ then
  $g\le1_{A\setminus\{a_k\mid k\ge1\}}$ pointwise.
  The set $[0,1]\setminus A$ is a dense open set, and each interval of
  constancy of $g$ is nondegenerate, so it meets some $I_k$.
  At a point of that intersection the upper bound is $0$, so the value of $g$ on
  each interval of constancy is at most $0$, that is, $g\le0$.
  Conversely, since $f_n\ge0$, if $g\le0$ then $g\le f_n$.
  Hence the set of lower bounds of $\{f_n\}$ is $\{g\in L\mid g\le0\}$, whose
  greatest element is $0$.

  That is, in the order of $L$ we have $f_n\downarrow0$, yet
  $I(f_n)\downarrow1/2$.
  In terms of pointwise convergence, $f_n\to1_A\ne0$, and $f_n$ does not
  decrease to $0$.
\end{example}

\noindent
What this example shows is \textbf{the difference between taking the infimum
inside $L$ and taking it outside}.
Taken inside $L$ we get $\inf_nf_n=0$, but taken inside
$\mathbb{R}^{[0,1]}$ we get $\inf_nf_n=1_A$.
What Daniell's continuity axiom requires is the latter, and the base space $X$
supplies $\mathbb{R}^X$, a container large enough for that purpose.

This difference is also at work in the definition of operator morphisms in
Section~\ref{chap:operator}
(Remark~\ref{rem:order-continuity-redundant}).
There too, it is exactly this that makes the statement fail if the order
continuity condition is read in the order of $M$.

\subsubsection{The direction of forgetting the base space does not work as it stands}

One might also consider removing the base space.
In what follows, an \textbf{abstract Riesz space} means an ordered vector space
that is a lattice (Remark~\ref{rem:riesz-terminology}).
The idea is to recover the base space from an abstract Riesz space $L$ by
\[
  \operatorname{Spec}(L):=\operatorname{Hom}_{\mathbf{Riesz}}(L,\mathbb{R}) .
\]
Assigning to each $f\in L$ the function
$\widehat{f}:\operatorname{Spec}(L)\to\mathbb{R};x\mapsto x(f)$ yields a left
adjoint of the functor forgetting the base space.
The adjunction itself does hold. What fails is rather that the unit
$L\to\widehat{L}$ be an isomorphism, that is, that the recovery can actually be
carried out.

There is also the problem that $\operatorname{Spec}(L)$ is too large, but the
fatal point is the following: \textbf{$\operatorname{Spec}$ collapses on exactly
the objects that matter.}

\begin{proposition}\label{prop:spec-collapses}
  Let $L$ be the abstract Riesz space obtained from the set of all Lebesgue
  integrable functions on $[0,1]$ by dividing by the null functions (that is,
  $L^1[0,1]$). Then
  \[
    \operatorname{Hom}_{\mathbf{Riesz}}(L,\mathbb{R})=\{0\} .
  \]
\end{proposition}
\begin{proof}
  Let $\varphi$ be a Riesz homomorphism.
  If $f\ge0$ then $f=f\vee0$ gives
  $\varphi(f)=\varphi(f)\vee0\ge0$, so $\varphi$ is positive.
  A positive linear functional on a Banach lattice is norm bounded
  (\cite{AbramovichAliprantis}), so from
  $L^1[0,1]^*=L^\infty[0,1]$ there is a $g\in L^\infty$ with $g\ge0$ such that
  $\varphi(f)=\int_0^1fg$.

  Suppose $\varphi\ne0$. Then one can take a set $S$ of positive measure on
  which $g>0$.
  Split $S$ into two sets $S_1,S_2$ of positive measure and put
  $f:=1_{S_1}$ and $h:=1_{S_2}$.
  Their supports are disjoint, so $f\wedge h=0$, and since $\varphi$ is a Riesz
  homomorphism,
  \[
    0=\varphi(f\wedge h)=\varphi(f)\wedge\varphi(h)=\min\{\varphi(f),\varphi(h)\} .
  \]
  But $\varphi(f)=\int_{S_1}g>0$ and $\varphi(h)=\int_{S_2}g>0$, a
  contradiction.
  Hence $\varphi=0$.
\end{proof}

\noindent
So $\operatorname{Spec}$ is not ``too large''; it
\textbf{vanishes the moment one divides by the null functions}.
And if one is to do measure theory seriously, dividing by the null functions
cannot be avoided.

Note that treating order continuous functionals on abstract Riesz spaces and
their extensions is itself a classical subject, and anyone pursuing the
direction of this subsection will find \cite{LuxemburgZaanen} and
\cite{FremlinTRS} to be the places to look.

\subsubsection{The correct form, and further prospects}

\begin{remark}[the correct form]
  The idea of ``forgetting the base space and recovering it when needed'' has a
  classical answer.
  For a Daniell closed integration system, the quotient $\widetilde{L}$ of
  Subsection~\ref{sec:ae} is a Banach lattice whose norm is additive on the
  positive cone, that is, an \textbf{abstract Lebesgue space} (AL-space).
  And by Kakutani's representation theorem \cite{Kakutani} (see also
  \cite{Yosida} for the representation of Riesz spaces as function spaces),
  every AL-space is isomorphic to some $L^1(\mu)$.
  That is, the recovery of the base space is already complete here.

  As stated in Subsection~\ref{sec:prior-work} of
  Section~\ref{chap:intro}, however, the more one moves in this direction the
  closer the main theorem of this paper comes to a restatement of an abstract
  completion.
  \textbf{Abstraction and concrete realization stand in a trade-off.}

  The category $\mathbf{Int}_{\mathrm{op},\mathrm{ae}}$ of
  Section~\ref{chap:operator} has in fact already arrived at the end of this
  direction.
  That category distinguishes objects only up to the pre-AL space
  $(\widetilde{L},\widetilde{I})$
  (Subsection~\ref{sec:what-op-ae-sees}), and what the Daniell extension gives
  is its Banach completion.
  That is, both ends of the trade-off sit side by side inside this paper.
\end{remark}

We mention two directions of generalization.

\begin{itemize}
  \item \textbf{Changing the object of values.}
    If $f_n$ is taken to be Banach-space valued while the dominating functions
    $\varphi_n$ remain real-valued, then the proofs of this paper go through
    almost unchanged and one obtains the Bochner integral.
    This is the very subject of \cite{Mikusinski} and is known, but it is a
    technical bridge towards stochastic integration.
  \item \textbf{Internalizing the definition itself.}
    This is the direction of taking the Riesz space and the reals inside a topos
    rather than in the category of sets.
    The ``base space'' is then replaced by the topos itself.
    There is earlier work in this direction.
    \cite{Borger2010} treats measures and integrals inside a monoidal closed
    category and defines the space of simple functions by a universal property.
    The same author's \cite{BorgerHabil} is broader in scope, but we have not
    obtained its text.
    That both are theorems different from this paper's was stated in
    Subsection~\ref{sec:prior-work}.
    There is also \cite{CoquandSpitters} as a constructive, point-free theory of
    integration.
    For presheaf toposes, and especially for presheaves indexed by a filtration,
    a relation to martingales is to be expected, since conditional expectation
    has the shape of an operator morphism of
    Section~\ref{chap:operator} (see also \cite{Simpson}).
\end{itemize}

\noindent
In either direction, the main theorem of this paper retains its meaning only in
so far as it is ``accompanied by a concrete realization''.
The question of how far to abstract, and at what point concreteness is given up,
is one that the framework of this paper raises and leaves open.

\end{document}